\documentclass[lettersize,onecolumn]{IEEEtran}
\ifCLASSINFOpdf
\else
\fi
\usepackage {xr}
\makeatletter

\usepackage{amsmath,amsfonts}
\usepackage{algorithmic}
\usepackage{algorithm}
\usepackage{array}
\usepackage[caption=false,font=normalsize,labelfont=sf,textfont=sf]{subfig}
\usepackage{textcomp}
\usepackage{stfloats}
\usepackage{url}
\usepackage{verbatim}
\usepackage{graphicx}
\RequirePackage{amsthm,amsmath,amsfonts,amssymb}
\RequirePackage[numbers,sort&compress]{natbib}
\RequirePackage[colorlinks,linkcolor=blue,citecolor=blue,urlcolor=blue]{hyperref}
\RequirePackage{graphicx}
\usepackage{comment}

\theoremstyle{plain}

\newtheorem{theorem}{Theorem}[section]
\newtheorem{lemma}[theorem]{Lemma}

\theoremstyle{remark}

\theoremstyle{plain}
\newtheorem{corollary}[theorem]{Corollary}
\newtheorem{proposition}{Proposition}

\usepackage{graphicx} 
\usepackage{amsmath, amsthm, amsfonts,amssymb}
\usepackage{bm}
\usepackage{indentfirst}
\usepackage{amsmath}
\usepackage{mathrsfs}
\usepackage{subfiles}
\usepackage{verbatim}
\usepackage{float}
\usepackage{caption}
\usepackage{subcaption}
\usepackage{algorithm}
\usepackage{algorithmic}
\usepackage{longtable}
\usepackage{multirow}
\usepackage{multicol}
\usepackage{xcolor}
\usepackage{latexsym}
\usepackage{amsthm}
\usepackage{amsthm,amsmath,amssymb}
\usepackage{mathrsfs}

\newcommand{\blue}

\begin{document}
%
\title{Optimal Estimation of Shared Singular Subspaces across Multiple Noisy Matrices}
%
%
%

\author{Zhengchi~Ma
        and~Rong~Ma
\thanks{Zhengchi Ma is with the Department of Electrical \& Computer Engineering at Duke University, Durham, NC 27708 USA (email: zhengchi.ma@duke.edu).}
\thanks{Rong Ma is with the Department of Biostatistics at Harvard T.H. Chan School of Public Health, Boston, MA 02115 USA, and the Eric and Wendy Schmidt Center at the Broad Institute of MIT and Harvard, Cambridge, MA 02142 USA (email: rongma@hsph.harvard.edu).}}

%
%

\markboth{IEEE Transactions on Information Theory,~Vol.~72, No.~5, May~2026}%
{Shell \MakeLowercase{\textit{Zhengchi Ma and Rong Ma}}: Optimal Estimation of Shared Singular Subspaces across Multiple Noisy Matrices}
%



\maketitle

\begin{abstract}
Estimating singular subspaces from noisy matrices is a fundamental problem with wide-ranging applications across various fields. Driven by the challenges of data integration and multi-view analysis, this study focuses on estimating shared singular subspaces across multiple matrices within a low-rank matrix denoising framework. A common approach for this task is to perform singular value decomposition on the stacked matrix (Stack-SVD), which concatenates all the matrices. We establish that Stack-SVD achieves minimax rate-optimality when the true singular subspaces of the noisy matrices are identical, whereas a popular alternative approach based on SVD of concatenated singular vector matrices (Average-SVD) can be sub-optimal. We then tackle the more complex scenario where the true singular subspaces are only partially shared across matrices. For various cases of partial sharing, we rigorously characterize the conditions under which Stack-SVD remains effective, achieves minimax optimality, or fails to deliver consistent estimates, offering theoretical insights into its practical applicability. To address the limitations of Stack-SVD in scenarios with partial sharing, we propose novel estimators and an efficient algorithm designed to identify both shared and unshared singular vectors. We further prove that these methods attain minimax rate-optimality under partial sharing. Extensive simulations and real-world data applications demonstrate the advantages of our proposed approach.
\end{abstract}


%
\IEEEpeerreviewmaketitle

\section{Introduction}
%
%
%
%

\IEEEPARstart{S}{pectral} methods, commonly referred to as a collection of algorithms built upon the singular value decomposition (SVD) of some properly designed matrices constructed from data, are widely utilized across various disciplines such as statistics, probability theory, and statistical machine learning \cite{chen2021spectral}. As an important theoretical framework that guarantees the performance of spectral methods, the perturbation theory of matrix singular subspaces has been extensively studied \citep{daviskahan1970,wedin1972,caizhang2018}.
For example, \cite{wedin1972} provided a uniform perturbation bound for both left and right singular subspaces, and \cite{caizhang2018} offered separate rate-optimal perturbation bounds for these subspaces under the same perturbation. Subsequent research has extended the estimation and perturbation analysis of singular subspaces to more complex scenarios, including settings with geometric or computational constraints \citep{ma2021optimal,structuredpca2021}, challenging distance measures \citep{cape2019two}, non-Gaussian or heteroskedastic noise \citep{zhang2022heteroskedastic}, missing data \citep{cai2021subspace}, and tensor data \citep{zhang2018tensor}.

Closely related to singular subspace estimation, low-rank matrix denoising is a fundamental task in applications such as image processing \citep{image1}, recommendation systems \citep{recommendation1}, and statistical genetics \citep{genetics1}. The objective is to recover the underlying signal matrix from its noisy observation by leveraging the intrinsic low-rank structure of data matrices. Given its critical role in many applications, SVD-based methods have been extensively studied from the theoretical perspective. For instance,  under the low-rank matrix denoising model, \cite{Capitaine2007} explore the high-dimensional asymptotic behavior of singular values and vectors; \cite{donoho2013threshold} focus on singular-value thresholding algorithms for recovering the noiseless matrix.

Recently, the problem of identifying shared and distinct patterns of variability across multiple data matrices has garnered significant attention across various research domains. For example, \cite{spectralmethod} employs rotational bootstrap methods and random matrix theory to decompose the observed spectrum into joint, individual, and noise subspaces. \cite{DCDLF} investigates the uncorrelatedness between distinct latent factors derived from different data views.  \cite{JIVE} introduces the Joint and Individual Variation Explained (JIVE) framework, to quantify the amount of joint variation between multiple data types. Extensions of JIVE, such as the angle-based JIVE (AJIVE) \citep{aJIVE}, the robust angle-based JIVE (RaJIVE) \citep{RaJIVE} and the Interpretive JIVE \cite{interpretiveJIVE}, further enhance its applicability. Supervised extensions of JIVE methods include Supervised Joint and Individual Variation Explained (SJIVE) \cite{sJIVE} and sparse exponential family sJIVE (sesJIVE) \cite{sesJIVE}. \cite{swain2025better} demonstrate that cross and joint covariance matrices detect shared signals between high-dimensional variables more efficiently than self-covariance matrices. \cite{baharav2025stacked} employ random matrix theory to demonstrate that optimally weighting datasets before stacking yields superior estimation of shared singular subspace compared to stacking individually estimated singular vectors.
In cases where the data exhibits partially shared patterns across matrices, a number of methodologies have been proposed. Notable approaches include the partially-joint structure identification \citep{partiallyjointidentify}, Structural Learning and Integrative DEcomposition of multi-view data (SLIDE) \citep{SLIDE}, hierarchical nuclear norm penalization \citep{Yi2022HierarchicalNN}, and Data Integration Via Analysis of Subspaces (DIVAS) \citep{DIVAS}. Furthermore, some asymptotic results for shared subspace estimation based on averages are obtained in \cite{fan2017distributed}. Other related topics include partial least square \cite{garthwaite1994interpretation} or canonical correlation analysis \cite{bao2019canonical}, which seeks associations between two matrices by maximizing their covariance or correlation, and high-dimensional hypothesis testing across multiple matrices \cite{ding2024two}.
Despite these  achievements, most existing research focuses on a specific estimation algorithm and its theoretical properties, while the fundamental limits of the estimation problem—such as the optimal rates of estimation under complete or partial information sharing—remain largely unexplored.

Here we examine the estimation of shared singular subspaces across multiple matrices, focusing on a class of methods,  referred to as the Stack-SVD method and exploring its properties in various scenarios. Stack-SVD estimates the shared singular subspace by concatenating noisy individual matrices and computing the top singular subspace of the resulting stacked matrix. 
For example, consider two rank-\(r\) noiseless matrices, \(X_1\in\mathbb{R}^{n\times p_1}\) and \(X_2\in\mathbb{R}^{n\times p_2}\), that share some left singular subspaces. Suppose we observe  $Y_1=X_1+Z_1$ and $Y_2=X_2+Z_2$, where $Z_1$ and $Z_2$ are additive noises. By stacking these noisy matrices to form \((Y_1~Y_2)\in\mathbb{R}^{n\times (p_1+p_2)}\), Stack-SVD estimates the shared left singular subspace of \(X_1\) and \(X_2\) using the top \(r\) left singular vectors of the stacked noisy matrix, sorted by descending singular values.

Stack-SVD is widely employed in applications, particularly in scenarios where multiple matrices are expected to share singular subspaces. For example, in multi-modal single-cell clustering, multiple matrices are often combined by stacking, followed by applying SVD or PCA for joint dimensionality reduction \citep{SVDclustering}. Similarly, in Multi-Omics Factor Analysis (MOFA+) \cite{MOFA+}, spectral factor analysis is performed on stacked data matrices from different modalities to identify shared cellular characteristics. Scanorama \citep{Scanorama}, an algorithm for integrating multiple heterogeneous single-cell datasets, also assumes shared singular subspaces by applying SVD to the stacked data matrix, to align datasets across conditions.
In electronic health record (EHR) data analysis, multiple clinical concepts are frequently stacked into a single matrix, facilitating the identification of patterns and relationships across clinical variables such as diagnoses and laboratory results. Techniques from natural language processing, such as the Positive Pointwise Mutual Information (PPMI) matrix, are commonly used in this context. Applying SVD to the PPMI matrix is conceptually analogous to performing SVD on a stacked matrix and inherently assumes a shared singular subspace across matrices \citep{EHRPPMI}.

We show that the convergence rate of Stack-SVD is minimax optimal when the noiseless matrices have identical singular subspaces. An alternative method, Average-SVD (also known as Distributed PCA or Principal Angle Analysis), involves computing leading singular vectors of the average matrix of the individually estimated subspaces (defined in Section \ref{suboptimal.sec}). Our analysis indicates the advantage of Stack-SVD over several alternative methods, including the Average-SVD estimators that are also widely used \cite{aJIVE,fan2017distributed}, especially when there is substantial signal imbalance; see Sections \ref{suboptimal.sec}, \ref{simulation} and \ref{realdata}. In cases where individual unshared singular vectors are present, we show that some modified Stack-SVD estimators achieve the optimal rates. Furthermore, we propose an algorithm to distinguish between shared and unshared singular vectors within the stacked noisy matrix, enabling practical implementation. Our findings deepen the understanding of subspace estimation under noise, particularly in settings involving both shared and unshared singular vectors across data matrices, and provide insights for applications in various fields. 

We consider the following two specific settings. For the case of identical singular subspace, we consider $k$ noiseless low-rank matrices $ X_i \in \mathbb{R}^{n \times p_i} $ $(i=1,2,...,k)$  sharing the same left singular subspace, whose SVD can be written as $X_1 = U \Sigma_1 V_1^T, ~ X_2 = U \Sigma_2 V_2^T, ~ ...~ X_k = U \Sigma_k V_k^T,$
where $U$ is the shared left singular subspace while $\Sigma_i$ and $V_i$ are the respective singular values and right singular subspaces for $X_i$. The observations are noisy versions of these matrices, given by  the matrix denoising model \citep{donoho2013threshold,structuredpca2021}:
\begin{equation}\label{model}
Y_1 = X_1 + Z_1, \quad Y_2 = X_2 + Z_2,\quad ...\quad Y_k = X_k + Z_k,
\end{equation}
where $Z_i\in\mathbb{R}^{n\times p_i} $ $(i=1,2,...,k)$ are additive noise matrices with zero-mean independent sub-Gaussian entries. 

For the case of partially shared singular subspace, where each noiseless matrix $X_i$ may contain  some unique singular vectors, we consider a two-matrix case for brevity, although the results can also be generalized to $k$ matrices. Specifically, we assume
$$
X_1 = (U_r~U_{1*}) \Sigma_1 V_1^T, \quad X_2 = (U_r~U_{2*}) \Sigma_2 V_2^T,$$
where $U_r\in \mathbb{O}(n,r)$ contains the shared $r$ left singular vectors, whereas $U_{1*}$ and $U_{2*}$ are individual unshared subspaces. $\Sigma_1$ and $\Sigma_2$ are diagonal matrices with non-zero positive entries. Note that the singular values in $\Sigma_1$ and $\Sigma_2$ are not necessarily in the decreasing order -- the shared and unshared vectors might be shuffled after sorting the singular values decreasingly. Again, we observe the noisy version of $X_1$ and $X_2$ as in (\ref{model}). The above model is closely related to  AJIVE  \cite{aJIVE,RaJIVE}: $$
X_i=J_i+A_i,~ i=1,...,k,
$$
where $J_i$ contains the joint structure and $A_i$ is the individual structure. Here $J_i$ and $A_j$ are required to be orthogonal for all $i,j=1,...,k$. Indeed, by rewriting $\Sigma_i=\begin{pmatrix}
    \Sigma_{11} & \\ & \Sigma_{12}
\end{pmatrix}$ and $V_i^T=\begin{pmatrix}
    V_{11}^T\\V_{12}^T
\end{pmatrix}$, the signal matrices can be expressed as
$$X_i = (U_r~U_{i*}) \begin{pmatrix}
    \Sigma_{11} & \\ & \Sigma_{12}
\end{pmatrix}\begin{pmatrix}
    V_{11}^T\\V_{12}^T
\end{pmatrix}=U_r\Sigma_{11}V_{11}^T+U_{i*}\Sigma_{12}V_{12}^T:=J_i+A_i,$$
which reduces to above AJIVE model. 
However, existing work \cite{aJIVE,RaJIVE} focuses on estimating the component matrices $\{J_i\}$ and $\{A_i\}$ from $\{Y_i\}$, whereas in this study, our goal is to estimate the shared  subspace $U_r$. We will focus on elucidating the fundamental limits of subspace estimation across diverse scenarios, analyzing the statistical properties of the  Stack-SVD procedure, and enhancing its performance in challenging regimes.

Our main contributions can be summarized as follows.
\begin{itemize}
    \item We establish the minimax optimal rates for the singular subspace estimation problem when the true subspaces of different matrices are completely shared. We show that Stack-SVD is minimax rate-optimal in this case, while Average-SVD is sub-optimal. We also demonstrate the minimax optimality of individual SVD estimators under imbalanced signals and dimensions. Our analysis reveals several phase transition phenomena in the problem as a function of the underlying signal-to-noise ratio (SNR), highlighting how the interplay among multiple matrices determines the fundamental limits of estimation.
    \item We establish the minimax optimal rates and identify the optimal estimators for the shared singular subspace when multiple matrices have partially shared singular subspaces. Focusing on two important scenarios of partial sharing where the unshared singular subspaces are mutually orthogonal, we provide theoretical insights on the performance of the Stack-SVD procedure, including its advantages in capturing individually non-identifiable singular vectors, and its limitations in excluding unshared singular vectors. We  extend our analyses to the case when the unshared  subspaces are not mutually orthogonal.
    \item To address the limitation of Stack-SVD in the presence of unknown, possibly unshared singular vectors, we propose an efficient algorithm for detecting and discriminating the shared and unshared singular vectors between two matrices. We also provide theoretical guarantees for the proposed method, justifying its practical advantages.
    \item Numerical results based on extensive simulations  are obtained to verify our theoretical findings. Applying our proposed methods to multiple single-cell omics datasets, we demonstrate the advantages of our methods in integrating and representing the shared latent cell-type structures captured by different sequencing technologies.
\end{itemize}

The rest of the paper is organized as follows. We introduce some mathematical notations. Section \ref{identical} concerns the  optimal estimation of the fully shared singular subspaces across multiple matrices, whereas Section \ref{unshared} focuses on the estimation in the presence of unshared subspaces across matrices. Section \ref{algorithm} concerns our proposed algorithms for tracing shared and unshared singular vectors. Sections \ref{simulation} and \ref{realdata} contain our numerical results from simulations and real data analysis. In Section \ref{discussion}, we discuss the broader potential of our work.

\subsection{Notation}
For $a,b\in\mathbb{R}$, let $a\wedge b=\min(a,b),a\vee b=$ $\max(a,b).$ For a matrix $A\in\mathbb{R}^{n\times p}$, write the SVD as $A=U\Sigma V^T$, where $\Sigma=\operatorname{diag}\{\sigma_{1}(A),\sigma_{2}(A),\cdots\}$ with the singular values $\sigma_1(A)\geq\sigma_{2}(A)\geq\cdots\geq0$ in descending order unless particularly specified. When there is no confusion, we also use $\mathrm{diag}(\cdots)$ to denote block diagonal matrices when the entries are matrices. We use $\sigma_{\min}(A),\sigma_{\max}(A)$ to denote the smallest and largest non-trivial singular values of $A$. We use $\sigma_{(i)}(A)$ to denote the singular value of $A$ corresponding to the $i$-th left singular vector of $A$ and $\sigma_i(A)$ denote the $i$-th largest singular value of $A$. For matrix norms , $\|A\|=\sigma_{\max}(A)$ is the spectral norm and  $\|A\|_F=\sqrt{\sum\sigma_i^2(A)}$ is the Frobenius norm. The $\sin\Theta$ distance between matrices A and B is denoted as $\|\sin\Theta(A,B)\|$. For two scalar sequences $\{a_n\}_{n\geq 1}$ and $\{b_n\}_{n\geq 1}$, we denote $a_n \gtrsim b_n(a_n \lesssim b_n)$ if there exist a universal constant $C$ such that $a_n\geq Cb_n(a_n\leq Cb_n)$ and $a_n \asymp b_n$ if both $a_n \gtrsim b_n$ and $a_n \lesssim b_n$. We use $\mathbb{O}(n,p)$ to denote the class of orthonormal matrices in $\mathbb{R}^{n\times p}$. Lastly, $C, c, c_1,c_2,...$ are universal constants that may vary.

\section{Rate-Optimal Estimation of Fully Shared Singular Subspace via Stack-SVD}\label{identical}

\subsection{Minimax upper and lower bounds}

In this section, we establish the minimax upper and lower bounds for estimating the shared singular subspace when multiple signal matrices share identical subspaces, thereby demonstrating the rate-optimality of the Stack-SVD procedure. We also analyze the SNR conditions necessary for the consistency of Stack-SVD estimators and compare with the traditional SVD estimators based on individual matrices. To streamline our presentation, we start with the two-matrix case and then extend our results to multiple matrices.

Throughout our theoretical analysis, we assume the entries of the noise matrices $Z_1$ and $Z_2$ are drawn i.i.d from a zero-mean sub-Gaussian distribution.  Specifically, the noise distribution class $\mathscr{G}_\tau$ is defined as follows: for some  constant $c>0$,
$\mathscr{G}_\tau = \{Z \in \mathbb{R}: \mathbb{E} Z=0, \mathbb{E}Z^2=\tau^2, \mathbb{E} \exp(tZ)\leq \exp(\delta^2 t^2), \forall t \in \mathbb{R} \}$,
where the sub-Gaussian parameter $\delta$ is a constant that depends on the variance $\tau$.    Without loss of generality, we take the variance of noise as $\tau=1$.

For the low-rank signal matrices $X_i$ in the observation model (\ref{model}), we consider the following parameter space,
\begin{align}\label{eq:paraspace1}
\mathscr{F}_{r,\gamma} = \left\{\begin{array}{c}X=(X_1~X_2) \in \mathbb{R}^{n \times (p_1+p_2)}:\\
       \text{rank}(X_i)=r, X_i = U \Sigma_i V_i^T\in\mathbb{R}^{n\times p_i},  i=1,2,\\
      \mathop{\min}\limits_{1 \leq i \leq r}\{\sigma_{(i)}^2(X_1)+\sigma_{(i)}^2(X_2)\} \geq \gamma^2
\end{array} \right\},
    \end{align}
where $U\in \mathbb{O}(n,r)$, $V_i\in \mathbb{O}(p_i,r)$, and $\Sigma_i=\text{diag}(\sigma_{(1)}(X_i),...,\sigma_{(r)}(X_i))$. Here $\sigma_{(i)}(X_j)$ denotes the singular value of $X_j$ associated with the $i$-th column of $U$, which means $\sigma_{(i)}(X_j)$ are not necessarily in a decreasing order. In particular, the parameter $\gamma$ indicates the overall signal strength of the stacked matrix $X$, which plays an important role in the subsequent minimax risk analysis. Unlike the singular subspace estimation based on a single matrix, where the overall signal strength is determined by the minimum singular value of the low-rank signal matrix \citep{cai2021subspace}, for estimation across multiple matrices, our analysis reveals  that the overall signal strength would depend on all the singular values across multiple matrices in a non-trivial way. The following result provides a theoretical guarantee for the performance of Stack-SVD by establishing the risk upper bounds and minimax lower bounds for two distance metrics between subspaces.

\begin{theorem}\label{upperlowerbound}
    Suppose $Y_i,i=1,2,$ are generated from (\ref{model}), where the noise matrices $Z_i$ are $i.i.d.$ generated from $\mathscr{G}_\tau$. Define $\hat U$ as the Stack-SVD estimator, whose columns are the first $r$ left singular vectors of $Y=(Y_1~Y_2)\in\mathbb{R}^{n\times (p_1+p_2)}$, and denote $X=(X_1~X_2)\in\mathbb{R}^{n\times (p_1+p_2)}$. Then
    there exists a constant $c>0$ that only depends on $\tau$,  such that:
    \begin{equation}\label{up1.share}
        \mathop{\sup}\limits_{X\in \mathscr{F}_{r,\gamma}}\mathbb{E}\|\sin\Theta(U,\hat{U})\|^2 \leq \frac{cn(\gamma^2+(p_1+p_2))}{\gamma^4}\wedge 1,
    \end{equation}
\begin{equation}\label{up2.share}
        \mathop{\sup}\limits_{X\in \mathscr{F}_{r,\gamma}}\mathbb{E}\|\sin\Theta(U,\hat{U})\|_F^2 \leq \frac{cnr(\gamma^2+(p_1+p_2))}{\gamma^4}\wedge r.
    \end{equation}
  If we further assume $r \leq \frac{n}{2} \wedge \frac{p_1}{16}\wedge \frac{p_2}{16}$, and either $p_1\asymp p_2$ or $\gamma^2\gtrsim p_1+p_2$ holds, then there exists a small constant $c>0$ that only depends on $\delta$, such that:
    \begin{equation}\label{low1.share}
        \mathop{\inf}\limits_{\Tilde{U}}\mathop{\sup}\limits_{X\in \mathscr{F}_{r,\gamma}} \mathbb{E}\|\sin\Theta(U,\Tilde{U})\|^2 \geq c\bigg(\frac{n(\gamma^2+(p_1+p_2))}{\gamma^4} \wedge 1\bigg),
    \end{equation}
\begin{equation}\label{low2.share}
        \mathop{\inf}\limits_{\Tilde{U}}\mathop{\sup}\limits_{X\in \mathscr{F}_{r,\gamma}} \mathbb{E}\|\sin\Theta(U,\Tilde{U})\|_F^2 \geq c\bigg(\frac{nr(\gamma^2+(p_1+p_2))}{\gamma^4} \wedge r\bigg).
    \end{equation}
\end{theorem}

An estimator is considered minimax optimal if its worst-case risk achieves the minimax lower bound up to a constant, characterizing the fundamental difficulty of the statistical problem. Formally, for a risk function $R(\theta, \Tilde{\theta})$, true parameters $\theta\in \mathscr{F}$ and a class of estimators $\mathscr{D}$, the minimax risk is defined as the infimum of the maximum risk: $\inf_{\Tilde{\theta} \in \mathcal{D}} \sup_{\theta \in \mathscr{F}} R(\theta, \Tilde{\theta})$. An estimator $\hat{\theta}$ is minimax rate-optimal if its maximum risk achieves the minimax risk up to a constant, that is, $\sup_{\theta \in \mathscr{F}} R(\theta, \hat{\theta}) \leq c \inf_{\Tilde{\theta} \in \mathcal{D}} \sup_{\theta \in \mathscr{F}} R(\theta, \Tilde{\theta})$ for some constant $c>0$. This means that no other estimator can perform better in convergence rate than $\hat\theta$ in the worst-case, and $\hat{\theta}$ is therefore a best possible estimator against the most adversarial parameter.

The first part of 
Theorem \ref{upperlowerbound} concerns the risk upper bounds (\ref{up1.share}) and (\ref{up2.share}) for the Stack-SVD estimator $\hat U$ over the parameter space $\mathscr{F}_{r,\gamma}$, whereas the second part establishes the minimax risk lower bounds (\ref{low1.share}) and (\ref{low2.share}) for estimating of $U$ over $\mathscr{F}_{r,\gamma}$ by any estimator. The minimax lower bounds delineate the fundamental limit on the estimation accuracy achieved by any estimator. Comparing these lower and upper bounds, we conclude that the Stack-SVD estimator achieves the optimal rate of convergence whenever the dimensions $p_1$ and $p_2$ are comparable ($p_1\asymp p_2$), or the signal strength is sufficiently large ($\gamma^2\gtrsim p_1+p_2$), thereby justifying its efficacy in estimating the shared singular subspaces in many applications. Note that Theorem \ref{upperlowerbound} does not require the rank $r$ to be finite. Our proof of Theorem \ref{upperlowerbound} relies on a novel lower-bound argument for estimating singular subspaces from stacked matrices, which directly links the fundamental limits to the properties of the individual component matrices. Unlike prior lower-bound arguments, such as those in \cite{caizhang2018,cai2021subspace}, our approach addresses a more intricate parameter space, highlighting the contributions and roles of the individual submatrices. See Section \ref{discussionminimax} for more discussion.

\subsection{Comparison with individual subspace estimators}

Theorem \ref{upperlowerbound} establishes the rate-optimality of Stack-SVD when the dimensions of the matrices  $p_1$  and  $p_2$  are of comparable size or overall signal strength is sufficiently large. However, when the orders of $p_1$ and  $p_2$  differ substantially, with one being much larger than the other, we find that there is no need to integrate the two matrices. Instead, under some signal constraints, the minimax optimal rate can  be achieved by individual SVD on the smaller dimension matrix.

\begin{corollary}\label{smalloptimal}
    Suppose the conditions of Theorem \ref{upperlowerbound} hold, and that there exists a constant $c>0$ such that $\sigma_r^2(X_1)\ge c\mathop{\min}\limits_{1 \leq i \leq r}\{\sigma_{(i)}^2(X_1)+\sigma_{(i)}^2(X_2)\}$. Define $\hat U^{(1)}\in \mathbb{O}(n,r)$ whose columns are the top $r$ left singular vectors of $Y_1$. Then, we have for constant $c>0$ that only depends on $\delta$, such that $$
        \mathop{\sup}\limits_{X\in \mathscr{F}_{r,\gamma}}\mathbb{E}\|\sin\Theta(U,\hat{U}^{(1)})\|^2 \leq \frac{cn(\gamma^2+p_1)}{\gamma^4}\wedge 1,\qquad 
        \mathop{\sup}\limits_{X\in \mathscr{F}_{r,\gamma}}\mathbb{E}\|\sin\Theta(U,\hat{U}^{(1)})\|_F^2 \leq \frac{cnr(\gamma^2+p_1)}{\gamma^4}\wedge r.$$
    If we further assume $r \leq \frac{n}{2} \wedge \frac{p_1}{16}\wedge \frac{p_2}{16}$ and  $p_1\le Cp_2$, then for constant $c>0$  that only depends on $\delta$, we have $$
        \mathop{\inf}\limits_{\Tilde{U}}\mathop{\sup}\limits_{X\in \mathscr{F}_{r,\gamma}} \mathbb{E}\|\sin\Theta(U,\Tilde{U})\|^2 \geq c\bigg(\frac{n(\gamma^2+p_1)}{\gamma^4} \wedge 1\bigg),\qquad
        \mathop{\inf}\limits_{\Tilde{U}}\mathop{\sup}\limits_{X\in \mathscr{F}_{r,\gamma}} \mathbb{E}\|\sin\Theta(U,\Tilde{U})\|_F^2 \geq c\bigg(\frac{nr(\gamma^2+p_1)}{\gamma^4} \wedge r\bigg).$$

\end{corollary}

Corollary \ref{smalloptimal} demonstrates that when the dimensions  $p_1$  and  $p_2$  are imbalanced, and the matrix (say, $Y_1$) with the smaller dimension exhibits a sufficiently large signal strength, then leading $r$ singular vectors ($\hat U^{(1)}$) of that matrix achieves the minimax optimal rates: 
$$
\mathop{\sup}\limits_{X\in \mathscr{F}_{r,\gamma}}\mathbb{E}\|\sin\Theta(U,\hat{U}^{(1)})\|_F^2 \asymp\mathop{\inf}\limits_{\Tilde{U}}\mathop{\sup}\limits_{X\in \mathscr{F}_{r,\gamma}} \mathbb{E}\|\sin\Theta(U,\Tilde{U})\|_F^2 \asymp \frac{nr(\gamma^2+p_1)}{\gamma^4} \wedge r.
$$
When  $p_1 \asymp p_2$  and the signal strength of one matrix $Y_1$ is particularly large so that $$\sigma_{r}^2(X_1)\gtrsim\mathop{\min}\limits_{1 \leq i \leq r}\{\sigma_{(i)}^2(X_1)+\sigma_{(i)}^2(X_2)\},$$ then both $\hat U^{(1)}$ and Stack-SVD estimator $\hat U$ achieve the minimax optimal rates. The results are further illustrated in the simulation studies in Section \ref{simulation}. Conversely, when the SNR of one matrix $Y_2$ is particularly low with a very large $p_2$, stacking the two matrices would reduce the overall SNR, particularly when $\gamma^2\ll p_1+p_2$. In such cases, the individual SVD based solely on $Y_1$ can yield a better estimate, provided that $\gamma^2\gg p_1$.

However, Corollary \ref{smalloptimal} only focuses on the special case of imbalanced SNR, as represented by the condition $$\sigma_r^2(X_1)\ge c\mathop{\min}\limits_{1 \leq i \leq r}\{\sigma_{(i)}^2(X_1)+\sigma_{(i)}^2(X_2)\},$$ and does not imply that individual SVD estimator is always preferable than Stack-SVD. To better elucidate the  advantages of integrating multiple data matrices through Stack-SVD, we compare the minimum SNR conditions for the consistency of Stack-SVD with those of the individual SVD estimators. On the one hand, by Theorem \ref{upperlowerbound}, there is a phase transition in the consistency of the Stack-SVD estimator in terms of the SNR $\mathop{\min}\limits_{1 \leq i \leq r}\{\sigma_{(i)}^2(X_1)+\sigma_{(i)}^2(X_2)\}/{\tau^2}$, which  is just the signal $\mathop{\min}\limits_{1 \leq i \leq r}\{\sigma_{(i)}^2(X_1)+\sigma_{(i)}^2(X_2)\}$ when assuming $\tau=1$, with the critical point at
$${\mathop{\min}\limits_{1 \leq i \leq r}\{\sigma_{(i)}^2(X_1)+\sigma_{(i)}^2(X_2)\}} \asymp \sqrt{n(n+p_1+p_2)}.$$
More specifically, when $\mathop{\min}\limits_{1 \leq i \leq r}\{\sigma_{(i)}^2(X_1)+\sigma_{(i)}^2(X_2)\}\gg \sqrt{n(n+p_1+p_2)}$, the right-hand-side of (\ref{up1.share}) and (\ref{up2.share}) converges to $0$, thereby ensuring the consistency of $\hat U$; when $\mathop{\min}\limits_{1 \leq i \leq r}\{\sigma_{(i)}^2(X_1)+\sigma_{(i)}^2(X_2)\}\lesssim \sqrt{n(n+p_1+p_2)}$, the minimax lower bounds (\ref{low1.share}) and (\ref{low2.share}) indicate the impossibility of consistent estimation.
On the other hand, for the subspace estimation based on a single data matrix, the results in \cite{caizhang2018} state that the phase transition for consistently estimating the singular subspaces of $X_1$ and $X_2$  happens at
${\sigma_r^2(X_1)} \asymp \sqrt{n(n+p_1)},$ and $\quad{\sigma_r^2(X_2)} \asymp \sqrt{n(n+p_2)},$
respectively. In other words, consistent estimation of the shared singular subspace based on $X_i$ alone is possible only if ${\sigma_r^2(X_i)} \gg \sqrt{n(n+p_i)}$.

It can be seen that the minimum SNR condition can be more easily achieved with Stack-SVD. Specifically, we elaborate this point by considering an interesting scenario where the SNR ${\sigma_r^2(X_1)}$ in $Y_1$ required by $\hat U^{(1)}$ is below its own  critical point $\sqrt{n(n+p_1)}$, whereas the SNR in $Y_2$ required by $\hat U^{(2)}$ is above its critical point $\sqrt{n(n+p_2)}$. In this case, only one of the two matrices can lead to consistent subspace estimation, but the Stack-SVD based on both matrices may still lead to consistent estimation. For example, if $p_2\gtrsim p_1$, $n\gtrsim p_1$, and $${\sigma_r^2(X_1)}=c\sqrt{n(\sqrt{n}+p_1)}\ll \sqrt{n(n+p_1)},\qquad {\sigma_r^2(X_2)}\gg\sqrt{n(n+p_2)},$$
then  we have
$$
    {\mathop{\min}\limits_{1 \leq i \leq r}\{\sigma_{(i)}^2(X_1)+\sigma_{(i)}^2(X_2)\}}\geq {(\sigma_{r}^2(X_1)+\sigma_{r}^2(X_2))} \gg
    \sqrt{n(n+p_1+p_2)}.
$$
 As such, there is no need to determine which matrix among the two has the stronger SNR. Consistent (and optimal) estimation is automatically achieved by Stack-SVD that integrates both matrices.

\begin{figure}[h]
    \centering
    \begin{minipage}{0.45\textwidth}
        \centering
        \includegraphics[width=\linewidth]{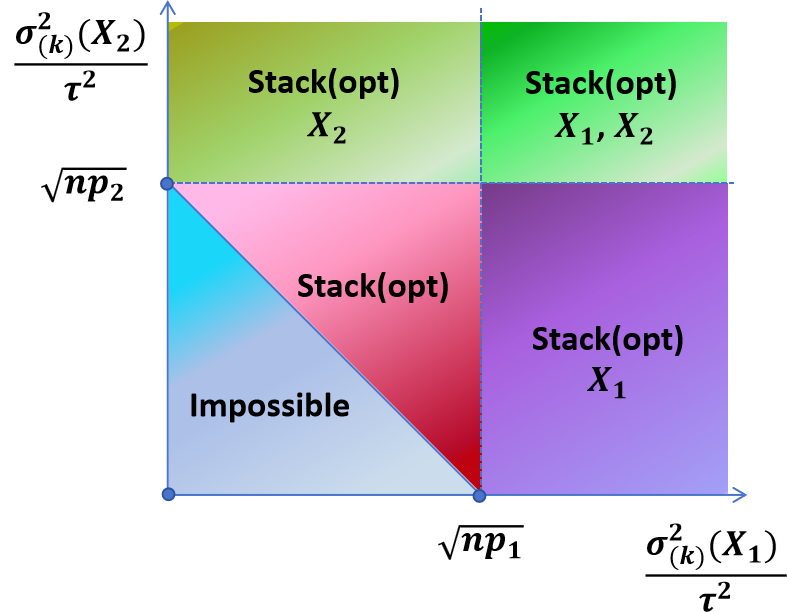}  
    \end{minipage}
    \begin{minipage}{0.45\textwidth}
        \centering
        \includegraphics[width=\linewidth]{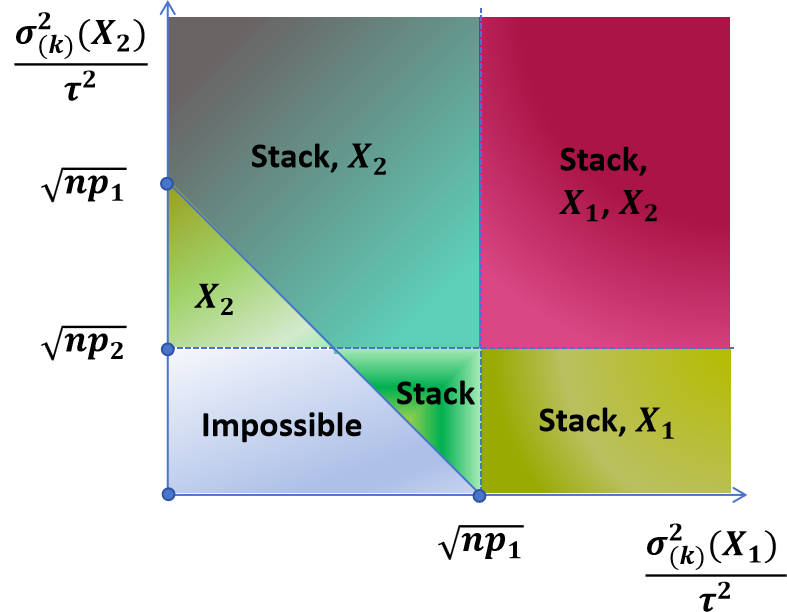}
    \end{minipage}
    \caption{Phase diagrams  when \( p_1 \asymp p_2 \gg n \) (Left) and when \( p_1 \gg p_2 \gg n \) (Right). Each region is labeled by the name of consistent estimators: "$X_i$" indicates $\hat U^{(i)}$ is consistent, "Stack" indicates $\hat U$ is consistent, "Stack(opt)" indicates $\hat U$ is minimax rate-optimal, and "impossible" indicates no consistent estimator exists.}\label{phase}
    \vspace{1em}
\end{figure}

The above observations regarding Stack-SVD can be more systematically demonstrated through the phase diagrams as shown in Figure \ref{phase}. Suppose here $k=\mathop{\arg\min}\limits_{1 \leq i \leq r}\{\sigma_{(i)}^2(X_1)+\sigma_{(i)}^2(X_2)\}$.
On the left of Figure \ref{phase}, when datasets are both high-dimensional ($\min(p_1,p_2)\gg n$) and $p_1$ and $p_2$ are of the same order, Stack-SVD is simultaneously rate-optimal whenever consistent estimation is possible (above the solid line) whereas individual SVD only achieve consistency over smaller regions (above or to the right of the dashed line). Importantly, there is a region (in red) in which Stack-SVD achieves optimal estimation, whereas neither $\hat U^{(1)}$ nor $\hat U^{(2)}$ is consistent. Moreover, on the right of Figure \ref{phase}, when $p_1\gg p_2$, the individual estimator $\hat U^{(2)}$ has strictly weaker SNR requirement (above the horizontal dashed line) compared to $\hat U^{(1)}$, whereas Stack-SVD still shows advantage when both $\hat U^{(1)}$ or $\hat U^{(2)}$ are inconsistent (green region). In practice, when it is unclear which dataset has higher SNR, Stack-SVD can be adopted a practical and oftentimes optimal solution.

\subsection{Extension to multiple matrices}
From Lemma \ref{singularvaluerelation}, it is evident that the above results can be extended to the estimation based on finite $k>2$ matrices. The following corollary demonstrates this generalization, whose proof can be derived in a manner analogous to the proofs of Proposition \ref{generalizecainzhanglowershare} and Theorem \ref{upperlowerbound}. 

\begin{corollary}\label{kmat}
For a finite integer $k>2$, suppose $Y_i,i=1,2,...,k,$ are generated from (\ref{model}), where the noise matrices $Z_i$ are $i.i.d.$ generated from $\mathscr{G}_\tau$. Define $\hat U$ as the Stack-SVD estimator, whose columns are the first $r$ left singular vectors of $Y=(Y_1~Y_2~...~Y_k)\in\mathbb{R}^{n\times \sum_{j=1}^kp_j}$, and denote $X=(X_1~X_2~...~X_k)\in\mathbb{R}^{n\times \sum_{j=1}^kp_j}$. Define
    \begin{align*}
\mathscr{F}_{r,\gamma}^k = \left\{\begin{array}{c}X=(X_1~X_2~...~X_k) \in \mathbb{R}^{n \times \sum_{j=1}^kp_j}:\\
       \textup{rank}(X_i)=r, X_i = U \Sigma_i V_i^T\in\mathbb{R}^{n\times p_i},\\
      \mathop{\min}\limits_{1 \leq s \leq r}\{\sum_{j=1}^k\sigma_{(s)}^2(X_j)\} \geq \gamma^2, i=1,2,...k
\end{array} \right\}.
    \end{align*}
where $U\in \mathbb{O}(n,r)$, $V_i\in \mathbb{O}(p_i,r)$, and $\Sigma_i=\text{diag}(\sigma_{(1)}(X_i),...,\sigma_{(r)}(X_i))$. Then there exists a constant $c>0$ that only depends on $\delta$,  such that:
    $$
        \mathop{\sup}\limits_{X\in \mathscr{F}^k_{r,\gamma}}\mathbb{E}\|\sin\Theta(U,\hat{U})\|^2 \leq \frac{cn(\gamma^2+(\sum_{j=1}^kp_j))}{\gamma^4}\wedge 1,
\qquad
        \mathop{\sup}\limits_{X\in \mathscr{F}^k_{r,\gamma}}\mathbb{E}\|\sin\Theta(U,\hat{U})\|_F^2 \leq \frac{cnr(\gamma^2+(\sum_{j=1}^kp_j))}{\gamma^4}\wedge r.$$
If we further assume $r \leq \frac{n}{k} \wedge \frac{p_1}{16}\wedge \frac{p_2}{16}\cdots\wedge \frac{p_k}{16}$, and either $p_1\asymp p_2\asymp...\asymp p_k$ or $\gamma^2\gtrsim \sum_{j=1}^kp_j$ holds, then there exists a small constant $c>0$ that only depends on $\delta$, such that: $$
        \mathop{\inf}\limits_{\Tilde{U}}\mathop{\sup}\limits_{X\in \mathscr{F}^k_{r,\gamma}} \mathbb{E}\|\sin\Theta(U,\Tilde{U})\|^2 \geq c\bigg(\frac{n(\gamma^2+(\sum_{j=1}^kp_j))}{\gamma^4} \wedge 1\bigg),\qquad
        \mathop{\inf}\limits_{\Tilde{U}}\mathop{\sup}\limits_{X\in \mathscr{F}^k_{r,\gamma}} \mathbb{E}\|\sin\Theta(U,\Tilde{U})\|_F^2 \geq c\bigg(\frac{nr(\gamma^2+(\sum_{j=1}^kp_j))}{\gamma^4} \wedge r\bigg).$$
\end{corollary}

The corollary indicates that when multiple matrices share the identical singular subspace, Stack-SVD still  achieves the minimax rate-optimal estimation.

\subsection{Sub-Optimality of Average-SVD} \label{suboptimal.sec}

In this section, we define and prove the sub-optimality of Average-SVD method, which is also referred to as AJIVE \cite{aJIVE_congma}, Distributed PCA \cite{fan2017distributed,chen2022distributed}, or Principal Angle Analysis \cite{aJIVE}. 

We first demonstrate the equivalence of Average-SVD and AJIVE.  Here suppose we have $k$ noisy matrices $Y_1,...,Y_k$. The signal matrices $X_1,...,X_k$ are all of rank $r$ and share the left singular subspace. First we take SVD and denote
$$
Y_1=U_1\Sigma_1 V_1^T,~...,~Y_k=U_k\Sigma_k V_k^T$$
 with singular values in $\Sigma_i$ decreasingly ordered. 
We take the top $r$ singular vectors of $U_i$ as $U_i^{(r)}$.
Then an average is taken
$$
A=\frac{U_1^{(r)}U_1^{(r)T}+...+U_k^{(r)}U_k^{(r)T}}{k}.
$$
The Average-SVD estimator $\hat{U}$ can be defined as the top  $r$ left singular vectors of the average matrix $A$.
The basic idea of AJIVE \cite{aJIVE_congma} is to take SVD to the concatenation of individual singular subspace, that is, the SVD of 
$
B=(U_1^{(r)}~U_2^{(r)}~...~U_k^{(r)}).
$
The AJIVE estimates the shared singular subspace by the top $r$ left singular vectors from the SVD of the concatenation matrix $B$.
Note that the idea of Average-SVD and AJIVE are actually equivalent. We define 
$$
C=BB^T=(U_1^{(r)}~U_2^{(r)}~...~U_k^{(r)})(U_1^{(r)}~U_2^{(r)}~...~U_k^{(r)})^T=U_1^{(r)}U_1^{(r)T}+...+U_k^{(r)}U_k^{(r)T}.
$$
The matrix $B$ and matrix $C=BB^T$ have the same left singular subspace. Also note that $C=kA$. The two matrices $C$ and $A$ are only different up to a constant, which means they have identical left singular subspace. Therefore, the Average-SVD and AJIVE are actually equivalent.

For the sub-optimality of AJIVE or Average-SVD, we begin by considering a rank-1 case where the matrices are defined as
$X_1 = \alpha  u v_1^T$ and  $X_2 = \beta  u v_2^T$, where $u$, $v_1$, $v_2$ are vectors and $\alpha$, $\beta$ are positive numbers. Let \( Y_1 \) and \( Y_2 \) represent the noisy version of \(X_1\) and \(X_2\), respectively, so that
$Y_1 = X_1 + Z_1$ and $Y_2 = X_2 + Z_2$, where \( Z_1 \) and \( Z_2 \) represent random noise.

After sorting the singular values in a decreasing order, the first left singular vectors of \( Y_1 \) and \( Y_2 \) are denoted as \( \hat{u}_1 \) and \( \hat{u}_2 \), respectively. To estimate the shared singular vector \( u \), the AJIVE or the Principal Angle Analysis estimator, denoted as \( \hat{u}_{ave} \), is defined as the top left singular vector of the stacked singular subspace formed by \( (\hat{u}_1 ~ \hat{u}_2) \). Equivalently, it is also the top eigenvector of the mean projection matrix 
$\frac{\hat{u}_1 \hat{u}_1^T + \hat{u}_2 \hat{u}_2^T}{2}$.

For the estimation of the shared vector \( u \) using $\hat{u}_{ave}$, we have the following lower bound result.

\begin{theorem}\label{suboptimalave}
Suppose we have balanced dimensions $(p_1\asymp p_2)$, imbalanced signals $(\beta \gg \alpha)$, and low rank conditions $r=1\leq \frac{n}{2}\wedge \frac{p_1}{16}$, $r=1\leq \frac{n}{2}\wedge \frac{p_2}{16}$. In addition, we require signals $\alpha$, $\beta$ are sufficiently large. When estimating within the parameter space
\begin{align*}
\mathscr{F}_{1,\gamma} = \left\{\begin{array}{c}X=(X_1~X_2) \in \mathbb{R}^{n \times (p_1+p_2)}:\\
       X_1=\alpha u v_1^T, X_2=\beta uv_2^T, \alpha^2+\beta^2 \geq \gamma^2,\\
       u\in \mathbb{O}(n,1), v_i\in \mathbb{O}(p_i,1), i=1,2
\end{array} \right\},
    \end{align*}
it holds that
    \begin{align*}
    \mathop{\sup}\limits_{X\in \mathscr{F}_{1,\gamma}}\mathbb{E}&||\sin\Theta(u,\hat{u}_{ave})||^2 \geq \left(c_1\sqrt{\frac{n(\alpha^2+p_1)}{\alpha^4}}-c_2\sqrt{\frac{n(\beta^2+p_2)}{\beta^4}}\right)^2.
    \end{align*}
\end{theorem}

Combining the minimax rate derived in \cite{caizhang2018}, the results above suggest that the convergence rate of Average-SVD is no smaller than the difference in convergence rates between the individual SVDs of $Y_1$ and $Y_2$. Consequently, as the difference between the rates becomes larger, Average-SVD is expected to perform sub-optimally. On the other hand, when this difference is small, or the two data matrices exhibit balanced signal strengths, the Average-SVD method, as explored in the recent work \cite{aJIVE_congma}, has been shown to achieve minimax optimality (over a smaller parameter space).

To demonstrate the sub-optimality of Average-SVD when the signal strengths in the different matrices are imbalanced, we consider the following scenario. Specifically, when the disparity between the signals is large, i.e., when \( \beta \gg \alpha \), even if dimensions are comparable \( (p_1 \asymp p_2) \) and both signals \( \alpha \) and \( \beta \) are individually large, we still observe that
$$
\left(c_1\sqrt{\frac{n(\alpha^2+p_1)}{\alpha^4}}-c_2\sqrt{\frac{n(\beta^2+p_2)}{\beta^4}}\right)^2 \asymp c\frac{n(\alpha^2+p_1)}{\alpha^4} 
\gg c\frac{n(\alpha^2+\beta^2+p_1+p_2)}{(\alpha^2+\beta^2)^2}.
$$
Here, the last term represents the rate for Stack-SVD. This analysis suggests that AJIVE can be suboptimal under imbalanced signal strengths, even when both signals are strong, whereas Stack-SVD remains optimal.

\subsection{Technical tool for minimax lower bounds in the fully shared subspace setting}\label{discussionminimaxfullyshared}

The process of stacking matrices imposes extra structural constraints on the stacked matrix. To facilitate the proof of Theorems \ref{upperlowerbound}, and to provide additional insights on the minimax lower bounds, we have the following proposition, which builds upon a general argument that explores a structured parameter space specifically designed for matrix-stacking.

\begin{proposition}\label{generalizecainzhanglowershare}
    Suppose $X \in \mathbb{R}^{p_1 \times p_2}$ has SVD 
    \begin{equation}\label{prop1SVD}
        X=U \sqrt{\Sigma_1^2+\Sigma_2^2} \bigg(\frac{\Sigma_1}{\sqrt{\Sigma_1^2+\Sigma_2^2}}V_1^T~\frac{\Sigma_2}{\sqrt{\Sigma_1^2+\Sigma_2^2}}V_2^T\bigg) ,
    \end{equation}
    where $\Sigma_1,\Sigma_2$ are rank $r$ positive diagonal matrices and $V_1 \in \mathbb{O}(p_{V_1},r),V_2 \in \mathbb{O}(p_{V_2},r)$ are respective orthonormal matrices with dimension $p_{V_1}+p_{V_2}=p_2$. Assume $r \leq \frac{p_1}{2} \wedge \frac{p_{V1}}{16} \wedge \frac{p_{V2}}{16}$ and either $p_{V1}\asymp p_{V2}$ or $\gamma^2\gtrsim p_{V1}+p_{V2}$ holds. The noise $ Z \overset{\text{i.i.d.}}{\sim}\mathscr{G}_\tau$. When estimating within 
     \begin{align}\label{eq:101}
\mathscr{H}_{r,\gamma}^{(s)} = \left\{\begin{array}{c}X \in \mathbb{R}^{p_1 \times p_2}:\\
       p_2=p_{V_1}+p_{V_2},\text{rank}(X)=r,\sigma_{\min}(X)\geq \gamma\\
      X~has~form~(\ref{prop1SVD}),U \in \mathbb{O}(p_1,r), \\V_i \in \mathbb{O}(p_{V_i},r),i=1,2
\end{array} \right\},
    \end{align}
    it holds for some constant $c>0$  that only depends on $\delta$ that 
    $$
        \mathop{\inf}\limits_{\Tilde{U}}\mathop{\sup}\limits_{X\in \mathscr{H}^{(s)}_{r,\gamma}} \mathbb{E}\|\sin\Theta(U,\Tilde{U})\|^2 \geq c\bigg(\frac{p_1(\gamma^2+(p_{V_1}+p_{V_2}))}{\gamma^4} \wedge 1\bigg),~\mathop{\inf}\limits_{\Tilde{U}}\mathop{\sup}\limits_{X\in \mathscr{H}^{(s)}_{r,\gamma}} \mathbb{E}\|\sin\Theta(U,\Tilde{U})\|_F^2 \geq c\bigg(\frac{p_1r(\gamma^2+(p_{V_1}+p_{V_2}))}{\gamma^4} \wedge r\bigg).
    $$    
\end{proposition}

Compared to Theorem 4 in \cite{caizhang2018}, Proposition \ref{generalizecainzhanglowershare} shows that when \( p_{V1} \) and \( p_{V2} \) are of comparable size, or when the signal strength is sufficiently large, the structural constraint we impose on the right singular subspace does not affect the final minimax lower bound.

\section{Estimation in the Presence of Unshared Singular Vectors}\label{unshared}
In many applications, different matrices may not share an identical singular subspace. Instead, each matrix possesses its unique singular vectors in addition to the shared subspace. This section explores the estimation of the shared subspace while accounting for the presence of unshared vectors.

\subsection{Minimax upper and lower bounds} \label{sec.3.1}

In this section, we analyze the scenario in which the subspaces of the signal matrices consist of both shared and unshared singular vectors. We first derive both upper bounds and matching minimax lower bounds under the condition that the unshared singular vectors are orthogonal to each other. Similar setup has been considered previously \cite{aJIVE,RaJIVE}, but the fundamental limit has not been established. In Section \ref{sec.nonorth}, we will extend our results to the general non-orthogonal settings.
We  first present a result clarifying the relationship between the singular subspace of the individual matrices and that of the stacked matrix, which is important in our analysis.
\begin{proposition}\label{vectorrelation}
    Let the SVDs of \( X_1 \), \( X_2 \), and their stacked version \((X_1~X_2)\) be given as 
$X_1 = U_1 \Sigma_1 V_1^T,  X_2 = U_2 \Sigma_2 V_2^T,$ and $(X_1~X_2) = U \Sigma V^T.$ \( U_1 \) and \( U_2 \) are allowed to share some identical singular vectors, and all distinct singular vectors in \( U_1 \) and \( U_2 \) are pairwise orthogonal. Selecting any column vector(s) from $U_1$ and $U_2$, represented as $U_1^*$ and $U_2^*$, the singular subspace spanned by $U_1^*$ and $U_2^*$ are both contained in the left singular subspace of the stacked matrix $(X_1~X_2)$. That is
$\text{span}\{U_1^*\}\subset \text{span}\{U\}, \text{span}\{U_2^*\}\subset \text{span}\{U\}.$
As a special case, when further assuming that the singular values of \( X_1 \), \( X_2 \), and \((X_1~X_2)\) are all distinct, then every singular vector in \( U_1 \) and \( U_2 \) is contained in the columns of \( U \).
\end{proposition}

Proposition \ref{vectorrelation} presents an intuition that with the orthogonality constraint on the unshared singular vectors, all singular vectors—-whether shared or unshared—will be present in the singular subspace of the stacked signal matrix $(X_1~X_2)$, albeit potentially in an intermixed order. For example, suppose we are interested in estimating the shared singular vectors $(u_1~u_2)$ in
   $$
X_1=(
    u_1 ~ u_2 ~ u_{1*})\mathrm{diag}(
    2\alpha,~\alpha,~\frac{3}{2}\sqrt{\alpha^2+\beta^2})V_1^T, \qquad X_2=(
    u_1 ~ u_2 ~ u_{2*})\mathrm{diag}(
    2\beta,~\beta,~\frac{1}{2}\sqrt{\alpha^2+\beta^2})V_2^T,
$$
where $\alpha,\beta>0$ and $u_{1*}^\top u_{2*}=0$.
It follows that the SVD of $(X_1~X_2)$ is
$$(X_1~X_2)=(    u_1 ~  u_{1*} ~ u_2 ~ u_{2*})\mathrm{diag}(2\sqrt{\alpha^2+\beta^2},~\frac{3}{2}\sqrt{\alpha^2+\beta^2},~\sqrt{\alpha^2+\beta^2},\\~\frac{1}{2}\sqrt{\alpha^2+\beta^2})V^T.$$
Comparing the orders of the singular vectors in $X_1$, $X_2$ and $(X_1~X_2)$, especially the order of $u_{1*}$ and $u_2$ before and after stacking, we see that the singular vectors of $X_1$ and $X_2$ can be intermixed after stacking, depending on the relations between the singular values. Thus, to reliably estimate the shared singular vectors, it is important to trace their locations in the SVD of the stacked matrix.

Before defining the parameter space, we introduce some notations. We say that the shared (left) singular vectors and unshared (left) singular vectors are in different \emph{vector types}. For the (left) singular vectors of the stacked matrix $(X_1~X_2)$, ordered based on the singular values of $(X_1~X_2)$, we say that there is a \emph{vector type switch} at $s$-th position if, the $s$-th and $(s+1)$-th singular vectors associated to the $s$-th and $(s+1)$-th largest singular values are in different vector types. Suppose there are $N$ vector type switches in $X=(X_1~X_2)$, and we denote $\mathbb{S}$ as the index set of the locations where the switch happens. The eigen-gap at the $i$-th switch that happens at $s$-th position is denoted as 
$
g_{i(s)}^2(X):=\sigma_s^2(X)-\sigma_{s+1}^2(X),~ s\in \mathbb{S},
$ 
where $i(s)\in\{1,2,...,N\}$ is the rank of $s$ in $\mathbb{S}$, sorted in increasing order.
In particular, when a shared vector appears to be the last singular vector of $X$, associated with the smallest singular value $\sigma_{\min}(X)$, we still consider that a vector type switch exists. The eigen-gap is thus defined as 
$g_N^2(X):=\sigma_{\min}^2(X)-0.
$
We should note that in this case, such an eigen-gap is not counted in the index set $\mathbb{S}$.

To effectively distinguish the shared and unshared singular vectors, eigen-gaps should be sufficiently large. Therefore, for any small constant $c>0$, we define the parameter space 
\begin{align}\label{eq:paraspace2}
\mathscr{H}_{r,t} = \left\{\begin{array}{c}X=(X_1~X_2):~\text{rank}(X_i)=r+r_{i*},\\ X_i = (U_r~U_{i*}) \Sigma_i V_i^T\in\mathbb{R}^{n\times p_i}, ~U_{1*}^T U_{2*}=0, \\i=1,2,~g_{i(s)}^2(X)>c\sigma_{s+1}^2(X), \forall s\in\mathbb{S},\\
      \mathop{\min}\limits_{1\le k\le N}g_k^2(X) \geq t^2
\end{array} \right\},
    \end{align}
where the $U_r \in \mathbb{O}(n,r)$ are the shared singular vectors while $U_{i*} \in \mathbb{O}(n,r_{i*})$ are respective unshared singular vectors for $X_1$ and $X_2$. Here, we do not necessarily require the diagonal of $\Sigma_i=\text{diag}(\sigma_{(1)}(X_i),...,\sigma_{(r+r_{i*})}(X_i))$ to have a decreasing order, which means after manually sorting in a decreasing order, the vectors in the shared $U_r$ and the unshared $U_{i*}$ might be shuffled in the left singular subspaces of individual $X_1$ and $X_2$. In particular, if the last singular vector in $(X_1~X_2)$ is a shared singular vector, then we have  $|\mathbb{S}|=N-1$, and the condition on $g_N^2(X)$ is still imposed through $\mathop{\min}\limits_{1\le k\le N}g_k^2(X)$.

Suppose the positions of shared singular vectors in the stacked noiseless matrix are known, and we denote their positions as an index set $\mathbb{J}$. We consider the estimator $\hat{U}^{\mathbb{J}}_r\in \mathbb{R}^{n\times r}$, whose columns are $r$ selected left singular vectors of $(Y_1~Y_2)$ indexed by $\mathbb{J}$ . We have the following upper and lower bounds.

\begin{theorem}\label{generalbounds}
    Suppose $Y_i,i=1,2,$ are generated from (\ref{model}), where the noise matrices $Z_i$ are $i.i.d.$ generated from $\mathscr{G}_\tau$. For any finite $r$, we have for some constant $c>0$ that only depends on $\delta$, such that
    \begin{equation}\label{up1}
        \mathop{\sup}\limits_{X\in \mathscr{H}_{r,t}}\mathbb{E}\|\sin\Theta(U_r,\hat{U}^{\mathbb{J}}_r)\|^2 \leq \frac{cn(t^2+(p_1+p_2))}{t^4}\wedge 1,
    \end{equation}
    \begin{equation}\label{up2}
        \mathop{\sup}\limits_{X\in \mathscr{H}_{r,t}}\mathbb{E}\|\sin\Theta(U_r,\hat{U}^{\mathbb{J}}_r)\|_F^2 \leq \frac{cnr(t^2+(p_1+p_2))}{t^4}\wedge r.
    \end{equation}
    If we further assume $r+r_{1*}+r_{2*} \leq \frac{n}{2}$, $r+r_{1*}\leq \frac{p_1}{16}$, $r+r_{2*}\leq \frac{p_2}{16}$, and either $p_1\asymp p_2$ or $t^2\gtrsim p_1+p_2$ holds, then for constant $c>0$  that only depends on $\delta$, we have 
    \begin{equation}\label{low1}
        \mathop{\inf}\limits_{\Tilde{U}_r}\mathop{\sup}\limits_{X\in \mathscr{H}_{r,t}} \mathbb{E}\|\sin\Theta(U_r,\Tilde{U}_r)\|^2 \geq c\bigg(\frac{n(t^2+(p_1+p_2))}{t^4} \wedge 1\bigg),
    \end{equation}
    \begin{equation}\label{low2}
        \mathop{\inf}\limits_{\Tilde{U}_r}\mathop{\sup}\limits_{X\in \mathscr{H}_{r,t}} \mathbb{E}\|\sin\Theta(U_r,\Tilde{U}_r)\|_F^2 \geq c\bigg(\frac{nr(t^2+(p_1+p_2))}{t^4} \wedge r\bigg).
    \end{equation} 
\end{theorem}

Theorem \ref{generalbounds} asserts that when both shared and unshared vectors are present, the rate-optimal estimation can be achieved, given that the location of the shared vectors can be identified. Compared with Theorem \ref{upperlowerbound}, Theorem \ref{generalbounds} indicates that the minimax-optimal estimation requires both a sufficiently large minimum eigen-gap at the switch location, and a minimum singular value condition. Indeed, under the current assumptions, the minimum singular value $\sigma_{\min}^2(U_{shared})$, among those associated to the singular vectors in $U_{shared}$, and the minimum eigen-gap $\mathop{\min}\limits_{1 \leq k \leq N}g_k^2$ have the same order  $$
\sigma_{\min}^2(U_{shared})\geq g_N^2\geq \mathop{\min}\limits_{1 \leq k \leq N}g_k^2 \geq c\sigma_s^2 \geq c\sigma_{\min}^2(U_{shared}),$$
where $s=\text{argmin}_k g_k^2$. The eigen-gap requirement at the switch points is necessary. When the gap is too small, it becomes difficult to effectively distinguish between shared and unshared singular vectors, resulting in a fundamental limitation in estimation due to non-identifiability. On the other hand, the eigen-gap condition is only required at the switch points. If the vector type remains unchanged, insufficient or even no eigen-gap between the singular vectors does not affect our theoretical results. 

We remark on the key difference between our results and the analysis of AJIVE \cite{JIVE}. Specifically, AJIVE aims to estimate the joint matrix rather than the singular subspace. Consequently, both the left and right singular subspaces contribute valuable information in AJIVE. To account for this, AJIVE employs a perturbation bound within the framework of Wedin’s bound \cite{wedin1972}. This bound is expressed as $\frac{\max(|A|, |B|)}{g}$,
where \( A \) captures information from the left singular subspace, \( B \) pertains to the right singular subspace, and \( g \) quantifies the minimum eigen-gap. 
While this perturbation bound is effective for the AJIVE framework, it proves less effective for our analysis, as we focus on estimating the shared (left) singular subspace alone. The Wedin perturbation bound is uniform for both left and right singular subspaces, which may lead to sub-optimal results when the perturbation affects the two subspaces differently. In contrast, similar to \cite{caizhang2018}, our upper bound is tailored specifically to one side of the singular subspace, leading to sharp results reflecting different roles of $n$ and $p_i$ in the rates of convergence. 
Furthermore, while our work establishes the minimax rate-optimality of the estimator, to the best of our knowledge, no previous work has provided a minimax lower bound for the JIVE framework.

The practical implications of Theorem \ref{generalbounds} are also significant. In many  applications, estimating the shared  subspace for multiple matrices typically involves the use of Stack-SVD, which concatenates all matrices and selects the top $r$ singular vectors from the noisy stacked matrix as the estimation. However, as indicated by Theorem \ref{generalbounds}, the shared and unshared singular vectors may become intermixed after stacking. To achieve optimal estimation, it is crucial to identify the correct index set $\mathbb{J}$. Relying solely on the top $r$ singular vectors from the stacked matrix $(Y_1~Y_2)$ for estimating the shared  subspace may result in biased or inconsistent estimation.

When the dimensions  $p_1$  and  $p_2$  are not of comparable sizes and the signal strength of the matrix with smaller dimension becomes particularly strong, similar results in parallel to Corollary \ref{smalloptimal} can be obtained, where the minimax optimal rates are achieved by different estimators. When one matrix exhibits a considerably larger signal strength and a relatively smaller or equal dimension, the individual SVD estimator of that matrix is also minimax optimal. Besides, when both the signal strength for stacked matrix ($t^2\gtrsim p_1+p_2$) and the individual signal strength of the matrix with smaller dimension are large, the oracle Stack-SVD $\hat{U}_r^{\mathbb{J}}$ and the individual SVD on the matrix with smaller dimension are both rate-optimal.

To obtain the minimax lower bounds in Theorem \ref{generalbounds}, we develop a general argument surrounding a structured parameter space tailored to matrix-stacking; see Section \ref{discussionminimax}.
For the upper bound results in Theorem \ref{generalbounds}, we leverage the following proposition, which generalizes the results in \cite{caizhang2018} to the case of partial subspace estimation.

\begin{proposition}\label{generalcainzhang}
    Suppose $X=U\Sigma V^T \in \mathbb{R}^{p_1\times p_2}$ has rank greater than $r$ and there exists some small constant $\epsilon>0$ such that $\sigma_r^2(X)>(1+\epsilon)\sigma_{r+1}^2(X)$. Let $U_r$ be the top $r$ left singular vectors of $X$ and $\hat U_r$ be the top $r$ left singular vectors of $Y=X+Z$, with $Z\in \mathscr{G}_\tau$. Then for some constant $c>0$  that only depends on $\delta$, we have
    $$
        \mathbb{E}\|\sin\Theta(U_r,\hat{U}_r)\|^2 \leq \frac{cp_1(\sigma_r^2(X)+ p_2)}{\sigma_r^4(X)}\wedge 1,
\qquad
        \mathbb{E}\|\sin\Theta(U_r,\hat{U}_r)\|_F^2 \leq \frac{cp_1r(\sigma_r^2(X)+ p_2)}{\sigma_r^4(X)}\wedge r.
    $$
\end{proposition}

Proposition \ref{generalcainzhang} plays a crucial role here and in our subsequent analysis. It highlights that the accurate estimation of any subspace relies on the separation of its singular values from those of the orthogonal subspace.

Further, Theorem \ref{generalbounds} provides an upper bound on the estimation error that depends on the smallest eigen-gap. This result allows us to generalize the error rate to the estimation of individual singular vectors, as in the following corollary.

\begin{corollary}\label{individualvectorrate}
    For observation $Y=X+Z$, where $X=U\Sigma V^T \in \mathbb{R}^{n\times p}$ is a rank $r$ signal matrix with singular value decreasingly ordered and the noise $Z\in \mathbb{R}^{n\times p}$ are i.i.d. generated from $\mathscr{G}_{\tau}$. For $i=1,..,r-1,$ denote the squared eigen-gap between the $i$-th and $i+1$-th left singular vectors as $g_i^2:=\sigma_i^2(X)-\sigma_{i+1}^2(X)$. Then for the estimation of $u_i$, the $i$-th left singular vector of $X$, the estimator $\hat{u}_i$, the $i$-th left singular vector of noised $Y$, we have the error bound $$\mathbb{E}\|\sin\Theta(u_i,\hat{u}_i)\|^2 \leq \frac{cn(\min\{g_i^2,g_{i-1}^2\}+p)}{\min\{g_i^2,g_{i-1}^2\}^2}\wedge 1,$$
    where we set $g_0^2=\infty$ and $g_{r+1}^2=0$.
\end{corollary}

Corollary \ref{individualvectorrate} establishes the convergence rate for estimating an individual singular vector of a general matrix. This rate is governed by the eigen-gap between the corresponding singular value and its immediate neighbors. Notably, this result applies to general matrix estimation and is not limited to the stacked SVD framework.

\subsection{Technical tool for minimax lower bounds in the partially shared subspace setting}\label{discussionminimax}

The stacking operation imposes complex structural constraints on the resulting matrix. While Proposition \ref{generalizecainzhanglowershare} in Section \ref{discussionminimaxfullyshared} establishes this for the special case of fully shared  subspaces, Proposition \ref{generalizacainzhanglower} then generalizes this result to the broader and more complex scenario of partially shared  subspaces.

\begin{proposition}\label{generalizacainzhanglower}
    Suppose $X \in \mathbb{R}^{p_1 \times p_2}$ has SVD of the form: \begin{align}\label{prop2SVD}
    X=\begin{pmatrix}
        U_r&U_{1*}&U_{2*}
    \end{pmatrix}\begin{pmatrix}
        \sqrt{\Sigma_{1r}^2+\Sigma_{2r}^2}&&\\&\Sigma_{1*}&\\&&\Sigma_{2*}
    \end{pmatrix} \begin{pmatrix}
        \frac{\Sigma_{1r}}{\sqrt{\Sigma_{1r}^2+\Sigma_{2r}^2}}V_{1r}^T & \frac{\Sigma_{2r}}{\sqrt{\Sigma_{1r}^2+\Sigma_{2r}^2}}V_{2r}^T\\
        V_{1*}^T & 0 \\ 0 & V_{2*}^T
    \end{pmatrix},
    \end{align} where $U_r \in \mathbb{O}(p_1, r)$, $U_{i*} \in \mathbb{O}(p_1,r_{i*})$, $V_{ir} \in \mathbb{O}(p_{Vi},r)$, $V_{i*} \in \mathbb{O}(p_{Vi},r_{i*})(i=1,2)$ are orthonormal matrices, $p_{V1}+p_{V2}=p_2$. Vectors in $U_r$ and $U_{i*}$ are all orthogonal to each other. $\Sigma_{ir} \in \mathbb{R}^{r\times r},\Sigma_{i*}\in \mathbb{R}^{r_{i*}\times r_{i*}}(i=1,2)$ are positive diagonal matrices.  $r+r_{1*}+r_{2*} \leq \frac{p_1}{2}$, $r+r_{1*}\leq \frac{p_{V1}}{16}$, $r+r_{2*}\leq \frac{p_{V2}}{16}$. Either $p_{V1}\asymp p_{V2}$ or $\gamma^2\gtrsim p_{V1}+p_{V2}$ holds. The noise $Z \overset{\text{i.i.d.}}{\sim}\mathscr{G}_\tau$. We similarly require a singular value gap condition between different vector types by taking vectors in $U_r$ as a vector type while taking vectors in $U_{i*}$ as another type. The same notation for gaps as in Theorem \ref{generalbounds} in the main text is used. When estimating within 
    \begin{align}\label{eq:102}
\mathscr{H}_{r,\gamma}^{(u)} = \left\{\begin{array}{c}X= \in \mathbb{R}^{p_1 \times p_2}:\\
       p_2=p_{V_1}+p_{V_2},\text{rank}(X)=r+r_{1*}+r_{2*},\\X~has~form~(\ref{prop2SVD}),
      U_r \in \mathbb{O}(p_1,r),\\U_{i*}\in \mathbb{O}(p_1,r_{i*}),U_{1*}^TU_{2*}=0, V_i \in \mathbb{O}(p_{V_i},r),\\
      g_{i(s)}^2(X)>c\sigma_{s+1}^2(X), \forall s\in\mathbb{S},\\
      \mathop{\min}\limits_{1 \leq k \leq r}\sigma_{(k)}^2(X) \geq \gamma^2, i=1,2
\end{array} \right\},
    \end{align}
    where the parenthesis $k$ in $\sigma_{(k)}(X)$ denotes the singular value in $X$ corresponding to the $k$-th vector in the shared $U_r$, for some constant $c>0$  that only depends on $\delta$, we have
    $$
    \mathop{\inf}\limits_{\Tilde{U}}\mathop{\sup}\limits_{X\in \mathscr{H}_{r,\gamma}^{(u)}} \mathbb{E}||\sin\Theta(U_r,\Tilde{U_r})||^2 \geq c(\frac{p_1(\gamma^2+(p_{V1}+p_{V2}))}{\gamma^4} \wedge 1),$$ and 
  $$
  \mathop{\inf}\limits_{\Tilde{U}}\mathop{\sup}\limits_{X\in \mathscr{H}_{r,\gamma}^{(u)}} \mathbb{E}||\sin\Theta(U_r,\Tilde{U_r})||_F^2 \geq c(\frac{p_1r(\gamma^2+(p_{V1}+p_{V2}))}{\gamma^4} \wedge r).$$
\end{proposition}
The implication of Proposition \ref{generalizacainzhanglower} is similar: by imposing the above constraints on the singular subspace that is not being estimated, the final minimax lower bounds remain unchanged, provided that the dimensions are of comparable size or the overall signal strength is large. Further, the structural constraints in the partially shared case are significantly more complex than those in the fully shared case.

\subsection{Two special cases of partial sharing}

To further elaborate the implications of Theorem \ref{generalbounds}, and the advantages and limitations of Stack-SVD, we analyze two special scenarios of partial sharing. 

\subsubsection{Scenario I: estimation with weak unshared signals}

When the singular values associated with the unshared vectors are relatively small compared to those of the shared vectors, it is reasonable to expect that the singular subspace derived from top $r$ Stack-SVD will still yield a reliable estimate for the shared vectors. Furthermore, since both matrices $X_1$	and $X_2$ contain the shared singular vectors, stacking the matrices enhances the signal strength of these shared vectors. In contrast, the unshared vectors, present in only one of the matrices, retain their original signal strength.

We consider the scenario when the singular values $\sigma_i(X)$ of the unshared singular vectors are all smaller than the singular values of the shared singular vectors. Specifically, recall the definition of $\mathbb{J}$ prior to Theorem \ref{generalbounds}; we consider that $$\min_{j\in \mathbb{J}} \sigma^2_{j}(X)\ge (1+\epsilon)\max_{j\notin \mathbb{J}} \sigma^2_{j}(X),$$ for some small $\epsilon>0$. In particular, if $X_1$ and $X_2$ have individual SVDs
$$X_1=\begin{pmatrix}
    U_r & U_{1*}
\end{pmatrix}\Sigma_1 V_1^T\in\mathbb{R}^{n\times p_1},\qquad
X_2=\begin{pmatrix}
    U_r & U_{2*}
\end{pmatrix}\Sigma_2 V_2^T\in\mathbb{R}^{n\times p_2},$$
where $U_r$ is the shared left singular subspace, $U_{1*}$ and $U_{2*}$ are the unshared left singular subspaces satisfying $U_{1*}^TU_{2*}=0$, and $(\Sigma_1, \Sigma_2)$ are diagonal matrices with decreasing singular values, then the above requirement can be implied by
\begin{align}\label{cond.s1}
\mathop{\max}\limits_{r+1 \leq k \leq \textup{rank}(X_1)}\sigma_{(k)}^2(X_1) \vee \mathop{\max}\limits_{r+1 \leq k \leq \textup{rank}(X_2)}\sigma_{(k)}^2(X_2)< c_1\mathop{\min}\limits_{1 \leq k \leq r}(\sigma_{(k)}^2(X_1)+\sigma_{(k)}^2(X_2)),
\end{align}
where $c_1\in(0,1)$ is some constant, $\sigma_{(k)}(X_i)$ is the singular value corresponding to the $k$-th singular vector of $X_i$. Our next theorem concerns the performance of Stack-SVD driven by the eigen-gap
\begin{align}\label{G}
 G_1(X)=\mathop{\min}\limits_{1 \leq k \leq r}\{\sigma_{(k)}^2(X_1)+\sigma_{(k)}^2(X_2)\}-\mathop{\max}\limits_{r+1 \leq k \leq \text{rank}(X_1)}\sigma_{(k)}^2(X_1) \vee \mathop{\max}\limits_{r+1 \leq k \leq \text{rank}(X_2)}\sigma_{(k)}^2(X_2).
\end{align}

\begin{theorem}\label{generalize}
    Suppose $Y_i\in \mathbb{R}^{n\times p_i},i=1,2,$ are generated from (\ref{model}), where the noise matrices $Z_i$ are $i.i.d.$ generated from $\mathscr{G}_\tau$. Define the parameter space
    \begin{align*}
\mathscr{H}^{(1)}_{r,t} = \left\{\begin{array}{c}X=(X_1~X_2):
       \textup{rank}(X_i)=r+r_{i*}, \\X_i = (U_r~U_{i*}) \Sigma_i V_i^T\in\mathbb{R}^{n\times p_i}, i=1,2, \\
    \text{(\ref{cond.s1}) holds},\quad U_{1*}^T U_{2*}=0, \quad G_1(X)\geq t^2
\end{array} \right\},
    \end{align*}
    where $\Sigma_i$ contains decreasing singular values and $G_1(X)$ is defined in (\ref{G}). Denote $\hat{U}_r$ as the Stack-SVD estimator, whose columns are the first $r$ left singular vectors of $Y=(Y_1~Y_2)\in\mathbb{R}^{n\times (p_1+p_2)}$. Then Equations (\ref{up1}) and (\ref{up2}) hold with $\mathscr{H}_{r,t}$ replaced by $\mathscr{H}^{(1)}_{r,t}$, and $\hat U_r^{\mathbb{J}}$ replaced by $\hat U_r$.
    If we further assume $r+r_{1*}+r_{2*} \leq \frac{n}{2}$, $r+r_{1*}\leq \frac{p_1}{16}$, $r+r_{2*}\leq \frac{p_2}{16}$, and either $p_1\asymp p_2$ or $\gamma^2\gtrsim p_1+p_2$ holds, then the minimax lower bounds (\ref{low1}) and (\ref{low2}) holds, with $\mathscr{H}_{r,t}$ replaced by $\mathscr{H}^{(1)}_{r,t}$.
\end{theorem}

Theorem \ref{generalize} indicates that the performance of top $r$ Stack-SVD is guaranteed and in fact minimax rate-optimal in the presence of weak unshared signals, provided that the shared signals are significantly more prominent than the unshared ones. 
Importantly,  for Stack-SVD, the conditions outlined in Theorem \ref{generalize} can be relaxed: we do not need to require that the singular values of $X_1$ and $X_2$, associated to $U_r$, are the largest ones for each matrix. By Propositions \ref{vectorrelation} and \ref{generalcainzhang},  as long as the singular values associated to $U_r$ are the largest ones in the stacked matrix $(X_1~X_2)$ and there is a sufficient eigen-gap between the  $r$-th and $(r+1)$-th singular values, the Stack-SVD estimator will still lead to consistent estimation.

\subsubsection{Scenario II: estimation with strong unshared signals}

We also consider an opposite scenario where the singular values of unshared vectors are greater than those of the shared vectors. Standard top $r$ Stack-SVD fails in this setting, but optimal estimation is still possible by correctly selecting vectors from the singular subspace of the stacked noisy matrix $(Y_1~Y_2)$. In this part,  we require that after stacking the singular values of the shared vectors are all smaller than the singular values of the unshared vectors. Specifically, we consider that $$\min_{j\notin \mathbb{J}} \sigma^2_{j}(X)\geq (1+\epsilon)\max_{j\in \mathbb{J}} \sigma^2_{j}(X)$$, for some small $\epsilon>0$. Similar to the previous case, we denote the individual SVDs of $X_1$ and $X_2$ as
$$
X_1=\begin{pmatrix}
    U_{1*} & U_{r}
\end{pmatrix}\Sigma_1 V_1^T
,\qquad
X_2=\begin{pmatrix}
    U_{2*} & U_{r}
\end{pmatrix}\Sigma_2 V_2^T,
$$
where $U_{1*}^TU_{2*}=0$, and $\Sigma_1$ and $\Sigma_2$ are diagonal matrices containing decreasing singular values. In particular, if we denote the number of unshared vectors in $X_1$ and $X_2$ as $d_1$ and $d_2$, then the above condition can be satisfied if 
\begin{equation}\label{cond.s}
\mathop{\max}\limits_{1 \leq k \leq r}\{\sigma_{(d_1+k)}^2(X_1)+\sigma_{(d_2+k)}^2(X_2)\} < c_1( \sigma_{d_1}^2(X_1) \wedge \sigma_{d_2}^2(X_2)),
\end{equation}
where $c_1\in(0,1)$ is some constant, $\sigma_{(k)}(X_i)$ is the  singular value corresponding to the $k$-th singular vector of $X_i$. 
Our next theorem concerns the minimax optimality of an alternative estimator $\hat U_r^{S}$, which consists of the left singular vectors associated with the $(d+1)$-th to $(d+r)$-th largest singular values of stacked matrix $(Y_1~Y_2)$.

\begin{theorem}\label{strongcontamination}
    Suppose $Y_i\in \mathbb{R}^{n\times p_i}, i=1,2,$ are generated from (\ref{model}), where the noise matrix $Z_i$ are i.i.d. generated from $\mathscr{G}_\tau$. Define the parameter space
    \begin{align*}
\mathscr{H}^{(2)}_{r,t} = \left\{\begin{array}{c}X=(X_1~X_2):
       \textup{rank}(X_i)=r+r_{i*}, \\X_i = (U_{i*}~U_r) \Sigma_i V_i^T\in\mathbb{R}^{n\times p_i}, i=1,2,\\
      \text{(\ref{cond.s}) holds},\quad U_{1*}^TU_{2*}=0,\quad  G_2(X)\geq t^2
\end{array} \right\},
    \end{align*}
    where $\Sigma_i$ contains decreasing singular values and $G_2(X)$ is given by
    $$
    G_2(X)=\sigma_{d_1}^2(X_1) \wedge \sigma_{d_2}^2(X_2)-\mathop{\max}\limits_{1 \leq k \leq r}\{\sigma_{(d_1+k)}^2(X_1)+\sigma_{(d_2+k)}^2(X_2)\}.$$
    Then Equations (\ref{up1}) and (\ref{up2}) hold with $\mathscr{H}_{r,t}$ replaced by $\mathscr{H}^{(2)}_{r,t}$, and $\hat U_r^{\mathbb{J}}$ replaced by $\hat U_r^S$.
    If we further assume $r+r_{1*}+r_{2*} \leq \frac{n}{2}$, $r+r_{1*}\leq \frac{p_1}{16}$, $r+r_{2*}\leq \frac{p_2}{16}$, and either $p_1\asymp p_2$ or $\gamma^2\gtrsim p_1+p_2$ holds, then the minimax lower bounds (\ref{low1}) and (\ref{low2}) holds, with $\mathscr{H}_{r,t}$ replaced by $\mathscr{H}^{(2)}_{r,t}$.
\end{theorem}

\subsection{Estimation with non-orthogonal unshared subspaces}\label{sec.nonorth}

So far we have been focusing on the cases with orthogonal unshared singular subspaces.
When the unshared singular subspaces are not orthogonal, we show that it can be transformed into an orthogonal scenario by studying the SVD of the stacked signal matrix $(X_1~X_2)$. In this case, the shared components remain unchanged, whereas the unshared subspaces will undergo rotation by a certain angle to achieve orthogonality. Generally, this rotation does not impact the estimation of the shared component. The linear transformations resulting from the non-orthogonality of the unshared subspaces only influence the unshared components, without inducing any interaction between shared and unshared components. Therefore, as long as similar gap conditions are still satisfied, the non-orthogonality doesn't compromise the optimal estimation of the shared subspace.

Our key insight comes from the following result, which concerns the SVD of the stacked signal matrix $(X_1~X_2)$ when the unshared subspaces are not orthogonal to each other.

\begin{theorem}\label{nonorthosvd}
   Let
    $$
    X_1=(
    U_r~ U_{1*})\mathrm{diag}(
    \Sigma_1,~\Sigma_{1*})V_1^T,\qquad
X_2=(
    U_r ~ U_{2*})\mathrm{diag}(
    \Sigma_2,\Sigma_{2*})V_2^T,$$
    be the SVDs of $X_1$ and $X_2$, where some left singular vectors in the unshared subspaces $U_{1*}\in \mathbb{O}(n,r_{1*})$ and $U_{2*}\in \mathbb{O}(n,r_{2*})$ are not orthogonal to each other. Suppose we have the eigen-decomposition
    $
    (U_{1*}~U_{2*})^T(U_{1*}~U_{2*})=\Gamma^T \Sigma \Gamma
    $. Then the following hold:
    \begin{enumerate}
   \item If $(U_{1*}~U_{2*})$ is not a column singular matrix,
    then the SVD of $(X_1~X_2)$ is
    $$
        (X_1~X_2)=(
        U_r~  (U_{1*}~U_{2*})SU^*)\mathrm{diag}(
        \sqrt{\Sigma_1^2+\Sigma_2^2},~\Sigma^*
    )\mathrm{diag}(
        I,~V^{*T})\begin{pmatrix}
        \frac{\Sigma_{1}}{\sqrt{\Sigma_{1}^2+\Sigma_{2}^2}}V_{1}^T & \frac{\Sigma_{2}}{\sqrt{\Sigma_{1}^2+\Sigma_{2}^2}}V_{2}^T\\
        V_{1*}^T & 0 \\ 0 & V_{2*}^T
    \end{pmatrix},$$
    where $S=\Gamma^T \Sigma^{-\frac{1}{2}} \Gamma $, and $U^*,\Sigma^*,V^*$ are defined from the SVD
    $
    S^{-1}\mathrm{diag}(
        \Sigma_{1*},~\Sigma_{2*})=U^*\Sigma^*V^{*T}
    $.
    \item If $(U_{1*}~U_{2*})$ is a column singular matrix and $rank(U_{1*}~U_{2*})=r^*$, there exist a transformation $L\in \mathbb{R}^{(r_{1*}+r_{2*})\times(r_{1*}+r_{2*})}$ such that $(U_{1*}~U_{2*})L=(\Tilde{U}^*~{\bf 0}_{n\times (r_{1*}+r_{2*}-r^*)})$, where $\Tilde{U}^*$ is a column normalized matrix. Similarly define $S=\Gamma^T \Sigma^{-\frac{1}{2}} \Gamma $ based on $\Tilde{U}^{*T}\Tilde{U}^{*}=\Gamma^T \Sigma \Gamma$. Then the SVD of $(X_1~X_2)$ is   \begin{align*}
        (X_1~X_2)=(
        U_r~  \Tilde{U}^*SU^* )\mathrm{diag}(
        \sqrt{\Sigma_1^2+\Sigma_2^2},~\Sigma^*)\mathrm{diag}(
        I,~V^{*T})\begin{pmatrix}
        \frac{\Sigma_{1}}{\sqrt{\Sigma_{1}^2+\Sigma_{2}^2}}V_{1}^T & \frac{\Sigma_{2}}{\sqrt{\Sigma_{1}^2+\Sigma_{2}^2}}V_{2}^T\\
        V_{1*}^T & 0 \\ 0 & V_{2*}^T
    \end{pmatrix},
    \end{align*}
    where $U^*,\Sigma^*,V^*$ come from the SVD
    {\small $
    (S^{-1}~{\bf 0}_{r^*\times(r_{1*}+r_{2*}-r^*)})L^{-1}\mathrm{diag}(
              \Sigma_{1*},~\Sigma_{2*})=U^*\Sigma^*V^{*T}.
    $}
    \end{enumerate}
\end{theorem}

Theorem \ref{nonorthosvd} demonstrates that, in the SVD of $(X_1~X_2)$, compared to that obtained in the orthogonal case, while there may be rotations within the subspace related to the unshared singular vectors, the shared subspace and their corresponding singular values remain unaffected by such non-orthogonality. Without merging the shared and unshared signals, the stacking preserves the structure of the shared subspace in its SVD, making it feasible to estimate the shared  subspace. 

In fact, optimal estimation is achievable as long as the gap conditions are still met. Recall $\mathbb{S}$ defined in Section \ref{sec.3.1}. We consider the parameter space
\begin{align*}
\mathscr{S}_{r,t} = \left\{\begin{array}{c}X=(X_1~X_2):\\
       \text{rank}(X_i)=r+r_{i*}, X_i = (U_r~U_{i*}) \Sigma_i V_i^T\in\mathbb{R}^{n\times p_i}\\
     g_{i(s)}^2(X)>c\sigma_{s+1}^2(X), \forall s\in\mathbb{S},\\ \mathop{\min}\limits_{1 \leq k \leq N}g_k^2(X) \geq t^2, i=1,2
\end{array} \right\},
    \end{align*}
where $U_{i*} \in \mathbb{O}(n,r_{i*}), i=1,2$, are not necessarily orthogonal to each other.
In particular, given the different forms of SVD, the eigen-gap condition $g_{i(s)}^2(X)>c\sigma_{s+1}^2(X)$ may also vary from case to case. Specifically, based on Theorem \ref{nonorthosvd}, for the shared subspace, the corresponding singular values are still in the diagonals of $\sqrt{\Sigma_1^2+\Sigma_2^2}$; however, the singular values for the unshared vectors in this case are the singular values of the matrix $
    S^{-1}\mathrm{diag}(
        \Sigma_{1*},~\Sigma_{2*})$ when $(U_{1*}~U_{2*})$ is column non-singular and $(S^{-1}~{\bf 0}_{r^*\times(r_{1*}+r_{2*}-r^*)})L^{-1}\mathrm{diag}(
              \Sigma_{1*},~\Sigma_{2*})$ when $(U_{1*}~U_{2*})$ is not column non-singular. 
          
          To achieve minimax optimal estimation, we still consider the oracle estimator $\hat U_r^{\mathbb{J}}$ in which the locations of the shared singular vectors are correctly identified. The following theorem establishes its minimax optimality.

\begin{theorem}\label{nonorthobounds}
 Suppose $Y_i,i=1,2,$ are generated from (\ref{model}), where the noise matrices $Z_i$ are $i.i.d.$ generated from $\mathscr{G}_\tau$. For any finite $r$, the upper bounds (\ref{up1}) and (\ref{up2}) in Theorem \ref{generalbounds} hold, with $\mathscr{H}_{r,t}$ replaced by $\mathscr{S}_{r,t}$.
If we further assume $r+r_{1*}+r_{2*} \leq \frac{n}{2}$, $r+r_{1*}\leq \frac{p_1}{16}$, $r+r_{2*}\leq \frac{p_2}{16}$, and either $p_1\asymp p_2$ or $\gamma^2\gtrsim p_1+p_2$ holds, then the minimax lower bounds (\ref{low1}) and (\ref{low2}) in Theorem \ref{generalbounds} still hold, with $\mathscr{H}_{r,t}$ replaced by $\mathscr{S}_{r,t}$.
\end{theorem}

Theorem \ref{nonorthobounds} demonstrates that the presence of non-orthogonal unshared singular subspace will not significantly alter the landscape regarding the optimal estimation of the shared singular subspace. Here follows an analysis of a detailed non-orthogonal case to clarify the underlying principles and make the concepts more accessible. We consider a rank-$2$ example, where $X_1$ and $X_2$ are both rank-$2$ matrices with one shared left singular vector ($u$) and one unshared vector ($u_1$ or $u_2$), 
$$X_1=(
    u ~ u_1)\mathrm{diag}(
    \alpha,~\sigma_a)V_1^T
, ~
X_2=(
    u ,~u_2)\mathrm{diag}(
    \beta,~\sigma_b)V_2^T.
$$ Assume $u_1$ and $u_2$ are not orthogonal to each other, that is, $u_1^Tu_2 \neq 0$.
We use the first left singular vector of the stacked noisy $(Y_1~Y_2)$ as the estimate of u.

In constructing the transformation, we consider 
$$
    R=\begin{pmatrix}
        1&\\&S
    \end{pmatrix}=\begin{pmatrix}
        1 & 0 & 0\\
        0 & 1 & -\frac{u_1^Tu_2}{\sqrt{1-(u_1^Tu_2)^2}}\\
        0 & 0 & \frac{1}{\sqrt{1-(u_1^Tu_2)^2}}
    \end{pmatrix}, \qquad
    R^{-1}=\begin{pmatrix}
        1&\\&S^{-1}
    \end{pmatrix}=\begin{pmatrix}
         1 & 0 & 0\\
        0 & 1 & u_1^Tu_2\\
        0 & 0 & \sqrt{1-(u_1^Tu_2)^2}
    \end{pmatrix}.$$
We denote the SVD
$$
T=\begin{pmatrix}
    \sigma_a & u_1^Tu_2 \sigma_b \\ 0 & \sqrt{1-(u_1^Tu_2)^2}\sigma_b
\end{pmatrix}=\begin{pmatrix}
    u_{11} & u_{12}\\ u_{21}& u_{22}
\end{pmatrix}\begin{pmatrix}
    \sigma_1(T) & 0 \\
    0 & \sigma_2(T)
\end{pmatrix}V_s^T.
$$
After applying this transformation, we have the SVD of the stacked noiseless matrix $(X_1~X_2)$:
$$
    \begin{pmatrix}
    u & u_1u_{11}-\frac{u_{21}((u_1^Tu_2) u_1-u_2)}{\sqrt{1-(u_1^Tu_2)^2}} & u_1u_{12}-\frac{u_{22}((u_1^Tu_2) u_1-u_2)}{\sqrt{1-(u_1^Tu_2)^2}}
\end{pmatrix}\\
\cdot\mathrm{diag}(
    \sqrt{\alpha^2+\beta^2},~\sigma_1(T),~\sigma_2(T))V^T.$$
Therefore, by calculation of $\sigma_1(T)$ and $\sigma_2(T)$ combined with our theories, we have that, if the signal gap constraint is satisfied by setting unshared signals to be weaker than shared signals, i.g. there exists constants $0<c_1<1$ s.t.$$
    \frac{\sigma_a^2+\sigma_b^2+\sqrt{(\sigma_a^2+\sigma_b^2)^2-4(1-|u_1^Tu_2|^2)\sigma_a^2\sigma_b^2}}{2} < c_1\sqrt{\alpha^2+\beta^2}
    $$
    for any $\gamma \leq \sqrt{\alpha^2+\beta^2}$, when estimating with the first left singular vector of $(Y_1~Y_2)$, it holds that $$
    \mathbb{E}\|\sin\Theta(\hat{u},u)\|^2 \leq \frac{cn(\gamma^2+(p_1+p_2))}{\gamma^4}\wedge 1
    .$$
Conversely, when the signal strength of unshared vector is stronger, the gap condition will become
$$
    \sqrt{\alpha^2+\beta^2}<c_1\frac{\sigma_a^2+\sigma_b^2-\sqrt{(\sigma_a^2+\sigma_b^2)^2-4(1-|u_1^Tu_2|^2)\sigma_a^2\sigma_b^2}}{2},
    $$
where the same bound holds for estimating with the second singular vector of $(Y_1~Y_2)$.

The results presented above indicate that, even when the unshared singular vectors are not  orthogonal, the estimation can still attain rate-optimality by appropriately selecting the singular vectors, as long as the singular value gap conditions are satisfied. The required gap conditions are influenced by the degree of non-orthogonality present. Furthermore, this non-orthogonality affects only the unshared component, without impacting the shared component or leading to the fusion of information between the shared and unshared portions.

\section{Tracing Shared Singular Vectors}\label{algorithm}

As suggested by the previous discussions, achieving optimal estimation requires accurately identifying the shared and unshared singular vectors from the stacked matrix. In this section, we propose an algorithm designed to distinguish the shared and unshared singular vectors across multiple matrices. This is accomplished by carefully comparing the singular subspaces of the stacked matrix with those of the individual matrices. Specifically, we assume that the unshared singular vectors are mutually orthogonal and that the singular values in the signal matrices are sufficiently well-separated. 

We first present the algorithm for the case of two matrices and then extend it to the general case involving $k$ matrices. Suppose $X_1$ and $X_2$ have ranks $r_1$ and $r_2$, respectively, with $r<\min(r_1,r_2)$ shared left singular vectors. Denote the number of unshared left singular vectors as $k_1 = r_1 - r$ and $k_2 = r_2 - r$ for $X_1$ and $X_2$, respectively. Let $\hat{U}_1$ and $\hat{U}_2$ represent the top $r_1$ and $r_2$ left singular vectors of $Y_1$ and $Y_2$, with the $i$-th and $j$-th singular vectors denoted as $\hat{u}_{1i}$ and $\hat{u}_{2j}$. The top $r + k_1 + k_2$ left singular subspace of the stacked matrix $(Y_1~Y_2)$ is denoted as $\hat{U}$, with the $i$-th vector as $\hat{u}_i$. The noiseless counterparts of these singular vectors and subspaces can be defined analogously, without the hats. For the moment, we assume $k_1$ and $k_2$ are known while presenting the following Algorithm \ref{alg:tracing}. Later, we will introduce a method to estimate them. 

\begin{algorithm}[!h]
    \caption{Identifying Shared Singular Vectors}
    \label{alg:tracing}
    \renewcommand{\algorithmicrequire}{\textbf{Input:}}
    \renewcommand{\algorithmicensure}{\textbf{Output:}}
    \begin{algorithmic}[1]
        \REQUIRE $Y_1$, $Y_2$, $k_1$, $k_2$, $r$  
        \ENSURE $\widehat{\mathbb{J}}$    
        
        \STATE  Compute top $r_1$, $r_2$, $r+k_1+k_2$ left singular vectors, denoted as $\hat{U}_1$, $\hat{U}_2$, $\hat{U}$, for $Y_1$, $Y_2$ and $(Y_1~Y_2)$ respectively.
        
        \FOR{each $i \in 1,...,r_1$}
            \STATE $d_{1i}=\mathop{\min}\limits_{j \in [1,r_2]}\|\sin\Theta(\hat{u}_{1i},\hat{u}_{2j})\|^2$
        \ENDFOR

        \FOR{each $j \in 1,...,r_2$}
            \STATE $d_{2j}=\mathop{\min}\limits_{i \in [1,r_1]}\|\sin\Theta(\hat{u}_{2j},\hat{u}_{1i})\|^2$
        \ENDFOR

        \STATE Take the last $r_1-k_1$ and $r_2-k_2$ index sets, $I$ and $J$, for $d_{1i}$ and $d_{2j}$ respectively. $I=\{i_1,...,i_{r_1-k_1}\}$, $J=\{j_1,...,j_{r_2-k_2}\}$, satisfying that $d_{1i_k}$ is the $k$-th smallest value among $d_{1i}$ while $d_{2j_k}$ is the $k$-th smallest value among $d_{2j}$

        \STATE Compute index sets:\\
        $K_1=\{k: k=\mathop{\arg\min}\limits_{k \in 1,...,r+k_1+k_2}\|\sin\Theta(\hat{u}_{1i},\hat{u}_k)\|^2, i \in I\}$ \\
        $K_2=\{k: k=\mathop{\arg\min}\limits_{k \in 1,...,r+k_1+k_2}\|\sin\Theta(\hat{u}_{2j},\hat{u}_k)\|^2, j \in J\}$ \\
        $\widehat{\mathbb{J}}=K_1 \cup K_2$

        \RETURN $ \widehat{\mathbb{J}}$
    \end{algorithmic}
\end{algorithm}

The idea behind Algorithm \ref{alg:tracing} is simple: based on the assumption that the individual SVDs of $Y_1$ and $Y_2$ and the SVD of $(Y_1~Y_2)$ provide sufficiently accurate estimations of their respective individual singular vectors under some conditions, we can leverage this  to identify shared and unshared vectors. Specifically, if a vector $u$ in the  subspace of $Y_1$ is shared, there will be a corresponding vector in the singular subspace of $Y_2$ that is  close to $u$. Conversely, if a vector $u$ in the  subspace of $Y_1$ is unshared, by the orthogonality assumption, all vectors in the  subspace of $Y_2$ will be largely different from $u$.
The following theorem concerns the consistency of this approach, proved under the assumption that the unshared vectors are mutually orthogonal and the singular values in both signal matrices are sufficiently large and separated. 

\begin{theorem}\label{algothrm}
    Suppose the unshared singular vectors are mutually orthogonal and there exists a small constant $\epsilon_1>0$ such that $\sigma_{k}\ge (1+\epsilon_1)\sigma_{k+1}$ for all singular values in $(X_1~X_2)$. We denote $\mathbb{J}$ as the index set for the shared singular vectors in the singular subspace of $(X_1~X_2)$. Suppose $(\alpha^{c_1} \wedge \beta^{c_1})\ge Cn$, for some constant $c_1<2$ and sufficiently large $C>0$, and $(\alpha \wedge \beta)^2\ge C_2(p_1 \vee p_2)$ for some $C_2>0$, where $\alpha=\sigma_{\min}(X_1)$ and $\beta=\sigma_{\min}(X_2)$. Then for sufficiently large $(n,p_1,p_2)$, there exists some constant $0<\epsilon<2-c_1$, such that
    $$
    \mathbb{P}\{\widehat{\mathbb{J}}=\mathbb{J}\}\geq 1-o({\alpha^{-\epsilon/2}})-o({\beta^{-\epsilon/2}})-o({\gamma^{-\epsilon/2}}),$$ where $\gamma=\alpha \wedge \beta$.
\end{theorem}

It is important to emphasize that by requiring sufficiently large signal strength, we implicitly ensure that the singular values are well separated and each singular vector is sufficiently identifiable. This is because the minimum gap between singular values and the minimum signal strength are inherently in the same order according to previous analysis.

Our analyses reveal that for unshared vectors in one individual matrix, the distances between these vectors and any left singular vectors in the other individual matrix tend to be relatively large—approaching \(1\) as the signal strength increases. In contrast, for shared vectors in one individual matrix, there exists a corresponding vector in the other individual matrix with a relatively small distance—approaching \(0\) as the signal strength becomes sufficiently large. This observation enables us to determine the number of unshared vectors, denoted as \(k_1\) and \(k_2\), by analyzing the distances \(d_{1i}\) and \(d_{2j}\). Algorithm \ref{alg:k1k2} outlines a method for estimating \(k_1\) and \(k_2\) by identifying the largest gap in these distance measurements.

\begin{algorithm}[!h]
    \caption{Estimating $k_1$ and $k_2$}
    \label{alg:k1k2}
    \renewcommand{\algorithmicrequire}{\textbf{Input:}}
    \renewcommand{\algorithmicensure}{\textbf{Output:}}
    \begin{algorithmic}[1]
        \REQUIRE 
        $Y_1$, $Y_2$;   
        \ENSURE $\hat{k}_1$ $\hat{k}_2$    
        \STATE Similarly calculate $d_{1i},~i=1,...,r_1$ and $d_{2j},~j=1,...,r_2$ as in Algorithm \ref{alg:tracing}.
        \STATE  Sort $d_{1i},~i=1,...,r_1$ and $d_{2j},~j=1,...,r_2$ in decreasing order respectively. The ordered sequences are denoted as $d_{1(i)}$ and $d_{2(j)}$, with values in the parenthesis representing its order.

        \STATE Set $g^{(1)}_0=1-d_{1(1)}$, $g^{(1)}_{r_1}=d_{1(r_1)}-0$; $g^{(2)}_0=1-d_{2(1)}$, $g^{(2)}_{r_2}=d_{2(r_2)}-0$
        
        \FOR{each $i \in 1,...,r_1-1$}
            \STATE $g^{(1)}_i=d_{1(i)}-d_{1(i+1)}$
        \ENDFOR

        \FOR{each $j \in 1,...,r_2-1$}
            \STATE $g^{(2)}_j=d_{2(j)}-d_{2(j+1)}$
        \ENDFOR

        \STATE Compute: $\hat k_1=\mathop{\arg\max}\limits_{k \in 0,...,r_1}g^{(1)}_k$, \quad $\hat k_2=\mathop{\arg\max}\limits_{k \in 0,...,r_2}g^{(2)}_k$

        \RETURN $\hat k_1$ $\hat k_2$
    \end{algorithmic}
\end{algorithm}

The consistency of the estimators $\hat k_1$ and $\hat k_2$, obtained in Algorithm \ref{alg:k1k2}, is established in Theorem \ref{findkalg} below. Together with Theorem \ref{algothrm}, these results provide a rigorous theoretical justification for our proposed algorithms.

\begin{theorem}\label{findkalg}
    Under the assumption of Theorem \ref{algothrm}, for sufficiently large $(n,p_1,p_2)$, we have $$P(\hat k_1=k_1, \hat k_2=k_2)\ge 1-o({\alpha^{-\epsilon/2}})-o({\beta^{-\epsilon/2}})-o({\gamma^{-\epsilon/2}}).$$
\end{theorem}

Note that Algorithms \ref{alg:tracing} and \ref{alg:k1k2} are designed  to help identify the shared singular vectors. To identify the unshared singular vectors, we can either look at the complement index set $$\widehat{\mathbb{J}}^c=\{1,2,...,r+k_1+k_2\}\setminus \widehat{\mathbb{J}},$$ or modify the step $8$ in Algorithm \ref{alg:tracing} by selecting the top  $k_1$ and $k_2$ index sets and tracing them in the stacked matrix. For the estimation of $\mathbb{J}$ in the presence of non-orthogonal unshared singular vectors, see Section \ref{discussion} for more discussion.

In Algorithm \ref{alg:tracing}, either $K_1$ or $K_2$ alone is sufficient to provide a consistent estimate of $\hat{J}$ under our theoretical assumptions. The algorithm employs their union ($K_1 \cup K_2$) for practical robustness, increasing the likelihood of capturing the complete shared subspace, even at the cost of including a few unshared singular vectors. Both Algorithm \ref{alg:tracing} and Algorithm \ref{alg:k1k2} rely on the assumption that the unshared components between matrices are perfectly orthogonal. In practice, however, this assumption may not hold exactly—non-orthogonal unshared components can exhibit partial alignment, implying that vectors labeled as “unshared” may still contain a non-negligible amount of common information. As a result, the algorithm may adaptively include such vectors in the estimated shared subspace, resulting $K_1\ne K_2$. This behavior does not indicate a failure of the algorithm; rather, it reflects its ability to recover empirically shared signals that extend beyond the theoretical definition of exact sharing.

Finally, there are several ways to extend the above methods to the case involving  $k $ matrices. One way is to select the pairwise shared singular vectors and then integrate the results. Another approach for generalizing the algorithm to trace shared singular vectors across  $k$  noisy matrices is to redefine the distance used in Algorithm \ref{alg:tracing} accordingly as
\begin{align*}d_{si}=\max &\{\mathop{\min}\limits_{j \in [1,r_1]}\|\sin\Theta(\hat{u}_{si},\hat{u}_{1j})\|^2,...,\mathop{\min}\limits_{j \in [1,r_{s-1}]}\|\sin\Theta(\hat{u}_{si},\hat{u}_{(s-1)j})\|^2,\\ &\qquad\qquad\mathop{\min}\limits_{j \in [1,r_{s+1}]}\|\sin\Theta(\hat{u}_{si},\hat{u}_{(s+1)j})\|^2,...,\mathop{\min}\limits_{j \in [1,r_k]}\|\sin\Theta(\hat{u}_{si},\hat{u}_{kj})\|^2\},\end{align*}
for any $s=1,2,...,k$. The implication is straightforward. If a vector is identified as a shared vector, then for any other matrix, there must exist a vector that is close to it. The subsequent steps in Algorithms \ref{alg:tracing} and \ref{alg:k1k2} can be correspondingly generalized to handle $k$ matrices with ease. Under conditions analogous to those in Theorem \ref{algorithm}, but in the context of $k$ matrices, the performance of the generalized algorithms can also be guaranteed by similar results.

\section{Simulation Studies}\label{simulation}


\subsection{Assessing numerical performance of Stack-SVD}
We first demonstrate the empirical advantages of Stack-SVD when the singular subspace is completely shared. We  randomly generate three $n$-dimensional orthonormal vectors to serve as the shared left singular vectors. Then we generate random orthonormal matrices of dimensions $3 \times p_1$ and $3 \times p_2$ for the right singular subspaces of $X_1$ and $X_2$, respectively. We set the corresponding singular values as $$\Sigma_1=\operatorname{diag}\{\alpha, \frac{\alpha}{2}, \frac{\alpha}{4}\},\qquad \Sigma_2=\operatorname{diag}\{\beta, \frac{\beta}{2}, \frac{\beta}{4}\}.$$
The noise matrices $Z_1$ and $Z_2$ are generated from iid $N(0,1)$ and we run $500$ simulations for each setting. In each case, the empirical estimation error is reported as $\|\sin\Theta(U,\hat{U})\|^2$. Both the mean and standard deviation of such distances are recorded.

We set $n=10$ and $n=20$ respectively, and also divide the scenarios into balanced and imbalanced dimension groups, with $(p_1, p_2)$ being $(20, 20)$ and $(300,400)$ for balanced dimension group, and $(20, 300)$ and $(50,500)$ for imbalanced dimension group. For the signal strength, while the overall signal strength $\sqrt{\alpha^2 + \beta^2}$ remains approximately constant, the individual signal strengths $(\alpha, \beta)$ transition from balanced values $(50, 50)$ to progressively imbalanced values $(37, 60)$ and $(10, 70)$ in relatively low-dimensional settings and from balanced values $(100, 100)$ to progressively imbalanced values $(74, 120)$ and $(20, 140)$ in relatively high-dimensional settings. Also, for the imbalanced dimension settings, we reverse the above signal strengths in $X_1$ and $X_2$ to evaluate the performance. We compare the performance of the following estimators: Individual SVD ("SVD-$X_1$" and "SVD-$X_2$"), Stack-SVD and Average-SVD. The results are summarized in Table \ref{simuexample1table}. 

Several observations can be made from Table \ref{simuexample1table}. First, by maintaining the overall signal strength $\sqrt{\alpha^2 + \beta^2}$ at an approximately constant level, the estimation error for Stack-SVD remains  small. Second, Stack-SVD exhibits the greatest stability, generally presenting the smallest standard deviation among the estimation methods. Third, Stack-SVD generally outperforms individual SVD in most cases. However, there are scenarios where an individual SVD may yield a smaller estimation error. This outcome is primarily driven by highly imbalanced signal strengths across matrices, where a strong performance by one individual-SVD estimator is often accompanied by poor estimation from the estimator of another matrix; see, for example, the last row of Table \ref{simuexample1table}. In this case, even if one of the individual-SVD estimators is precise, we should note that in practice it requires the users  to determine which individual-SVD estimator performs better. In comparison, Stack-SVD provides a simple solution with competitive performance. Lastly, while Average-SVD performs adequately under balanced signal strength, their estimates become less accurate when the signal strength is imbalanced.

\begin{table}[]
\centering
\caption{Average estimation errors ($\|\sin\Theta(U,\hat{U})\|^2$) and their standard deviations based on $500$ simulations. The smallest estimation error in each row is highlighted in bold.}
\label{simuexample1table}
\scalebox{1}{
\begin{tabular}{c|c|cccc}
$(n,p_1,p_2)$                 & $(\alpha,\beta)$ & SVD-$X_1$    & SVD-$X_2$    & Stack-SVD    & \multicolumn{1}{l}{Average-SVD} \\ \hline \hline
\multirow{3}{*}{(10,20,20)}   & (50,50)          & 0.049(0.024) & 0.050(0.025) & \textbf{0.027(0.012)} & \textbf{0.027(0.012)}                    \\
                              & (37,60)          & 0.095(0.047) & 0.034(0.017) & \textbf{0.027(0.012)} & 0.035(0.016)                    \\
                              & (10,70)          & 0.806(0.191) & \textbf{0.025(0.012)} & 0.026(0.012) & 0.321(0.116)                    \\ \hline
\multirow{3}{*}{(20,300,400)} & (100,100)        & 0.039(0.013) & 0.045(0.016) & \textbf{0.021(0.008)} & \textbf{0.021(0.008)}                    \\
                              & (74,120)         & 0.088(0.031) & 0.028(0.010) & \textbf{0.021(0.008)} & 0.030(0.011)                    \\
                              & (20,140)         & 0.911(0.104) & \textbf{0.019(0.006)} & 0.021(0.008) & 0.388(0.090)                    \\ \hline
\multirow{5}{*}{(10,20,300)}  & (50,50)          & 0.051(0.025) & 0.146(0.090) & \textbf{0.050(0.025)} & 0.052(0.031)                    \\
                              & (37,60)          & 0.098(0.048) & 0.080(0.044) & 0.050(0.026) & \textbf{0.047(0.024)}                    \\
                              & (60,37)          & \textbf{0.035(0.017)} & 0.381(0.236) & 0.050(0.025) & 0.130(0.102)                    \\
                              & (10,70)          & 0.808(0.182) & \textbf{0.049(0.026)} & \textbf{0.049(0.026)} & 0.332(0.126)                    \\
                              & (70,10)          & \textbf{0.025(0.013)} & 0.914(0.098) & 0.049(0.025) & 0.388(0.088)                    \\ \hline
\multirow{5}{*}{(20,50,500)}  & (100,100)        & 0.030(0.010) & 0.050(0.015) & \textbf{0.020(0.007)} & \textbf{0.020(0.007)}                    \\
                              & (74,120)         & 0.056(0.019) & 0.030(0.009) & \textbf{0.020(0.007)} & 0.022(0.008)                    \\
                              & (120,74)         & \textbf{0.020(0.007)} & 0.125(0.043) & \textbf{0.020(0.007)} & 0.038(0.014)                    \\
                              & (20,140)         & 0.781(0.187) & \textbf{0.020(0.006)} & \textbf{0.020(0.006)} & 0.302(0.116)                    \\
                              & (140,20)         & \textbf{0.015(0.005)} & 0.927(0.090) & 0.020(0.007) & 0.398(0.084)                   
\end{tabular}}
\end{table}

\subsection{Estimating individually non-identifiable singular vectors}\label{simu2}
Our next experiment focuses on the interesting scenario where the shared singular vectors are not individually identifiable in each data matrix, whereas  Stack-SVD is still able to capture them accurately. This experiment is designed as
$$X_1=(
    u ~ u_1)\mathrm{diag}(
    \alpha,~\alpha)V_1^T ,\qquad X_2=(
    u ~ u_2)\mathrm{diag}(
    \beta,\beta)V_2^T.$$ The dimension settings include both low and high,  balanced and imbalanced dimensions. For signal strength $(\alpha, \beta)$, we consider various settings from $(10, 15)$, $(20,20)$ to $(50, 60)$. The estimation errors for both individual SVD and Stack-SVD are reported in Table \ref{weaktable} based on $500$ simulations.

\begin{table}[]
\caption{Average estimation errors ($\|\sin\Theta(U_r,\hat{U}_r)\|^2$) and their standard deviations based on $500$ simulations. The SVD on individual $X_i$ are denoted as SVD-$X_i$.}\label{weaktable}
\centering
\scalebox{1}{
\begin{tabular}{c|c|ccc}
$(n,p_1,p_2)$                   & $(\alpha,\beta)$ & SVD-$X_1$ & SVD-$X_2$ & Stack-SVD    \\ \hline \hline
\multirow{3}{*}{(10,100,100)}   & (10,15)          & 0.541(0.323) & 0.502(0.337) & 0.114(0.126) \\
                                & (20,20)          & 0.485(0.359) & 0.492(0.347) & 0.025(0.017) \\
                                & (50,60)          & 0.470(0.362) & 0.481(0.352) & 0.003(0.002) \\ \hline
\multirow{3}{*}{(20,1000,1000)} & (10,15)          & 0.856(0.144) & 0.683(0.255) & 0.553(0.241) \\
                                & (20,20)          & 0.574(0.311) & 0.580(0.326) & 0.113(0.063) \\
                                & (50,60)          & 0.493(0.358) & 0.503(0.367) & 0.006(0.003) \\ \hline
\multirow{3}{*}{(10,100,500)}   & (10,15)          & 0.558(0.312) & 0.554(0.318) & 0.222(0.221) \\
                                & (20,20)          & 0.493(0.351) & 0.527(0.332) & 0.040(0.026) \\
                                & (50,60)          & 0.485(0.358) & 0.506(0.346) & 0.004(0.003) \\ \hline
\multirow{3}{*}{(20,1000,2000)} & (10,15)          & 0.851(0.153) & 0.756(0.195) & 0.637(0.228) \\
                                & (20,20)          & 0.569(0.304) & 0.608(0.281) & 0.159(0.081) \\
                                & (50,60)          & 0.508(0.355) & 0.502(0.356) & 0.007(0.003)
\end{tabular}}
\end{table}

As expected, individual-SVD estimators are not able to provide reliable estimates, with most estimation errors above $0.5$. In contrast, Stack-SVD reliably estimates the shared vector, with estimation error significantly decreasing as the signal strength increases.

\subsection{Assessing the proposed signal-tracing algorithms}

Finally, we conduct a simulation study to evaluate the algorithm for tracing the shared and unshared singular vectors. Specifically, we examine three simulation settings with increasing complexity. We denote $U_1$ and $U_2$ as the left singular subspace corresponding to $X_1$ and $X_2$ respectively. The shared vectors are denoted as $u_i$ while the unshared vectors are denoted as $u_{i*}$.
\begin{itemize}
    \item Setting 1: 
    $U_1=(u~u_{1*}), U_2=(u~u_{2*}~u_{3*})$, with one shared and three unshared vectors.
    \item Setting 2: $U_1=(u_1~u_2~u_{1*}~u_{2*}), U_2=(u_1~u_2~u_{3*}~u_{4*})$, with two shared and four unshared vectors.
    \item Setting 3: $U_1=(u_1~u_2~u_3~u_{1*}~u_{2*}~u_{3*}), U_2=(u_1~u_2~u_3~u_{4*}~u_{5*}~u_{6*})$, with three shared  and six unshared vectors.
\end{itemize}
The unshared vectors are orthogonal to each other.
Upon stacking two matrices, the shared and unshared singular vectors become shuffled within the left singular subspace of the stacked signal matrix. Here we define the $G_i:=\sigma_i-\sigma_{i+1}$ as the singular value gaps in the noiseless stacked matrix $(X_1~X_2)$.  The minimum $\min_{1\le i\le N}G_i$ measures overall SNR. In this case, given our focus on the interactions between shared and unshared signals rather than distinguishing signals from noise like previous examples, we set the minimum singular value of the stacked signal matrix to be sufficiently large. A situation is classified as a success if all unshared vectors are accurately detected without any errors. To evaluate the effectiveness of our approach, we estimate the success rates across different settings and various signal strengths, based on $1,000$ simulation trials.

\begin{table}
\centering
\caption{Accuracy for The Algorithm under Different Settings Based on $1000$ Simulations}
\label{algorithmsimulation}
\scalebox{1}{
\begin{tabular}{cc||cc||cc}
$\mathop{\min}\limits_{1 \leq i \leq N}G_i$  & Setting 1 & $\mathop{\min}\limits_{1 \leq i \leq N}G_i$  & Setting 2 & $\mathop{\min}\limits_{1 \leq i \leq N}G_i$  & Setting 3 \\ \hline  \hline
1  & 0.622  & 1  & 0.187  & 1  & 0.089  \\ \hline
3  & 0.981  & 5  & 0.782  & 5  & 0.727  \\ \hline
5  & 1.000  & 10 & 0.815  & 10 & 0.798  \\ \hline
7  & 1.000  & 15 & 1.000  & 20 & 0.999  \\ \hline
10 & 1.000  & 20 & 1.000  & 25 & 1.000 
\end{tabular}}
\end{table}

As illustrated in Table \ref{algorithmsimulation}, the algorithm demonstrates reliable performance across all cases when the minimum singular value gap is sufficiently high. Additionally, when considering the complexity of the various scenarios, there is a noticeable trend that more complicated cases tend to necessitate greater minimum gap conditions to ensure accurate detection. 

\section{Application to Single-Cell Data Integration}\label{realdata}

\begin{figure}[h!]
    \centering
    \caption{Comparison of the Silhouette Index, Calinski-Harabasz Index, and Neighborhood Purity. The embeddings of data from 10x Chromium (v3), 10x Chromium (v2) and Drop-seq are labeled as Data 1, Data 2 and Data 3.}  
    \label{realdata2}
    \includegraphics[width=\textwidth]{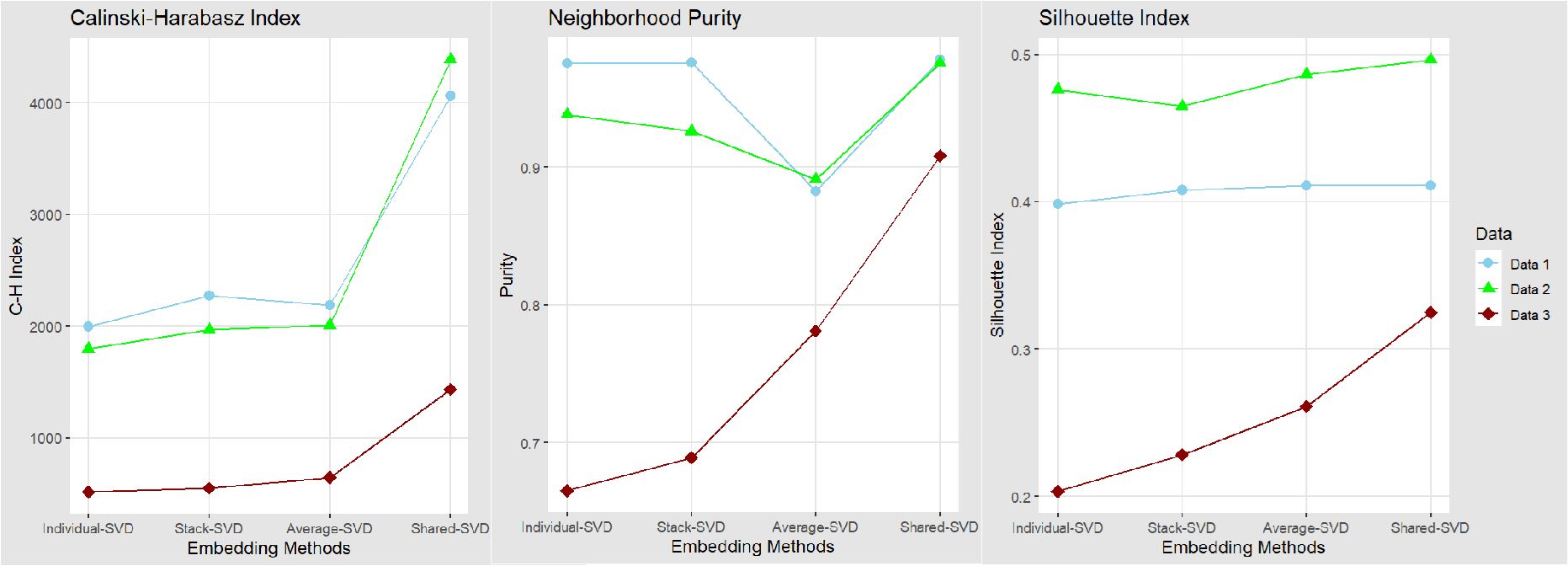} 
\end{figure}

We apply our proposed methods to jointly analyze three single-cell RNA-seq datasets of human peripheral blood mononuclear cells (PBMCs), each generated using a distinct sequencing technology \citep{pbmcsca}. Specifically, the datasets correspond to three technologies: 10x Chromium v3 (Dataset 1), 10x Chromium v2 (Dataset 2), and Drop-seq (Dataset 3). Each dataset consists of a gene expression matrix $Y_i\in\mathbb{R}^{g\times c_i}$ containing thousands of cells (whose numbers are denoted by $c_i$ for each dataset) and $g$ genes. For these datasets, since they are all related to PBMCs, we expect they contain the same family of cell types, governed by the similar genetic programs. In particular, we reason that the latent gene structures captured by the leading left singular vectors of each gene expression matrix is at least partially shared across all datasets.
Following standard preprocessing and feature selection pipelines using the R package \texttt{Seurat}, we obtain three gene-by-cell matrices, each with $g=1000$ genes and $c_1=3222$, $c_2=3362$, and $c_3=3500$ cells respectively. We then apply different methods to estimate the shared left singular subspaces ${U}_r$ across the three datasets that captures the latent gene structures. For each dataset, for each shared subspace estimator $\hat U_r$, we compute the joint low-dimensional cell embeddings as $Y_i^\top \hat{U}_r$ for $i = 1, 2, 3$.
We expect that more accurate estimation of the shared latent subspace ${U_r}$ will produce improved cell embeddings where different cell types are more distinctly clustered. We evaluate clustering quality using three metricsm the Silhouette index, Calinski-Harabasz (CH) index, and Neighborhood Purity index, computed either from the original publication’s cell type annotations \citep{pbmcsca} or from clustering results. The CH index measures the ratio of between-cluster to within-cluster dispersion.The Silhouette index compares each point’s similarity to its own cluster versus others. And Neighborhood Purity quantifies how often a point’s nearest neighbors belong to the same cluster. The results are summarized in Figure \ref{realdata2}.

We find that the shared-signal tracing algorithm (Shared-SVD) consistently outperforms other methods across all three evaluation criteria for each of the three embedding tasks. Specifically, Shared-SVD achieves a significantly higher CH Index compared to both Individual-SVD and Stack-SVD for cell embeddings. For neighborhood purity, the cell embeddings based on Shared-SVD show substantial improvements over those based on Individual-SVD, Stack-SVD and Average-SVD in the Drop-seq data, with notable gains also observed in the other two datasets. Besides, the comparative performance of Average-SVD to other methods is not stable measured in Neighborhood Purity. Additionally, Shared-SVD outperforms the other methods when evaluated using the Silhouette Index. Overall, these results highlight the superior performance of the shared-signal tracing algorithm in effectively estimating shared information across matrices.

\begin{figure}[h!]
    \centering
    \caption{Comparison of the Silhouette Index for three embedding methods. Two plots correspond to the cell embedding of two data matrices.}
    \label{Sil}
    \includegraphics[width=\textwidth]{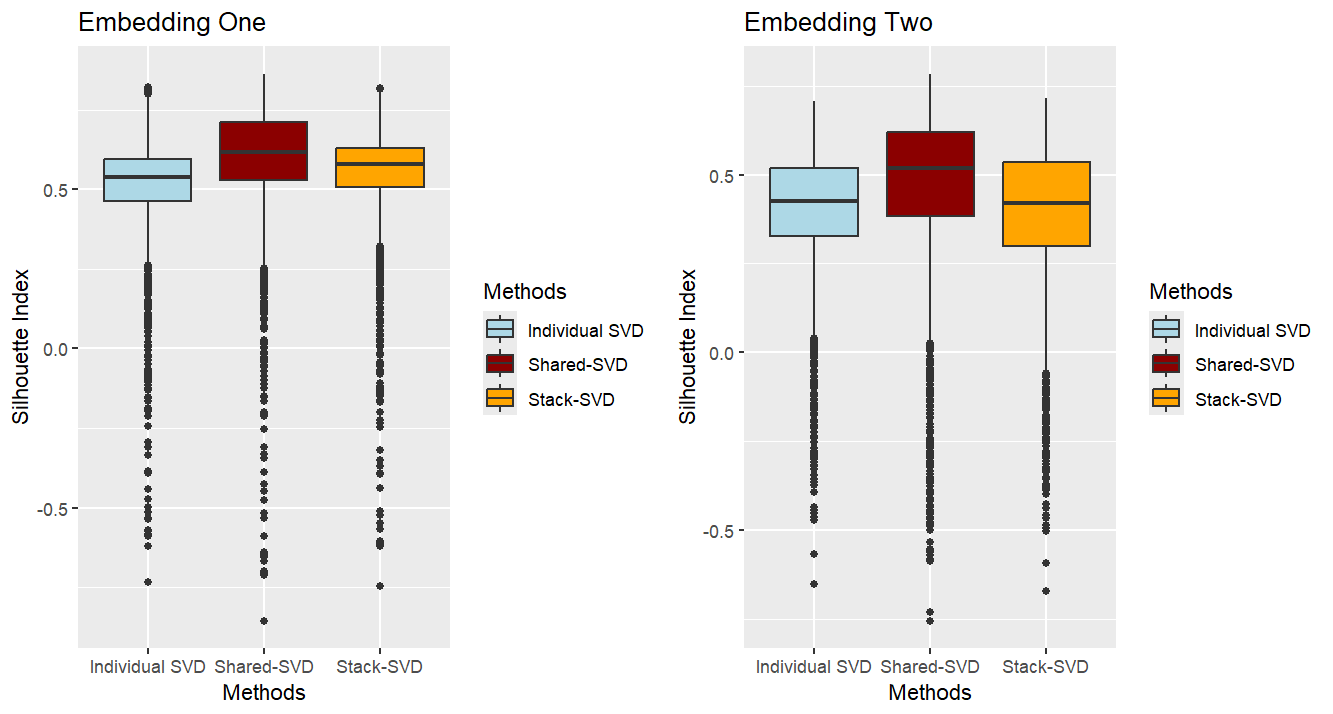}  

\end{figure}

\begin{figure}[h!]
    \centering
    \caption{Comparison of the Silhouette Index for three embedding methods grouped by different cell types. Two plots correspond to the cell embedding of two data matrices.} 
    \label{diffcellSil}
    \includegraphics[width=\textwidth]{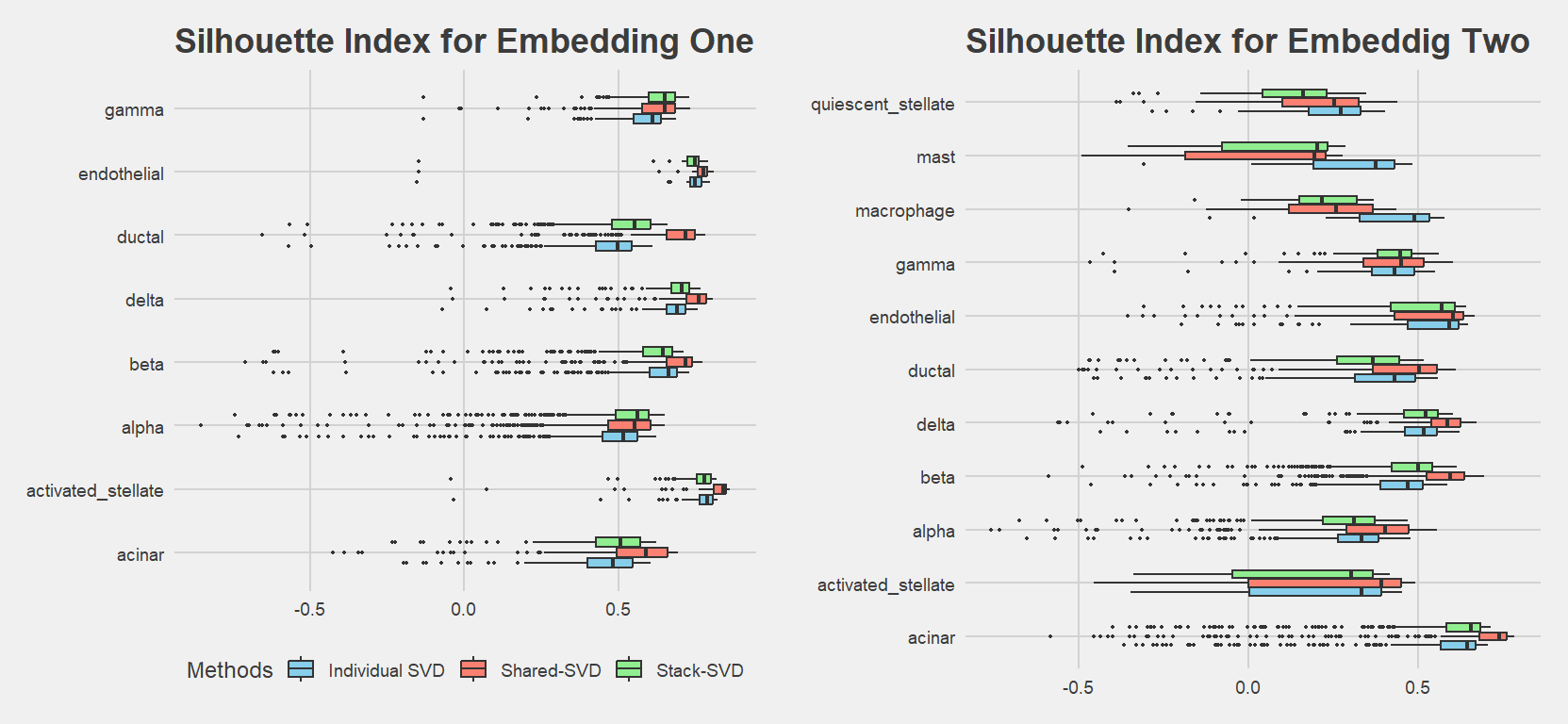}

\end{figure}

We also analyze the "panc8" dataset from the R package \texttt{SeuratData} \citep{panc8}. Figure \ref{Sil} presents a comparison of the Silhouette index across different embedding methods, while Figure \ref{diffcellSil} shows a boxplot of the Silhouette index across different cell types. Similar findings can be observed from the two figures that Shared-SVD outperforms other methods.

\section{Discussion}\label{discussion}

A limitation of our study is that the minimax optimality results for Stack-SVD  are established only when the dimensions of the individual matrices are comparable or when the signal strength is sufficiently large. However, our analysis suggests that when dimensions are not comparable, individual SVD estimators based on the smaller matrix can still achieve the optimal rate under certain signal constraints. An interesting direction for future research is to explore the minimax landscape when matrix dimensions differ significantly and to develop adaptive, minimax-optimal estimators that do not rely on signal constraints and remain effective in such scenarios.

Here, we propose a heuristic procedure to decide whether to use {Individual SVD} or {Stack-SVD} in practice. 
The key idea is to test whether the leading singular values of an observed matrix $Y$ significantly exceed the theoretical maximum expected under pure noise. 
According to {Bai--Yin's law} \cite{bai1993limit, Vershynin2012}, for an $n \times p$ random noise matrix $Z$ with independent, mean-zero, unit-variance entries, where $n/p\to(0,1]$ without loss of generality,
then it holds that $$
\sigma_{\max}(Z) = \sqrt{n} + \sqrt{p} + o(\sqrt{n})$$ almost surely.
Thus, if the top singular values of $Y$ exceeds this noise threshold, the matrix likely contains a strong signal; if it is near the threshold, the signal is weak.
Because the noise variance $\tau^2$ is unknown, we estimate it using the robust method of \cite{gavish2014}, 
$\hat{\tau} = \frac{y_{\text{med}}}{\sqrt{p\, \mu_\beta}}$,
where $y_{\text{med}}$ is the median singular value of $Y$, $\beta = n/p$, and $\mu_\beta$ is the median of the Marčenko--Pastur distribution.
Our {adaptive decision rule} is as follows: 1) For each data matrix $Y_k$, compute the threshold
    $
    T_k = \hat{\tau}_k(\sqrt{n} + \sqrt{p_k})
    $; and 2) Compare $\sigma_1(Y_k)$ to $T_k$: If only one matrix satisfies $\sigma_1(Y_k) \gg T_k$, use {Individual SVD}; if both satisfy the inequality, use {Stack-SVD}.
Intuitively, when $\sigma_1(Y_k)$ far exceeds the noise threshold, the matrix carries meaningful signal and should be included; otherwise, it likely adds noise and should be excluded. 
This approach naturally extends to settings with more than two matrices.

For the algorithm detecting shared and unshared singular vectors, we focus on the case where unshared vectors are orthogonal. The algorithm relies on the distance between vectors and remains effective even in the presence of mild non-orthogonality among unshared vectors. This is because mildly non-orthogonal unshared vectors generally still exhibit a relatively large $\sin \Theta$ distance. When the non-orthogonality is more pronounced, the algorithms still yield a result such that $\mathbb{J} \subset \widehat{\mathbb{J}}$. Investigating the degree of non-orthogonality that can be tolerated and extending it to handle non-orthogonal unshared singular vectors are intriguing for future research.

In this work, we focus exclusively on leveraging the shared singular subspace. However, when unshared singular vectors are not exactly orthogonal, there may be latent shared information embedded within them. For instance, if two unshared vectors form a small angle, their near alignment suggests potential shared structure. To simplify the problem, we do not account for such shared information within the unshared singular subspaces. A promising avenue for future exploration would involve redefining the notions of shared and unshared information and developing methods to capture shared structure arising from non-orthogonal unshared vectors.

\section{Acknowledgment}
The authors are grateful to the associate editor and the two anonymous referees for their valuable comments and suggestions, which have significantly enhanced the clarity and quality of this work. Rong Ma would like to thank Tavor Baharav, Phillip Nicol, Yue M. Lu, Cong Ma, Joshua Agterberg, TJ Channagiri and Ou Ren for helpful discussions. Most of this work is done while Zhengchi Ma was affiliated with Nankai University. This arXiv version is typeset in a single-column format for readability, while the official two-column version appears in the IEEE Transactions on Information Theory.

\appendix
\section{Proof of Theoretical Results}
\subsection{Proof of theorems and corollaries}
\begin{proof}[PROOF OF THEOREM \ref{upperlowerbound}]
    We first have the lemma:
\begin{lemma}\label{singularvaluerelation}
    For any finite integer $k\geq 2$, Let $X_1,X_2,...,X_k$ are $k$ rank $r$ matrices with identical left singular subspace $U$. $X$ is the noiseless stacked matrix, that is, $X=(X_1,X_2,...,X_k)$. Then we have $$\sigma_r^2(X)=\mathop{\min}\limits_{1 \leq i \leq r} \sum_{j=1}^k \sigma_{(i)}^2(X_j),
    $$
    where $\sigma_{(i)}(X_j)$ denotes the  singular value of $X_j$ corresponding to the $i$-th left singular vector.
\end{lemma}
Then conduct SVD
$$
    (
        X_1~X_2)=(
        U\Sigma_1V^T_1~U\Sigma_2V^T_2)=U(
        \Sigma_1~\Sigma_2)\mathrm{diag}(
        V_1^T ,~V_2^T=U\sqrt{\Sigma_1^2+\Sigma_2^2}\begin{pmatrix}
        \frac{\Sigma_1}{\sqrt{\Sigma_1^2+\Sigma_2^2}}V_1^T & \frac{\Sigma_2}{\sqrt{\Sigma_1^2+\Sigma_2^2}}V_2^T
    \end{pmatrix}.$$
Thus, the left singular subspace of $(X_1~X_2)$ is $U$, with the minimum singular value being
\[
\sigma_r^2(X)=\mathop{\min}\limits_{1 \leq i \leq r} \{\sigma_{(i)}^2(X_1)+\sigma_{(i)}^2(X_2)\}
\] by Lemma \ref{singularvaluerelation}. By same procedure with Theorem 3 in \cite{caizhang2018}, 
$$
        \mathbb{E}\|\sin\Theta(U,\hat{U})\|^2 \leq \frac{cn\big(\mathop{\min}\limits_{1 \leq i \leq r} \{\sigma_{(i)}^2(X_1)+\sigma_{(i)}^2(X_2)\}+(p_1+p_2)\big)}{\big(\mathop{\min}\limits_{1 \leq i \leq r} \{\sigma_{(i)}^2(X_1)+\sigma_{(i)}^2(X_2)\}\big)^2}\wedge 1.$$
Also by noting 
$$
\mathop{\min}\limits_{1 \leq i \leq r} \{\sigma_{(i)}^2(X_1)+\sigma_{(i)}^2(X_2)\}\geq \gamma^2,$$
we have 
$$
        \mathop{\sup}\limits_{X\in \mathscr{F}_{r,\gamma}}\mathbb{E}\|\sin\Theta(U,\hat{U})\|^2 \leq \frac{cn(\gamma^2+(p_1+p_2))}{\gamma^4}\wedge 1.$$
For minimax lower bound, we note the parameter space (\ref{eq:paraspace1}) is equal to (\ref{eq:101}).
Thus we complete the proof by Proposition \ref{generalizecainzhanglowershare}.
\end{proof}

\begin{proof}[Proof of Theorem \ref{suboptimalave}]
    For brevity, we denote the value $\hat{u}_1^T\hat{u}_2=\eta \in [-1,1]$. Take SVD of $(\hat{u}_1~\hat{u}_2)$, we have, $$
(\hat{u}_1~\hat{u}_2)=(\frac{1}{\sqrt{2}\sqrt{1+\eta}}\hat{u}_1+\frac{1}{\sqrt{2}\sqrt{1+\eta}}\hat{u}_2 \quad \frac{1}{\sqrt{2}\sqrt{1-\eta}}\hat{u}_1-\frac{1}{\sqrt{2}\sqrt{1-\eta}}\hat{u}_2)\begin{pmatrix}
        \sqrt{1+\eta} & 0 \\ 0 & \sqrt{1-\eta}
    \end{pmatrix}\begin{pmatrix}
        \frac{1}{\sqrt{2}} & \frac{1}{\sqrt{2}} \\ \frac{1}{\sqrt{2}} & -\frac{1}{\sqrt{2}}
    \end{pmatrix}.$$
Then Average-SVD is $$
\hat{u}=\frac{1}{\sqrt{2}\sqrt{1+\eta}}\hat{u}_1+\frac{1}{\sqrt{2}\sqrt{1+\eta}}\hat{u}_2.$$
Also
$$
\|\sin \Theta(\hat{u}_1,\hat{u})\|=\sqrt{1-\sigma_{\min}^2\left(\frac{\hat{u}_1^T(\hat{u}_1+\hat{u}_2)}{\sqrt{2}\sqrt{1+\eta}}\right)}=\frac{\sqrt{1-\eta}}{\sqrt{2}},$$
$$
\|\sin \Theta(\hat{u}_1,\hat{u}_2)\|=\sqrt{1-\sigma_{min}^2(\hat{u}_1^T\hat{u}_2)}=\sqrt{1+\eta}\sqrt{1-\eta}.$$
Combining inequalities above, 
$$
\|\sin\Theta(\hat{u}_1,\hat{u})\|=\frac{1}{\sqrt{2}\sqrt{1+\eta}}\|\sin\Theta(\hat{u}_1,\hat{u}_2)\|.$$
Thus
\begin{align*}&\|\sin\Theta(u,\hat{u})\|\geq  \|\sin\Theta(u,\hat{u}_1)\|-\frac{1}{\sqrt{2}\sqrt{1+\eta}}\|\sin\Theta(u,\hat{u}_1)\|-\frac{1}{\sqrt{2}\sqrt{1+\eta}}\|\sin\Theta(u,\hat{u}_2)\|\\&\qquad\geq \frac{\sqrt{2}\sqrt{1+\eta}-1}{\sqrt{2}\sqrt{1+\eta}}\|\sin\Theta(u,\hat{u}_1)\|-\frac{1}{\sqrt{2}\sqrt{1+\eta}}\|\sin\Theta(u,\hat{u}_2)\|\geq \frac{\sqrt{2}-1}{\sqrt{2}}\|\sin\Theta(u,\hat{u}_1)\|-\frac{1}{\sqrt{2}}\|\sin\Theta(u,\hat{u}_2)\|.\end{align*}
For large signal, by Theorem 3 and 4 in \cite{caizhang2018}, we have 
$$
    \mathbb{E}\|\sin\Theta(u,\hat{u}_2)\| \leq c\sqrt{\frac{n(\beta^2+p_2)}{\beta^4} },$$
 and $$\mathop{\sup}\limits_{X\in \mathscr{F}_{1,\gamma}}\mathbb{E}\|\sin\Theta(u,\hat{u}_1)\| \geq \mathop{\inf}\limits_{\Tilde{u}} \mathop{\sup}\limits_{X\in \mathscr{F}_{1,\gamma}} \mathbb{E}\|\sin\Theta(u,\Tilde{u})\|
    \geq c\sqrt{\frac{n(\alpha^2+p_1)}{\alpha^4}}.$$
Therefore, we have $$
    \mathbb{E}\|\sin\Theta(u,\hat{u})\|\geq  \frac{\sqrt{2}-1}{\sqrt{2}}\mathbb{E}\|\sin\Theta(u,\hat{u}_1)\|-c_2\sqrt{\frac{n(\beta^2+p_2)}{\beta^4} }.$$
Taking $\sup$,  $$
   \mathop{\sup}\limits_{X\in \mathscr{F}_{1,\gamma}}\mathbb{E}\|\sin\Theta(u,\hat{u})\|\geq c_1\sqrt{\frac{n(\alpha^2+p_1)}{\alpha^4}}-c_2\sqrt{\frac{n(\beta^2+p_2)}{\beta^4}}. $$
By the assumption that $\beta \gg \alpha$ and $\alpha,\beta$ are large enough, $$
c_1\sqrt{\frac{n(\alpha^2+p_1)}{\alpha^4}}-c_2\sqrt{\frac{n(\beta^2+p_2)}{\beta^4}}>0.$$
Therefore, combining inequalities above, we complete the proof.
\end{proof}

\begin{proof}[PROOF OF THEOREM \ref{generalbounds}]
Take SVD:
\begin{align*}
    (X_1~X_2)
    &=\begin{pmatrix}
        U_r & U_{1*} & U_{2*}
    \end{pmatrix}\mathrm{diag}( 
    \sqrt{\Sigma_{1r}^2+\Sigma_{2r}^2},~\Sigma_{1*},~\Sigma_{2*})\begin{pmatrix}
        \frac{\Sigma_{1r}}{\sqrt{\Sigma_{1r}^2+\Sigma_{2r}^2}} & 0 & \frac{\Sigma_{2r}}{\sqrt{\Sigma_{1r}^2+\Sigma_{2r}^2}} & 0 \\
        0 & I & 0 & 0 \\ 0& 0 & 0 & I
    \end{pmatrix} V^T\\
    &=\begin{pmatrix}
        U_r & U_{1*} & U_{2*}
    \end{pmatrix}\mathrm{diag}( 
    \sqrt{\Sigma_{1r}^2+\Sigma_{2r}^2},~\Sigma_{1*},~\Sigma_{2*}) V_*^T,
\end{align*}
where in the last step we merge the two orthogonal matrices.
Then sort the singular values decreasingly, which is equivalent to introducing a permutation matrix $P$ (with only a single $1$ in each row and column, and other entries are all $0$). Denoting the decreasing $\Sigma$ as $\Sigma_D$, we have
$$
(X_1~X_2)=(
        U_r ~ U_{1*} ~ U_{2*})P\Sigma_{D}P^T V_*^T.$$
We say shared vectors (in $U_r$) and unshared vectors (in $U_{1*},U_{2*}$) are of  different vector types. For the left singular vectors, after permuted into a decreasing order, we define all the consecutive vectors with same vector type as a vector group. The $i$-th shared group is denoted as $U_{share(i)}$ while the $i$-th unshared group is denoted as $U_{unshare(i)}$. For some number of groups $k$, there might be four possibilities for the notation of the left singular subspace of $(X_1~X_2)$ after sorting singular values in decreasing order:
$$
(U_{share(1)}~U_{unshare(1)}~...~U_{share(k)}~U_{unshare(k)}),$$ $$
(U_{share(1)}~U_{unshare(1)}~...~U_{unshare(k-1)}~U_{share(k)}),$$ $$
(U_{unshare(1)}~U_{share(1)}~...~U_{unshare(k)}~U_{share(k)}),$$ $$
(U_{unshare(1)}~U_{share(1)}~...~U_{share(k-1)}~U_{unshare(k)}).$$
We provide the proof for the first scenario while others scenarios are analogous.
Then the decreasing SVD is 
$$
(X_1~X_2)=(U_{share(1)}~U_{unshare(1)}~...~U_{share(k)}~U_{unshare(k)})\Sigma_D V^T,$$
denoted for brevity as $$(X_1~X_2)=(U_{s(1)}~U_{u(1)}~U_{s(2)}~U_{u(2)}~...~U_{s(k)}~U_{u(k)})\Sigma_D V^T.$$
Also, here we denote the minimum singular value associating with vectors in $U_{share(i)}$ (or $U_{unshare(i)}$) as $\sigma_{\min}(U_{share(i)})$ (or $\sigma_{\min}(U_{unshare(i)})$). Recall gap condition at each switch of vector types, $$g_{i(s)}^2(X)>c\sigma_{s+1}^2(X), \forall s\in\mathbb{S},$$
which is equivalent to $$
g_{i(s)}^2(X)=\sigma_s^2-\sigma_{s+1}^2>c\sigma_{s+1}^2(X),$$
 that is
 \(
 \sigma_{s}^2>(1+c)\sigma_{s+1}^2
 \), where $c>0$ is a small constant. It can be written as $$
\sigma_{\min}^2(U_{share(i)})>(1+\epsilon)\sigma_{\min}^2(U_{unshare(i)}),$$ and $$ \sigma_{\min}^2(U_{unshare(j)})>(1+\epsilon)\sigma_{\min}^2(U_{share(j+1)})$$
for all $i=1,2,...,k$ and $j=1,...,k-1$, where $\epsilon>0$ is a small constant. Denote the number of vectors in $U_{share(i)}$ and $U_{unshare(i)}$ as $n_{si}$ and $n_{ui}$ respectively. Denote the top $i$ left singular vectors of  $(X_1~X_2)$ as $U_{(i)}$ and left singular vectors of  $(Y_1~Y_2)$ at the corresponding positions as the same version with hat. Then 
\begin{align*}
&\|\sin\Theta(\hat{U}_r^{\mathbb{J}},U_r)\|^2\leq \|\hat{U}_r^{\mathbb{J}}\hat{U}_r^{\mathbb{J}T}-U_rU_r^T\|^2= \bigg\|\sum_{i=1}^kU_{s(i)}U_{s(i)}^T-\sum_{i=1}^k\hat{U}_{s(i)}\hat{U}_{s(i)}^T\bigg\|^2\\
&\leq c\|\sin\Theta((U_{s(1)}~U_{u(1)}~U_{s(2)}~U_{u(2)}~...~U_{s(k)}),(\hat{U}_{s(1)}~\hat{U}_{u(1)}~\hat{U}_{s(2)}~\hat{U}_{u(2)}~...~\hat{U}_{s(k)}))\|^2+c\bigg\|(\sum_{i=1}^{k-1}U_{u(i)}U_{u(i)}^T-\sum_{i=1}^{k-1}\hat{U}_{u(i)}\hat{U}_{u(i)}^T)\bigg\|^2 \\
&\leq c\|\sin\Theta((U_{s(1)}~U_{u(1)}~U_{s(2)}~U_{u(2)}~...~U_{s(k)}),(\hat{U}_{s(1)}~\hat{U}_{u(1)}\hat{U}_{s(2)}~\hat{U}_{u(2)}~...~\hat{U}_{s(k)}))\|^2\\
&\qquad+c\|\sin\Theta((U_{s(1)}~U_{u(1)}~U_{s(2)}~U_{u(2)}~...~U_{u(k-1)}),(\hat{U}_{s(1)}~\hat{U}_{u(1)}~\hat{U}_{s(2)}~\hat{U}_{u(2)}~...~\hat{U}_{u(k-1)}))\|^2\\
&\qquad+... +  c\|\sin\Theta((U_{s(1)}~U_{u(1)})),(\hat{U}_{s(1)}~\hat{U}_{u(1)}))\|^2 + c\|\sin\Theta(U_{s(1)},\hat{U}_{s(1)}\|^2 
\\&\leq  k\frac{cn(\sigma_{\min}^2(U_{shared(k)})+(p_1+p_2))}{\sigma_{\min}^4(U_{shared(k)})}\wedge 1,
\end{align*}
where $k$ is the number of shared or unshared groups. Recall the definition that, after sorting singular values in decreasing order, suppose there are $N$ vector type switches and suppose the $i$-th vector type switch happen between $s$ and $s+1$-th vector, we denote the $i$-th gap as $\sigma_s^2-\sigma_{s+1}^2:=g_{i(s)}^2$.
By the conditions $g_{i(s)}^2\geq c\sigma_{s+1}^2$, $c>0$, we know $$\sigma_s^2\geq g_{i(s)}^2= \sigma_s^2-\sigma_{s+1}^2\geq c\sigma_s^2,$$ which means $g_i^2$ and $\sigma_s^2$ are in the same order. 
(if the last singular vector in the stack is shared, we still denote it as a gap by  $g^2_{N+1}=\sigma_{\min}^2(U_{shared(k)})-0$). Note that 
$$
\sigma_{\min}^2(U_{shared(k)})\geq g_N^2 \geq \mathop{\min}\limits_{1 \leq i\leq N}g_i^2.$$ 
We have $$k \frac{cn(\sigma_{\min}^2(U_{shared(k)})+(p_1+p_2))}{\sigma_{\min}^4(U_{shared(k)})}\leq k \frac{cn(\sigma_{\min}^2(U_{shared(k)})+(p_1+p_2))}{(\mathop{\min}\limits_{1 \leq i\leq N}g_i^2)^2}.$$
Also when \(s=\mathop{\arg\min}\limits_{1 \leq i\leq N}g_i^2\), we know $$
\sigma_{\min}^2(U_{shared(k)})\leq \sigma_s^2\leq cg_s^2=c\mathop{\min}\limits_{1 \leq i\leq N}g_i^2.$$
Therefore $$k \frac{cn\big(\sigma_{\min}^2(U_{shared(k)})+(p_1+p_2)\big)}{\sigma_{\min}^4(U_{shared(k)})}\leq k \frac{cn\big(\mathop{\min}\limits_{1 \leq i\leq N}g_i^2+(p_1+p_2)\big)}{\big(\mathop{\min}\limits_{1 \leq i\leq N}g_i^2\big)^2}.$$
Times that vector type switch happens should be less than or equal to $r$, thus $k\leq r\leq c$.
Therefore, 
$$
\|\sin\Theta(\hat{U}_r^{\mathbb{J}},U_r)\|^2\leq \frac{cn\left(\mathop{\min}\limits_{1 \leq i\leq N}g_i^2+(p_1+p_2)\right)}{\left(\mathop{\min}\limits_{1 \leq i\leq N}g_i^2\right)^2}\wedge 1.$$
Finally considering that $$\mathop{\min}\limits_{1 \leq i\leq N}g_i^2\geq t^2$$ and the bound is a decreasing function, the result holds that
$$
        \mathop{\sup}\limits_{X\in \mathscr{H}_{r,t}}\mathbb{E}\|\sin\Theta(U_r,\hat{U}_r^{\mathbb{J}})\|^2 \leq \frac{cn(t^2+(p_1+p_2))}{t^4}\wedge 1
.$$
For minimax lower bound, recall the parameter space (\ref{eq:paraspace2}).
By previous analysis of the SVD for stacked matrix, we have
$$(X_1~X_2)
    =(U_r~U_{1*}~U_{2*})\mathrm{diag}( 
    \sqrt{\Sigma_{1r}^2+\Sigma_{2r}^2},~\Sigma_{1*},~\Sigma_{2*})V_*^T.$$
Also note that, under assumptions, the minimum singular value associated to the shared singular vectors (denoted as $\sigma_{\min}^2(U_{r})$), and the minimum eigen-gap $\mathop{\min}\limits_{1 \leq k \leq N}g_k^2$ have the same order: $$
\sigma_{\min}^2(U_{r})\geq g_N^2\geq \mathop{\min}\limits_{1 \leq k \leq N}g_k^2 \geq c\sigma_s^2 \geq c\sigma_{\min}^2(U_{r}),$$
where $s=\text{argmin}_k g_k^2$.
Thus the parameter space considered is equivalent to the following for constant $c_s$:
    \begin{align*}
\mathscr{H}_{r,\gamma}^{(u)} = \left\{\begin{array}{c}X= \in \mathbb{R}^{n \times (p_1+p_2)}:~\text{rank}(X)=r+r_{1*}+r_{2*},\\X~has~form~(\ref{prop2SVD}),
      U_r \in \mathbb{O}(n,r),U_{i*}\in \mathbb{O}(n,r_{i*}),\\U_{1*}^TU_{2*}=0, V_i \in \mathbb{O}(p_i,r), g_{i(s)}^2(X)>c\sigma_{s+1}^2(X), \\\forall s\in\mathbb{S}, \mathop{\min}\limits_{1 \leq k \leq r}\sigma_{(k)}^2(X) \geq c_s\gamma^2, i=1,2
\end{array} \right\}.
    \end{align*}
    where $\sigma_{(k)}(X)$ denotes the singular value in $X$ corresponding to the $k$-th vector in the shared $U_r$.
Therefore, by Proposition \ref{generalizacainzhanglower}, we have
$$
        \mathop{\inf}\limits_{\Tilde{U}_r}\mathop{\sup}\limits_{X\in \mathscr{H}_{r,\gamma}^{(u)}} \mathbb{E}\|\sin\Theta(U_r,\Tilde{U_r})\|_F^2 \geq c\bigg(\frac{nr(\gamma^2+(p_1+p_2))}{\gamma^4} \wedge r\bigg).
    $$
This completes our proof.
\end{proof}

\begin{proof}[PROOF OF THEOREM \ref{generalize}]
Similarly decompose:
    $$(X_1~X_2)=(
        U_r ~ U_{1*} ~ U_{2*})\mathrm{diag}( 
    \sqrt{\Sigma_{1r}^2+\Sigma_{2r}^2},~\Sigma_{1*},~\Sigma_{2*})V_*^T.$$
And for $$\phi_i:=\sqrt{\sigma_{(i)}^2(X_1)+\sigma_{(i)}^2(X_2)},\qquad (\Sigma_{1r}^2+\Sigma_{2r}^2)^{\frac{1}{2}}=\mathrm{diag}(\phi_1,...\phi_r).$$
By assumption (\ref{cond.s1}),
 the diagonal of $(\Sigma_{1r}^2+\Sigma_{2r}^2)^{\frac{1}{2}}$ are exactly at the top $r$ positions among singular values of $(X_1~X_2)$.
Also, 
$$
\sigma_r(\begin{pmatrix}
    X_1 & X_2
\end{pmatrix})=\mathop{\min}\limits_{1 \leq i \leq r}\sqrt{\sigma_{(i)}^2(X_1)+\sigma_{(i)}^2(X_2)}.$$
Note that the eigen-gap condition of Proposition \ref{generalcainzhang} ($\sigma_r^2>(1+\epsilon)\sigma_{r+1}^2$) is satisfied by
$$c_1\mathop{\min}\limits_{1 \leq k \leq r}(\sigma_{(k)}^2(X_1)+\sigma_{(k)}^2(X_2)) > (1+\epsilon)
\mathop{\max}\limits_{r+1 \leq k \leq \textup{rank}(X_1)}\sigma_{(k)}^2(X_1) \vee \mathop{\max}\limits_{r+1 \leq k \leq \textup{rank}(X_2)}\sigma_{(k)}^2(X_2),$$ with some small constant $\epsilon>0$.
 By Proposition \ref{generalcainzhang}, we have :
$$\mathbb{E}\|\sin\Theta(U_r,\hat{U})\|^2 \leq \frac{cn((\mathop{\min}\limits_{1 \leq i \leq r}\sqrt{\sigma_{(i)}^2(X_1)+\sigma_{(i)}^2(X_2)})^2+(p_1+p_2))}{(\mathop{\min}\limits_{1 \leq i \leq r}\sqrt{\sigma_{(i)}^2(X_1)+\sigma_{(i)}^2(X_2)})^4}\wedge 1.$$
Note $$\mathop{\min}\limits_{1 \leq i \leq r}\sqrt{\sigma_{(i)}^2(X_1)+\sigma_{(i)}^2(X_2)}$$ is no less than the minimum singular value gap $$
    t^2 \leq \mathop{\min}\limits_{1 \leq i \leq r}\{\sigma_{(i)}^2(X_1)+\sigma_{(i)}^2(X_2)\}-\mathop{\max}\limits_{r+1 \leq i \leq rank(X_1)}\sigma_{(i)}^2(X_1) \vee \mathop{\max}\limits_{r+1 \leq i \leq rank(X_2)}\sigma_{(i)}^2(X_2)
    \leq \mathop{\min}\limits_{1 \leq i \leq r}\sqrt{\sigma_{(i)}^2(X_1)+\sigma_{(i)}^2(X_2)},$$ we have $$
\mathbb{E}\|\sin\Theta(U_r,\hat{U})\|^2 \leq \frac{cn(t^2+(p_1+p_2))}{t^4}\wedge 1.$$
Then the upper bound result holds.
For minimax lower bound, note that when in  the proof of Theorem \ref{generalbounds} and Proposition \ref{generalizacainzhanglower}, we do not specify the singular value for the unshared vectors, which make it  general. We now only need to specify $\Sigma_{1*}^2$ and $\Sigma_{2*}^2$ with a smaller rate than $\Sigma_{1r}^2+\Sigma_{2r}^2$, then this case is a specific instance. For example, for strong signal($\gamma^2 \gtrsim p_{V1}\wedge p_{V2}$) in Proposition \ref{generalizacainzhanglower}, construct
 $$
\Sigma_{1*}=\frac{\sqrt{\alpha^2+\beta^2}}{2}I_{1*} ,\qquad \Sigma_{2*}=\frac{\sqrt{\alpha^2+\beta^2}}{2}I_{2*},$$
while for weak signal($p_{V1}\asymp p_{V2}$ and $\gamma^2 \lesssim p_{V1}\wedge p_{V2}$), construct
$$
\Sigma_{1*}=\frac{\gamma}{2}I_{1*},\qquad \Sigma_{2*}=\frac{\gamma}{2}I_{2*}.$$
Other part of the proof can be derived analogously as the proof of Theorem \ref{generalbounds}. Thus we have finished the proof of the theorem.
\end{proof}

\begin{proof}[PROOF OF THEOREM \ref{strongcontamination}]
Still we only have to consider the spectral norm for $\sin\Theta$ distance. Similar to the proof of Theorem \ref{generalize}, we decompose
    $$\begin{pmatrix}
        X_1 & X_2
    \end{pmatrix}=\begin{pmatrix}
         U_{1*} & U_{2*} & U_{r}
    \end{pmatrix}\mathrm{diag}(
    \Sigma_{1*},~\Sigma_{2*},~(\Sigma_{1r}^2+\Sigma_{2r}^2)^{\frac{1}{2}}
    )V^T
    :=\begin{pmatrix}
         U_{c}  & U_{r}
    \end{pmatrix}\mathrm{diag}(
          \Sigma_{c},~(\Sigma_{1r}^2+\Sigma_{2r}^2)^{\frac{1}{2}}
    )V^T.$$
Then by triangle inequality
\begin{align*}
\|\sin\Theta(U_r,\hat{U}_r^S)\|^2 \leq \|U_rU_r^T-\hat{U}_r^S\hat{U}_r^{ST}\|^2\leq 2\bigg\|(
        U_c ~ U_r)\begin{pmatrix}
        U_c^T\\
        U_r^T
    \end{pmatrix}-(\hat{U_c} ~ \hat{U}_r^S)\begin{pmatrix}
        \hat{U_c}^T\\
        \hat{U}_r^{ST}
    \end{pmatrix}\bigg\|^2+2\|U_cU_c^T-\hat{U_c}\hat{U_c}\|^2\\
    \leq 4\|\sin\Theta((U_c~U_r),(\hat{U_c}~\hat{U}_r^S))\|^2 + 4\|\sin\Theta(U_c,\hat{U_c})\|^2.\end{align*}
For the first term, by theorem 3 in \cite{caizhang2018},
$$
\mathbb{E}\|\sin\Theta((U_c~U_r),(\hat{U_c}~\hat{U}_r^S))\|^2 \leq \frac{cn(\sigma_{\min}^2(\begin{pmatrix}
    X_1 & X_2
\end{pmatrix})+(p_1+p_2))}{\sigma_{\min}^4(\begin{pmatrix}
    X_1 & X_2
\end{pmatrix})} \wedge 1.$$
For the second term, by Proposition \ref{generalcainzhang}, $$
\|\sin\Theta(U_c,\hat{U_c})\|^2 \leq \frac{cn(\sigma_{d}^2(\begin{pmatrix}
    X_1 & X_2
\end{pmatrix})+(p_1+p_2))}{\sigma_{d}^4(\begin{pmatrix}
    X_1 & X_2
\end{pmatrix})} \wedge 1,$$ 
where the condition of Proposition \ref{generalcainzhang} is satisfied because 
$$
\mathop{\max}\limits_{1 \leq k \leq r}\{\sigma_{(d_1+k)}^2(X_1)+\sigma_{(d_2+k)}^2(X_2)\} < c_1( \sigma_{d_1}^2(X_1) \wedge \sigma_{d_2}^2(X_2))$$
with $c_1 \in (0,1)$ implies
$$
( \sigma_{d_1}^2(X_1) \wedge \sigma_{d_2}^2(X_2))>(1+\epsilon)\mathop{\max}\limits_{1 \leq k \leq r}\{\sigma_{(d_1+k)}^2(X_1)+\sigma_{(d_2+k)}^2(X_2)\}
$$ for some small constant $\epsilon>0$.
Combining two terms while noting  
$$
\sigma_{d}^2(\begin{pmatrix}
    X_1 & X_2
\end{pmatrix}) \geq \sigma_{\min}^2(\begin{pmatrix}
    X_1 & X_2
\end{pmatrix}) = \mathop{\min}\limits_{1 \leq i \leq r}\{\sigma_{(d_1+i)}^2(X_1)+\sigma_{(d_2+i)}^2(X_2)\} \geq t^2,$$
we have 
$$
\mathbb{E}\|\sin\Theta(U_r,\hat{U}_r^S)\|^2 \leq \frac{cn(t^2+(p_1+p_2))}{t^4} \wedge 1.$$ 
Thus the upper bound holds. For minimax lower bound, we construct specific $\Sigma_{1*}$ and $\Sigma_{2*}$ similar to the proof of Theorem \ref{generalize}. However, we construct some strong unshared signals here.  For example, for strong signal ($\gamma^2 \gtrsim p_{V1}\wedge p_{V2}$) in Proposition \ref{generalizacainzhanglower}, we construct
$$
\Sigma_{1*}=2\sqrt{\alpha^2+\beta^2}I_{1*},\qquad\Sigma_{2*}=2\sqrt{\alpha^2+\beta^2}I_{2*},$$
while for weak signal ($\gamma^2 \lesssim p_{V1}\wedge p_{V2}$), we construct  $$
\Sigma_{1*}=2\gamma I_{1*},\qquad \Sigma_{2*}=2\gamma I_{2*}.$$
Here $\Sigma_{1*}$ and $\Sigma_{2*}$ does not affect the final outcome because the Theorem \ref{generalbounds} is general for all choices of them under some conditions. Thus we finish the proof by similar procedure with Theorem \ref{generalbounds}.
\end{proof}

\begin{proof}[PROOF OF THEOREM \ref{nonorthosvd}]

    For two matrices
$$X_1=(
    U_r ~ U_{1*})\mathrm{diag}(
    \Sigma_1,\Sigma_{1*})V_1^T,\qquad 
X_2=(
    U_r ~ U_{2*})\mathrm{diag}(
    \Sigma_2,~\Sigma_{2*})V_2^T,$$
    similar to the orthogonal unshared vectors case, we have the decomposition
    $$
    (X_1~X_2)=(
        U_r ~ U_{1*} ~ U_{2*})\mathrm{diag}(
        \sqrt{\Sigma_1^2+\Sigma_2^2},~\Sigma_{1*},~\Sigma_{2*})V^T,$$
$$V^T=\begin{pmatrix}
        \frac{\Sigma_{1}}{\sqrt{\Sigma_{1}^2+\Sigma_{2}^2}}V_{1}^T & \frac{\Sigma_{2}}{\sqrt{\Sigma_{1}^2+\Sigma_{2}^2}}V_{2}^T\\
        V_{1*}^T & 0 \\ 0 & V_{2*}^T
    \end{pmatrix}.$$
    This is not the SVD because the leftmost matrix is not orthonormal. We first suppose the sub-matrix $(U_{1*}~U_{2*})$ is not a singular matrix in columns. Then by definition, the matrix $(U_{1*}~U_{2*})^T(U_{1*}~U_{2*})$ is a positive definite matrix with eigen-decomposition $$
    (U_{1*}~U_{2*})^T(U_{1*}~U_{2*})=\Gamma^T \Sigma \Gamma ,
    $$
    where $\Sigma$ is a diagonal matrix with all positive entries in diagonal and $\Gamma \in \mathbb{O}(r_{1*}+r_{2*},r_{1*}+r_{2*})$.
    We construct the matrices
    $$
    S=\Gamma^T \Sigma^{-\frac{1}{2}} \Gamma,\quad R=\mathrm{diag}(
        I ,~S),\quad R^{-1}=\mathrm{diag}(
        I,~S^{-1}).$$
    Note that $$ S^T(U_{1*}~U_{2*})^T(U_{1*}~U_{2*})S=I.$$
    Therefore $(U_{1*}~U_{2*})S$ is an orthogonal matrix.
    Then we have
    \begin{align*}(X_1~X_2)=(
        U_r ~ U_{1*} ~ U_{2*})RR^{-1}\mathrm{diag}(
        \sqrt{\Sigma_1^2+\Sigma_2^2},~\Sigma_{1*},~\Sigma_{2*}
    )V^T\\
    =\begin{pmatrix}
        U_r~ & (U_{1*}~U_{2*})S
    \end{pmatrix}\mathrm{diag}(
        \sqrt{\Sigma_1^2+\Sigma_2^2},~S^{-1}\mathrm{diag}(
              \Sigma_{1*},~\Sigma_{2*}))V^T.\end{align*}
    Further, for non-singular matrix $S^{-1}\mathrm{diag}(
              \Sigma_{1*},~\Sigma_{2*})$, we have SVD
          $$
          S^{-1}\mathrm{diag}(
              \Sigma_{1*},~\Sigma_{2*})=U^*\Sigma^*V^{*T}.
          $$
    Then 
    $$
        (X_1~X_2)=\begin{pmatrix}
        U_r~ & (U_{1*}~U_{2*})SU^*
    \end{pmatrix}\mathrm{diag}(
        \sqrt{\Sigma_1^2+\Sigma_2^2},~\Sigma^*)V^T.$$
    When $(U_{1*}~U_{2*})$ is a singular matrix in columns with rank $r^*$, we can first conduct a column linear transformation by right multiplying the column elementary operation matrix $L\in \mathbb{R}^{(r_{1*}+r_{2*})\times(r_{1*}+r_{2*})}$ such that
    $$
    (U_{1*}~U_{2*})L=(u_1~u_2~...~u_{r^*}~0~...~0)=(\Tilde{U}^*~{\bf 0}_{n\times (r_{1*}+r_{2*}-r^*)})
    ,$$
    where $u_i$ are unit vectors. The $L$ can be easily constructed by first constructing an elementary operation matrix to move $r^*$ vectors that cannot be expressed as linear combinations of each other to the leftmost $r^*$ positions. Then it can be multiplied by another elementary operation matrix to make the last vectors zero by subtracting their linear combination by precious $r^*$ vectors.   Then we construct $S$ similarly for $\Tilde{U}^*$ as stated in our theorem. That is $S=\Gamma^T \Sigma^{-\frac{1}{2}} \Gamma $ based on $\Tilde{U}^{*T}\Tilde{U}^{*}=\Gamma^T \Sigma \Gamma$, where $\Gamma \in \mathbb{O}(r^*,r^*)$. We take 
    $$
    S_1=\begin{pmatrix}
        S\\
        {\bf 0}_{(r_{1*}+r_{2*}-r^*)\times r^*}
    \end{pmatrix}, ~ R_1=\mathrm{diag}(
        I,~\begin{pmatrix}
             S\\ {\bf 0}_{(r_{1*}+r_{2*}-r^*)\times r^*}
         \end{pmatrix}),
$$
$$
R_1^{-1*}=\mathrm{diag}(
        I,~\begin{pmatrix}
             S^{-1}& {\bf 0}_{r^* \times (r_{1*}+r_{2*}-r^*)}
         \end{pmatrix}), ~ R_1R_1^{-1*}=\mathrm{diag}(
        I,~I,~{\bf 0}_{(r_{1*}+r_{2*}-r^*)\times (r_{1*}+r_{2*}-r^*)}).$$
    Also note that
    $$
    (U_r~u_1~,...,~u_{r^*}~0~,...,~0)=(U_r~u_1~,...,~u_{r^*}~0~,...,~0)\mathrm{diag}(
        I,~I,~{\bf 0}_{(r_{1*}+r_{2*}-r^*)\times (r_{1*}+r_{2*}-r^*)}).$$
    So we have 
    \begin{align*}
        (X_1~X_2)&=(
        U_r ~ U_{1*} ~ U_{2*})\mathrm{diag}(
        I,~LL^{-1})\mathrm{diag}(
        \sqrt{\Sigma_1^2+\Sigma_2^2},~\Sigma_{1*},~\Sigma_{2*})V^T\\
        &=(
        U_r ~ u_1,...,u_{r^*} ~ 0,...,0)R_1R_1^{-1*}\mathrm{diag}(
        \sqrt{\Sigma_1^2+\Sigma_2^2},~L^{-1}\mathrm{diag}(
              \Sigma_{1*},~\Sigma_{2*}))V^T
              \\&=(
        U_r~  (u_1,...,u_{r^*})S)\mathrm{diag}(
        \sqrt{\Sigma_1^2+\Sigma_2^2},~(S^{-1}~{\bf 0}_{r^* \times (r_{1*}+r_{2*}-r^*)})L^{-1}\mathrm{diag}(
              \Sigma_{1*},~\Sigma_{2*}))V^T.\end{align*}
    Therefore the SVD is derived by denoting
    $$
    (S^{-1}~{\bf 0}_{r^* \times (r_{1*}+r_{2*}-r^*)})L^{-1}\begin{pmatrix}
              \Sigma_{1*} & \\
            & \Sigma_{2*}
          \end{pmatrix}=U^*\Sigma^*V^{*T}.
    $$
\end{proof}

\begin{proof}[PROOF OF THEOREM \ref{nonorthobounds}]
    The proof of the current theorem with the non-orthogonal unshared vectors is under the same structure with and very similar to the proof of Theorem \ref{generalbounds}, so we only provide the sketch of the proof here. Also, we only provide when $(U_{1*}~U_{2*})$ is singular, it is similar with the non-singular situation except with extra zero block matrices. For the upper bound result, we can directly have the result by similar technique. For the lower bound, we still first prove a result similar to Proposition \ref{generalizacainzhanglower}, where we correspondingly construct something different. For the strong signal regime where $\gamma^2 \gtrsim p_{V1}\wedge p_{V2}$, we construct the following matrices in this case:
    $$U_r=(U_{0}^T~0_{r\times (r_{1*}+r_{2*})})^T,~U_0\in \mathbb{O}_{p_1-r_{1*}-r_{2*},r},$$
    $$
    U_{1*}=(0_{r\times p_1-r_{1*}-r_{2*}}~I_{r_{1*}}~0_{r\times r_{2*}})^T,~U_{2*}=(0_{r\times p_1-r_{1*}-r_{2*}}~0_{r\times r_{1*}}~I_{r_{2*}})^T,$$
    $$
    X_{U_r}=\begin{pmatrix}
        U_r&U_{1*}&U_{2*}
    \end{pmatrix}\mathrm{diag}(
        \sqrt{\alpha^2+\beta^2}I_r,~\Sigma_{1*},~\Sigma_{2*})\mathrm{diag}(
        I,~V^{*T})\begin{pmatrix}
        \frac{\alpha}{\sqrt{\alpha^2+\beta^2}}V_{1r}^T & \frac{\beta}{\sqrt{\alpha^2+\beta^2}}V_{2r}^T\\
        V_{1*}^T & 0\\ 0 & V_{2*}^T
    \end{pmatrix},$$
    where the added $V^*$ is fixed and satisfies certain conditions because it should be derived from the SVD of the corresponding previously defined $S^{-1}\mathrm{diag}(
        \Sigma_{1*},~\Sigma_{2*})$. The added matrix will not affect the later steps because of its orthonormality.
    
    In the weak signal regime, where $p_{V1}\asymp p_{V2}$ and  $\gamma^2 \lesssim p_{V1}\wedge p_{V2}$,
    We construct Gaussian mixture measure
    \begin{align*}&\Bar{P}_{U_r,\gamma}(Y)=C_{U_r,\gamma}\int_{W \in \mathbb{R}^{(p_{V1}-r_*) \times r}:\sigma_{min}(W) \geq \frac{1}{2\sqrt{2}}}\frac{1}{(2\pi\tau^2)^{p_1p_2/2}}\bigg(\frac{p_{V1}}{2\pi}\bigg)^{\frac{r(p_{V1}-r_{1*})}{2}}\\&\exp\bigg(-\frac{p_{V1}\|W\|_F^2}{2}\bigg)\exp\bigg(-\bigg\|Y-2\gamma \begin{pmatrix}
    U_r~U_{1*}\Sigma_{1*}~U_{2*}\Sigma_{2*}
    \end{pmatrix}\begin{pmatrix}
        I&\\&V^{*T}
    \end{pmatrix}\begin{pmatrix}
        W^T& 0_{r_{1*}}& W^T &0_d& 0_{r_{2*}}\\ 0 & I_{r_{1*}} & 0 &0& 0\\0&0&0&0& I_{r_{2*}}
    \end{pmatrix}\bigg\|_F^2/2\tau^2\bigg)dW.\end{align*}
Here we still require $V^*$ to satisfy certain conditions. The later steps of the proof is the same as the proof of Proposition \ref{generalizacainzhanglower} because the added term $\mathrm{diag}(
    I,~V^{*T})$ will become a part of previously defined $P(Y_{1*})$ and $P(Y_{2*})$ in the proof of Proposition \ref{generalizacainzhanglower} , which does not change the later KL-divergence calculation. We can then finish the proof by same procedure.
\end{proof}

\begin{proof}[PROOF OF THEOREM \ref{algothrm}]
    For matrices $Y_1$, $Y_2$ and $(Y_1~Y_2)$, we have  by Proposition \ref{generalcainzhang},
   $$
    \mathbb{E}\|\sin\Theta(\hat{u}_{1i},u_{1i})\|^2 \leq c\frac{n(\alpha^2+p_1)}{\alpha^4} \leq o\bigg(\frac{1}{\alpha^\epsilon}\bigg), $$
    $$
    \mathbb{E}\|\sin\Theta(\hat{u}_{2j},u_{2j})\|^2 \leq c\frac{n(\beta^2+p_1)}{\beta^4} \leq o\bigg(\frac{1}{\beta^\epsilon}\bigg),
    $$
   $$
    \mathbb{E}\|\sin\Theta(\hat{u}_{i},u_{i})\|^2 \leq c\frac{n(\gamma^2+p_1+p_2)}{\gamma^4} \leq o\bigg(\frac{1}{\gamma^\epsilon}\bigg),
    $$
    for any $\epsilon<2-c_1$.
    Define $a\asymp \frac{1}{\alpha^{\epsilon/2}}$, $b\asymp \frac{1}{\beta^{\epsilon/2}}$ and $d\asymp \frac{1}{\gamma^{\epsilon/2}}$. By Markov inequality, we have
    $$
    \mathbb{P}(\|\sin\Theta(\hat{u}_{1i},u_{1i})\|^2 \geq a) \leq \frac{\mathbb{E}\|\sin\Theta(\hat{u}_{1i},u_{1i})\|^2}{a} \leq o\bigg(\frac{1}{\alpha^{\epsilon/2}}\bigg).$$
    Similarly for $Y_2$ and $(Y_1~Y_2)$, we have
    $$
    \mathbb{P}(\|\sin\Theta(\hat{u}_{2j},u_{2j})\|^2 \geq b) \leq o\bigg(\frac{1}{\beta^{\epsilon/2}}\bigg),$$
    $$\mathbb{P}(\|\sin\Theta(\hat{u}_{i},u_{i})\|^2 \geq d) \leq o\bigg(\frac{1}{\gamma^{\epsilon/2}}\bigg).$$
    Consider the event 
    $$
    \mathscr{E}=\{\|\sin\Theta(\hat{u}_{1i},u_{1i})\|^2 < a; ~\|\sin\Theta(\hat{u}_{2j},u_{2j})\|^2 < b; ~\|\sin\Theta(\hat{u}_{i},u_{i})\|^2 < d\}.$$
    We have $$\mathbb{P}(\mathscr{E}) \geq 1-o\bigg(\frac{1}{\alpha^{\epsilon/2}}\bigg)-o\bigg(\frac{1}{\beta^{\epsilon/2}}\bigg)-o\bigg(\frac{1}{\gamma^{\epsilon/2}}\bigg).$$ We then discuss the behavior of $d_{1i}$ under event $\mathscr{E}$. When the underlying vector of the corresponding $\hat{u}_{1i}$ is shared. We temporarily denote the corresponding shared vector in $Y_2$ as $\hat{u}_{1i}^{(2)}$, then we have $$d_{1i}=\mathop{\min}\limits_{j=1,..,r_2}\|\sin\Theta(\hat{u}_{1i},\hat{u}_{2j})\|^2\leq \|\sin\Theta(\hat{u}_{1i},\hat{u}_{1i}^{(2)})\|^2\leq 2\|\sin\Theta(\hat{u}_{1i},u_{1i})\|^2+2\|\sin\Theta(\hat{u}_{1i}^{(2)},u_{1i})\|^2\leq 2(a+b) \to 0.$$ When the underlying vector of the corresponding $\hat{u}_{1i}$ is unshared. For any $j=1,..,r_2$, we have \begin{align*}
        \|\sin\Theta(\hat{u}_{1i},\hat{u}_{2j})\| \geq \|\sin\Theta(\hat{u}_{2j},u_{1i})\|-\|\sin\Theta(\hat{u}_{1i},u_{1i})\|\geq \|\sin\Theta(u_{1i},u_{2j})\|-\|\sin\Theta(u_{2j},\hat{u}_{2j})\|-\|\sin\Theta(\hat{u}_{1i},u_{1i})\|\\
        \geq 1-\sqrt{a}-\sqrt{b} \to 1.\end{align*}
        Combining results, $$d_{1i}=\mathop{\min}\limits_{j=1,..,r_2}\|\sin\Theta(\hat{u}_{1i},\hat{u}_{2j})\|^2 \geq 1+a+b-2\sqrt{a}-2\sqrt{b}+2\sqrt{ab} \to 1.$$ As a consequence, as long as $\alpha,\beta$ are sufficiently large, we can always have the correct singular vectors with indices in $I$, by selecting the vectors corresponding to the last $r_1-k_1$ of $d_{1i}$, because there are  $r_1-k_1$ shared vectors and their corresponding $d_{1i}$ converge to $0$ while the $d_{1i}$ for unshared vectors converge to $1$. Analogous to $X_1$, we can select correct singular vectors in $X_2$ with indices $J$. 
Next we can start tracing the vectors in the stacked matrix. By previous analysis, we know that, for orthogonal unshared vectors case, the left singular vectors for noiseless $(X_1~X_2)$ is the union of left singular vectors of $X_1$ and $X_2$. For $X_1$ and $Y_1$, suppose we traced $\hat{u}_{1i}$, $i \in I$ and the corresponding noiseless $u_{1i}$ is a shared vector. Still under event $\mathscr{E}$, among all vectors $\hat{u}_k$ in the left singular subspace of $(Y_1~Y_2)$, there always exist a $\hat{u}_s$ that is close to its corresponding noiseless vector. i.e. $$\|\sin\Theta(\hat{u}_s,u_{1i})\|^2=\|\sin\Theta(\hat{u}_s,u_{s})\|^2 \leq d.$$ so we have 
$$
\|\sin\Theta(\hat{u}_{1i},\hat{u}_s)\|^2\leq 2\|\sin\Theta(\hat{u}_{1i},u_{1i})\|^2+2\|\sin\Theta(u_{1i},\hat{u}_s)\|^2
\leq 2(a+d) \to 0.$$
For other top $r+k_1+k_2$ vectors in the stacked noisy matrix, that is when  $k \neq s$, we have
$$\|\sin\Theta(\hat{u}_{1i},\hat{u}_k)\|^2\geq (\|\sin\Theta(u_{k},u_{1i})\| - \|\sin\Theta(\hat{u}_{1i},u_{1i})\| - \|\sin\Theta(\hat{u}_{k},u_{k})\|)^2\geq (1-\sqrt{a}-\sqrt{d})^2 \to 1.
$$
Consequently, when $\alpha,\beta$ are large enough, it is guaranteed that 
\[
s=\mathop{\arg\min}\limits_{k \in 1,...,r+k_1+k_2}\|\sin\Theta(\hat{u}_{1i},\hat{u}_{k})\|^2
,\]
which means that we successfully trace the correct vectors.
Then, under event $\mathscr{E}$, we have $\widehat{\mathbb{J}}=\mathbb{J}$ and $$\mathbb{P}(\mathscr{E}) \geq 1-o\big(\frac{1}{\alpha^{\epsilon/2}}\big)-o\big(\frac{1}{\beta^{\epsilon/2}}\big)-o\big(\frac{1}{\gamma^{\epsilon/2}}\big).$$
\end{proof}

\begin{proof}[PROOF OF COROLLARY \ref{smalloptimal}]
    For the upper bound results, we can first derive the initial results  by Theorem 3 in \cite{caizhang2018} that
    $$
        \mathbb{E}\|\sin\Theta(U,\hat{U})\|^2 \leq \frac{cn\tau^2(\sigma_r^2(X_1)+\tau^2p_1)}{\sigma_r^4(X_1)}\wedge 1.
    $$
    Then the upper bound results hold by noting that this is a decreasing function of $\sigma_r(X_1)$ and 
    $$
    \gamma^2 \leq \mathop{\min}\limits_{1 \leq i \leq r}\{\sigma_{(i)}^2(X_1)+\sigma_{(i)}^2(X_2)\}\leq c\sigma_r^2(X_1).
    $$
For the minimax lower bound result, note that $p_1$ is of equal or lower order than $p_2$, the result holds by the same procedure as in the proof of Theorem \ref{upperlowerbound} and Proposition \ref{generalizecainzhanglowershare}, removing the final steps requiring the same order of $p_1$ and $p_2$.
\end{proof}

\begin{proof}[PROOF OF COROLLARY \ref{individualvectorrate}]
In this proof, for a matrix $U$, we use $U_{[i:j]}$ to denote the sub-matrix of $U$ formed by its $i$-th to $j$-th column. By the essentially the same procedure in Proposition \ref{generalcainzhang}, for the estimation of the top $i$ singular vectors, $$
    \mathbb{E}\|\sin(\hat{U}_{[1:i]},U_{[1:i]})\|^2\leq \frac{cn(g_i^2+p_1+p_2)}{g_i^4}.$$
For the estimation of the top $i-1$ singular vectors, $$
    \mathbb{E}\|\sin(\hat{U}_{[1:(i-1)]},U_{[1:(i-1)]})\|^2\leq \frac{cn(g_{i-1}^2+p_1+p_2)}{g_{i-1}^4}.$$
Combining Equations above, it holds that
$$
    \|\sin(\hat{u}_i,u_i)\|\leq \|\hat{u}_i\hat{u}_i^T-u_iu_i^T\|\leq \left\|\hat{U}_{[1:i]}\hat{U}_{[1:i]}^T-U_{[1:i]}U_{[1:i]}^T\right\| +\left\|U_{[1:(i-1)]}U_{[1:(i-1)]}^T-\hat{U}_{[1:(i-1)]}\hat{U}_{[1:(i-1)]}^T\right\|.$$
Square and taking expectation gives the final result
$$
\mathbb{E}\|\sin(\hat{u}_i,u_i)\|^2\leq \frac{cn(g_i^2+p_1+p_2)}{g_i^4}+\frac{cn(g_{i-1}^2+p_1+p_2)}{g_{i-1}^4}\leq \frac{cn(\min\{g_i^2,g_{i-1}^2\}+p_1+p_2)}{\min\{g_i^2,g_{i-1}^2\}^2}.$$
\end{proof}

\subsection{Proof of propositions}
\begin{proof}[PROOF OF PROPOSITION \ref{generalizecainzhanglowershare}]
    First, we consider the situation when $\gamma^2\gtrsim \tau^2(p_{V1}+p_{V2})$ holds. For $X \in \mathscr{H}_{r,\gamma}^{(s)}$, we construct the following matrices:
    $$V_1^T=(I_r,0_{p_{V1}-r}),~ V_2^T=(I_r,0_{p_{V2}-r}),$$
    $$\Sigma_1=\alpha I_r,~\Sigma_2=\beta I_r,$$
    $$
    X_{U}=U\sqrt{\alpha^2+\beta^2}I_r\bigg(\frac{\alpha}{\sqrt{\alpha^2+\beta^2}}I_r~0_{p_{V1}-r}~\frac{\beta}{\sqrt{\alpha^2+\beta^2}}I_r~0_{p_{V2}-r}\bigg).
    $$
    We denote $P_{U}(Y)$ as the probability measure of $Y=X_{U}+Z$,
    with the entries of $Z \in \mathbb{R}^{p_1 \times p_2}$ follow $N(0,\tau^2)$ independently and identically. By the KL-divergence between two $p$-dimensional multivariate Gaussians:
    $$
    D(N(\mu_0,\Sigma_0)||N(\mu_1,\Sigma_1))=\frac{1}{2}\bigg(tr(\Sigma_0^{-1}\Sigma_1)+(\mu_1-\mu_0)^T\Sigma_1^{-1}(\mu_1-\mu_0)-p+\log(\frac{\det\Sigma_1}{\det\Sigma_0})\bigg).$$
    Set $\alpha^2+\beta^2=\gamma^2$, we have the following expression of KL-divergence:
    \begin{align*}
    D(P_{U}||P_{U'})=\frac{1}{2\tau^2}\bigg\|(U-U')\sqrt{\alpha^2+\beta^2}I_r\bigg(\frac{\alpha}{\sqrt{\alpha^2+\beta^2}}I_r~0_{p_{V1}-r}~\frac{\beta}{\sqrt{\alpha^2+\beta^2}}I_r~0_{p_{V1}-r}\bigg)\bigg\|_F^2  \\=\frac{1}{2\tau^2}tr((\alpha^2+\beta^2)(U-U')(U-U')^T)=\frac{\gamma^2}{2\tau^2}
    \|U-U'\|_F^2.\end{align*}
    Then we construct a ball of radius $\epsilon$ centered at $U_{0} \in \mathbb{O}(p_1,r)$:
    $$B(U_{0},\epsilon)=\{U_{r}':\|\sin\Theta(U_{0},U')\|_F^2 \leq \epsilon\}
    .$$
    Based on Lemma 1 in \cite{Cai2013}, for any $0<a<1$ and $0<\epsilon<\sqrt{2r}$, there exists $\{U'_{1},...,U'_{m}\} \subseteq B(U_{0},\epsilon)$, such that:
    $$
    m\geq (\frac{c}{a})^{r(p_1-r)} \quad 
    \underset{1\leq i\neq j \leq m}{\min}\|\sin\Theta(U'_{i},U'_{j})\|_F \geq a\epsilon.$$
    By Lemma 1 in \cite{caizhang2018}, there exists $O_i \in \mathbb{O}(r,r)$ such that
    $$
    \|U_{0}-U'_{i}O_i\|_F \leq \sqrt{2}\|\sin\Theta(U_{0},U'_{i})\|_F \leq \sqrt{2}\epsilon.$$
    Denote $U_i=U'_{i}O_i \in \mathbb{O}(p_1,r)$, we have
    $$
    \underset{1\leq i\neq j \leq m}{\max}D(P_{U_{i}}||P_{U_{j}})=\underset{1\leq i\neq j \leq m}{\max}\frac{\gamma^2}{2\tau^2}\|U_{i}-U_{j}\|_F^2\leq \frac{\gamma^2}{2\tau^2}\underset{1\leq i\neq j \leq m}{\max}2(\|U_{0}-U_{i}\|_F^2+\|U_{0}-U_{j}\|_F^2) \leq \frac{4\gamma^2\epsilon^2}{\tau^2}.$$
    Then, by Fano's Lemma \citep{Yu1997}, we have 
    $$
    \underset{\hat{U}}{\inf}\underset{U \in \mathbb{O}(p_1,r)}{\sup}\mathbb{E}_{P_{U}}\|\sin\Theta(\hat{U},U)\|_F^2\geq a^2\epsilon^2\bigg(1-\frac{4\gamma^2\epsilon^2/\tau^2+\log(2)}{r(p_1-r)\log(c/a)}\bigg).$$
    We particularly select $$a=c\exp(-(1+\log2)),\qquad \epsilon=\sqrt{\frac{r(p_1-r)\tau^2}{4(1+\log2)\gamma^2}} \wedge \sqrt{2r}.$$ When $$\sqrt{\frac{r(p_1-r)\tau^2}{4(1+\log2)\gamma^2}} \leq \sqrt{2r},$$ we have:
    \begin{align*}
        \underset{\hat{U}}{\inf}\underset{U \in \mathbb{O}(p_1,r)}{\sup}\mathbb{E}_{P_{U}}\|\sin\Theta(\hat{U},U)\|_F^2\geq c^2\exp(-2(1+\log2))\frac{r(p_1-r)\tau^2}{4(1+\log2)\gamma^2}\bigg(1-\frac{r(p_1-r)/(1+\log2)+\log2}{r(p_1-r)(1+\log2)}\bigg)\\
        \stackrel{\text{(i)}}{\geq} \frac{c^2\exp(-2(1+\log2)\tau^2}{4(1+\log2)^2\gamma^2}\bigg((1+\log2-\frac{1}{1+\log2})\frac{p_1r}{2}-\log2\bigg) \stackrel{\text{(ii)}}{\geq} c\frac{p_1r\tau^2}{\gamma^2},
    \end{align*}
    where (i) is because of the assumption $p_1 \geq 2r$ and (ii) is because  $$(1+\log2-\frac{1}{1+\log2})\frac{p_1r}{2}-\log2 \geq \frac{p_1r}{c_1}$$ for large enough constant $c_1>0$.
    Considering $\gamma^2 \geq C\tau^2(p_{V1}+p_{V2})$, 
    \begin{align}\label{diffstep1}
     \frac{p_1r\tau^2}{\gamma^2} 
    \geq c\frac{p_1r\tau^2(\gamma^2+\tau^2(p_{V1}+p_{V2}))}{\gamma^4}.
    \end{align}
    As a result, we have $$
    \underset{\hat{U}}{\inf}\underset{U \in \mathbb{O}(p_1,r)}{\sup}\mathbb{E}_{P_{U}}\|\sin\Theta(\hat{U},U)\|_F^2\geq c\bigg(\frac{p_1r\tau^2(\gamma^2+\tau^2(p_{V1}+p_{V2}))}{\gamma^4}
    \wedge r\bigg).$$
    By Lemma 6 in \cite{caizhang2018} and taking $\tau=1$ without loss of generality, we further have
    $$
    \underset{\hat{U}}{\inf}\underset{X \in \mathscr{H}_{r,\gamma}^{(s)}}{\sup}\mathbb{E}\|\sin\Theta(\hat{U},U)\|_F^2 \geq c\bigg(\frac{p_1r(\gamma^2+(p_{V1}+p_{V2}))}{\gamma^4}
    \wedge r\bigg).$$

    Next we consider the circumstance when $p_{V1}\asymp p_{V2}$ holds. In this case, when $\gamma^2 \geq \frac{\tau^2}{16}p_{V1}\wedge p_{V2}$, the proof can be finished as the same procedure in the previous case. We only need to substitute the Equation (\ref{diffstep1}) with
    $$
    \frac{p_1r\tau^2}{\gamma^2} \geq \frac{p_1r\tau^2(\gamma^2/2+\tau^2(p_{V1}\wedge p_{V2})/32)}{\gamma^4} \geq c\frac{p_1r\tau^2(\gamma^2+\tau^2(p_{V1}+p_{V2}))}{\gamma^4},$$
    considering $p_{V1}\asymp p_{V2}$  and $\gamma^2 \geq \frac{\tau^2}{16}p_{V1}\wedge p_{V2}$.
 Then we consider the case when $p_{V1}\asymp p_{V2}$ holds and $\gamma^2 < \frac{\tau^2}{16}p_{V1}\wedge p_{V2}$. Define the following class of density $P_Y$ : 
   $$
    \mathscr{P}_{U,\gamma}=\{P_Y:~Y=X+Z \in \mathbb{R}^{p_1\times p_2},~Z\sim N(0,\tau^2),X \in \mathscr{H}_{r,\gamma}^{(s)}, \text{left singular subspace of X is}~UO, O \in \mathbb{O}(r,r)\}.$$
     We provide the proof when $p_{V1}\leq p_{V2}$, and the proof when $p_{V2}\leq p_{V1}$ is symmetric. We construct the following Gaussian mixture measure
    \begin{align*}
    \Bar{P}_{U,\gamma}(Y)=&C_{U,\gamma}\int_{W \in \mathbb{R}^{p_{V1} \times r}:\sigma_{min}(W) \geq \frac{1}{2\sqrt{2}}}\frac{1}{(2\pi\tau^2)^{p_1p_2/2}}\bigg(\frac{p_{V1}}{2\pi}\bigg)^{\frac{p_{V1}r}{2}}\\
    &\exp\bigg(-\frac{\|Y-2\gamma U(W^T~W^T~{\bf 0}_{r \times (p_{V2}-p_{V1})})\|_F^2}{2\tau^2}\bigg)\exp\bigg(-\frac{p_{V1}\|W\|_F^2}{2}\bigg)dW,\end{align*}
    where $C_{U,\gamma}$ is a constant to make it a valid probability measure. It is easy to verify that
    $$2\gamma U(W^T~W^T~{\bf 0}_{r \times (p_{V2}-p_{V1})})=UO2\gamma \sqrt{\Sigma^2+\Sigma^2}\bigg(\frac{2\gamma\Sigma}{\sqrt{4\gamma^2\Sigma+4\gamma^2\Sigma}}V^T~\frac{2\gamma\Sigma}{\sqrt{4\gamma^2\Sigma+4\gamma^2\Sigma}}(V^T~{\bf 0}_{r \times (p_{V2}-p_{V1})})\bigg).$$
    Under the event $\{\sigma_{\min}(W) \geq \frac{1}{2\sqrt{2}}\}$, the minimum singular value corresponding to one of $U$ is greater than $\gamma$, which means that $\Bar{P}_{U,\gamma}$ is indeed a mixture density in $\mathscr{P}_{U,\gamma}$.
    For bounding KL-divergence, we have the following.
    \begin{lemma}\label{shareboundKL}
        Under the assumption of the theorem and $\gamma^2 \leq \frac{\tau^2p_{V1}}{16}$, for any $U,U' \in \mathbb{O}(p_1, r)$, we have $$
        D(\Bar{P}_{U,\gamma}||\Bar{P}_{U',\gamma})\leq \frac{p_{V1}128\gamma^4}{(8\gamma^2+\tau^2p_{V1})(\tau^2p_{V1}-8\gamma^2)}\|\sin\Theta(U,U')\|_F^2+C_{KL},$$
        where $C_{KL}$ is a uniform constant.
    \end{lemma}
    Also considering when $\gamma^2 \leq \frac{\tau^2p_{V1}}{16}$, we have $$\tau^2p_{V1}-8\gamma^2 \geq \frac{\tau^2p_{V1}}{2}.$$ Further we have
    $$
    D(\Bar{P}_{U,\gamma}||\Bar{P}_{U',\gamma}) \leq \frac{256\gamma^4}{\tau^2(8\gamma^2+\tau^2p_{V1})}\|\sin\Theta(U,U')\|_F^2+C_{KL}.
    $$
    We consider a ball of radius $0<\epsilon<\sqrt{2r}$ centered at $U_{0} \in \mathbb{O}(p_1,r)$:
$$B(U_{0},\epsilon)=\{U':\|\sin\Theta(U_{0},U')\|_F^2 \leq \epsilon\}.$$
    Again based on Lemma 1 in \cite{Cai2013}, for any $0<a<1$ and $0<\epsilon<\sqrt{2r}$, there exists $\{U_{1},...,U_{m}\} \subseteq B(U_{0},\epsilon)$, such that $$
    m\geq (\frac{c}{a})^{r(p_1-r)},\qquad
    \underset{1\leq i\neq j \leq m}{\min}\|\sin\Theta(U_{i},U_{j})\|_F \geq a\epsilon.$$ 
    So we have $$
        D(\Bar{P}_{U,\gamma}||\Bar{P}_{U',\gamma}) \leq \frac{256\gamma^4}{\tau^2(8\gamma^2+\tau^2p_{V1})}4\epsilon^2+C_{KL}.$$ 
    Then, by Fano's Lemma \citep{Yu1997}, we have 
    $$\underset{\hat{U}}{\inf}\underset{U \in \mathbb{O}(p_1,r)}{\sup}\mathbb{E}_{\hat{P}_{U,\gamma}}\|\sin\Theta(\hat{U},U)\|_F^2 \geq a^2\epsilon^2\bigg(1-\frac{\frac{1024\gamma^4}{\tau^2(8\gamma^2+\tau^2p_{V1})}\epsilon^2+C_{KL}+log(2)}{r(p_1-r)log(c/a)}\bigg).$$
    We particularly select $$a=c\exp(-(1+C_{KL}+log2))$$ and we choose $$\epsilon=\sqrt{\frac{r(p_1-r)\tau^2(8\gamma^2+\tau^2p_{V1})}{2\cdot1024\gamma^4}} \wedge \sqrt{2r},$$ when it holds that $$\sqrt{\frac{r(p_1-r)\tau^2(8\gamma^2+\tau^2p_{V1})}{2\cdot1024\gamma^4}} \leq \sqrt{2r},$$ we have the following inequality:
    \begin{align*}
    \underset{\hat{U}}{\inf}&\underset{U \in \mathbb{O}(p_1,r)}{\sup}\mathbb{E}_{\hat{P}_{U,\gamma}}\|\sin\Theta(\hat{U},U)\|_F^2 \geq c^2\exp(-2(1+C_{KL}+\log2))\frac{r(p_1-r)(\tau^2(4\gamma^2+\tau^2p_{V1}))}{2\cdot 1024\gamma^4}\\
    &\bigg(r(p_1-r)(1+C_{KL}+\log2)-r(p_1-r)/2-(C_{KL}+\log2)\bigg)/\bigg(r(p_1-r)(1+C_{KL}+\log2)\bigg)\geq c\frac{p_1r\tau^2(\gamma^2+\tau^2p_{V1})}{\gamma^4},
    \end{align*}
    where the last inequality holds because when $p_1\geq 2r,$ $$\frac{r(p_1-r)}{2}+(r(p_1-r)-1)(C_{KL}+\log2) \geq \frac{p_1r}{8}.$$
    With Lemma 6 in \cite{caizhang2018} and taking $\tau=1$ without loss of generality, also noting the same order of $p_{V1}$ and $p_{V2}$, we have
    $$
    \underset{\hat{U}}{\inf}\underset{X \in \mathscr{H}_{r,\gamma}^{(s)}}{\sup}\mathbb{E}_{\hat{P}_{U,\gamma}}\|\sin\Theta(\hat{U},U)\|_F^2\geq \underset{\hat{U}}{\inf}\underset{U \in \mathbb{O}(p_1,r)}{\sup}\mathbb{E}_{\hat{P}_{U,\gamma}}\|\sin\Theta(\hat{U},U)\|_F^2\geq c\bigg(\frac{p_1r(\gamma^2+(p_{V1}+p_{V2}))}{\gamma^4} \wedge r\bigg).$$
    Finally, under all circumstances, the theorem for the Frobenius norm holds. For the spectral norm, it holds simply by $$\|\sin\Theta(\hat{U},U)\| \geq \frac{1}{\sqrt{r}} \|\sin\Theta(\hat{U},U)\|_F.$$
\end{proof}

\begin{proof}[PROOF OF PROPOSITION \ref{vectorrelation}]
    This proposition can be easily shown by analyzing SVD of the stacked noiseless matrix. (i) When all vectors in $U_1$ and $U_2$ are different:
    The SVD is
$$
    (X_1~X_2)=(
        U_1~U_2)\mathrm{diag}(
        \Sigma_1,~\Sigma_2)\mathrm{diag}(
        V^T_1,~V^T_2).$$ 
    Here, after sorting singular values in decreasing order, the vectors in $U_1$ and $U_2$ may come into shuffled positions. (ii) When there are some shared vectors between $U_1$ and $U_2$: 
    Suppose
$$
    X_1=(U_r~U_{1*})\mathrm{diag}(
        \Sigma_{1r},~ \Sigma_{1*})V_1^T,$$ and  $$X_2=(U_r~U_{2*})\mathrm{diag}(
        \Sigma_{2r},~\Sigma_{2*})V_2^T, $$
    where $U_r$ contains shared vectors while $U_{1*}$ and $U_{2*}$ are unshared. Thus by similar decomposition as our previous analysis
   $$
    (X_1~X_2)=(U_r~U_{1*}~U_{2*})\mathrm{diag}(
        \sqrt{\Sigma_{1r}^2+\Sigma_{2r}^2} ,~\Sigma_{1*},~\Sigma_{2*})V_*^T.$$
Therefore, the singular vectors of individual noiseless matrices $X_1$, $X_2$ can all be found in the stacked noiseless matrix $(X_1~X_2)$. This completes the proof of Proposition \ref{vectorrelation}.
\end{proof}

\begin{proof}[PROOF OF PROPOSITION \ref{generalcainzhang}]

This proposition is a generalization of the Theorem 3 in \cite{caizhang2018}. Therefore, we follow their structure of proof. We only need to prove the right singular subspace statement as the proof for the left singular subspace is symmetric. Also only the proof for spectral norm is necessary by its relationship with corresponding Frobenius norm. For brevity, we denote the top $r$ right singular subspace as $V_r$ and the remaining vectors as $V_*$.  We have the lemma.
 \begin{lemma}\label{lemmageneralize}
Suppose $X \in \mathbb{R}^{p_1\times p_2}$ has rank greater than $r$ with right singular subspace $V^T=(
    V_r ~V_*)^T$. $Y=X+Z$ with $Z_{ij}\mathop{\sim}\limits^{iid} \mathscr{G}_\tau  $. Then there exists constants $C, c>0$ such that for any $x>0$
$$\mathbb{P}\left(\sigma_r^2(YV_r)\leq(\sigma_r^2(X)+p_1\tau^2)(1-x)\right)\leq C\exp\left\{Cr-c\left(\sigma_r^2(X)+p_1\tau^2\right)x^2\wedge x\right\},$$ and also
$$
\mathbb{P}\left(\sigma_{r+1}^2(Y)-\sigma_{r+1}^2(X)\geq p_1\tau^2(1+x)\right)\leq C\exp\left\{Cp_{2}-c (p_1x^2) \wedge (p_1x) \wedge (x^2\frac{p_1^2\tau^2}{\sigma_r^2(X)})\right\}.$$
Moreover, denote columns orthogonal to $V$ as $V_\perp$ and denote $(V_* \quad V_\perp )$ as $V_{r\perp}$, there exists $C_{\mathrm{gap}}, C, c$, such that whenever
$\sigma _{r}^{2}(X) \geq C_{\mathrm{gap}}p_{2}$, for any $x> 0$ we have
\begin{align}\label{randommatrix3}
\nonumber&\mathbb{P}\left(\|\mathbb{P}_{YV_r}YV_{r\perp}\|\leq x\right)\geq1-C\exp\left\{-c(\sigma_r^2(X)+p_1\tau^2)\right\}\\
&\quad-C\exp\left\{Cp_2-c\min(x^2,\sqrt{\sigma_r^2(X)+p_1\tau^2}x)\right\}
\end{align}
\end{lemma}

We first focus on the scenario that $\sigma_r^2(X)\geq C_{\mathrm{gap}}((p_1p_2)^{\frac12}+p_2)$ for some large constant $C_\mathrm{gap}>0$. By Lemma \ref{lemmageneralize}, for large constant $k>2$, set $$x=\frac{\sigma_r^2(X)}{k(\sigma_r^2(X)+p_1\tau^2)},$$ then there exists constants $C,c$ such that
$$\mathbb{P}\left(\sigma_r^2(YV_r)\leq\sigma_r^2(X)+p_1\tau^2-\frac{\sigma_r^2(X)}{k}\right)\leq C\exp\left\{Cr-c\min\left(k\sigma_r^2(X),\frac{\sigma_r^4(X)}{\sigma_r^2(X)+p_1\tau^2}\right)\right\}.$$
For large constant $k$, set $$x=\frac{\sigma_r^2(X)}{kp_1\tau^2},$$
then
$$
\mathbb{P}\left(\sigma_{r+1}^2(Y)\geq \sigma_{r+1}^2(X)+p_1\tau^2+\dfrac{1}{k}\sigma_r^2(X)\right)\leq C\exp\left\{Cp_2-c\min\left(\frac{\sigma_r^4(X)}{k^2p_1\tau^2},\frac{\sigma_r^2(X)}{k},\frac{\sigma_r^2(X)}{k^2}\right)\right\},$$ for constant $c$ depending on sub-Gaussian parameter $\delta$.
Also
$$
\mathbb{P}\left\{\|\mathbb{P}_{YV_r}YV_{r\perp}\|\geq x\right\}\leq C\exp\left\{Cp_2-c\min\left(x^2,x\sqrt{\sigma_r^2(X)+p_1\tau^2}\right)\right\}+C\exp\left\{-c(\sigma_r^2(X)+p_1\tau^2)\right\}.$$
When $C_{\mathrm{gap}}$ is large enough, it holds that $$\sigma_r^{4}(X)/(\sigma_{r}^{2}(X)+p_{1}\tau^2)\geq Cp_{2}.$$ Then
$$
c\min\left(\frac{\sigma_r^4(X)}{\sigma_r^2(X)+p_1\tau^2},k\sigma_r^2(X)\right)-Cr\geq c\frac{\sigma_r^4(X)}{\sigma_r^2(X)+p_1\tau^2}-Cr\geq c\frac{\sigma_r^4(X)}{\sigma_r^2(X)+p_1\tau^2},$$ and 
$$c\min\left(\frac{\sigma_r^4(X)}{k^2p_1\tau^2},\frac{\sigma_r^2(X)}{k},\frac{\sigma_r^2(X)}{k^2}\right)-Cp_2\geq c\frac{\sigma_r^4(X)}{\sigma_r^2(X)+p_1\tau^2}-Cp_2\geq c\frac{\sigma_r^4(x)}{\sigma_r^2(X)+p_1\tau^2},$$ taking $\tau=1$ without loss of generality. Also, we have
$$c\operatorname*{min}(x^{2},\sqrt{\sigma_{r}^{2}(X)+p_{1}\tau^2}x)-Cp_{2}=c(\sigma_{r}^{2}(X)+p_{1}\tau^2)-Cp_{2}\geq c\left(\sigma_{r}^{2}(X)+p_{1}\tau^2\right) \geq c\frac{\sigma_{r}^{4}(X)}{\sigma_{r}^{2}(X)+p_{1}\tau^2},\quad\mathrm{if}\:x=\sqrt{\sigma_{r}^{2}(X)+p_{1}\tau^2}.$$
To sum up, if we denote the event $Q$ as
$$
    Q=\Big\{\sigma_{r}^{2}(YV_r)\geq\sigma_{r}^{2}(X)+p_{1}\tau^2-\frac{\sigma_{r}^{2}(X)}{k}, \sigma_{r+1}^{2}(Y)\leq \sigma_{r+1}^{2}(X)+p_{1}\tau^2+\frac{1}{k}\sigma_{r}^{2}(X),  |\mathbb{P}_{YV_r}YV_{r\perp}\|\leq\sqrt{\sigma_{r}^{2}(X)+p_{1}\tau^2}\Big\},$$
when $\sigma_r^{2}(X)\geq C_{\mathrm{gap}}((p_{1}p_{2})^{\frac{1}{2}}+p_{2}),$ for some large constant $C_\mathrm{gap}>0$,
$$\mathbb{P}\left(Q^c\right)\leq C\exp\left\{-c\frac{\sigma_r^4(X)}{\sigma_r^2(X)+p_1\tau^2}\right\}.$$
Under event $Q$,we can apply Proposition 1 in \cite{caizhang2018}
$$\|\sin\Theta(\hat{V},V_r)\|^2\leq\frac{\sigma_r^2(YV_r)\|\mathbb{P}_{YV_r}YV_{r\perp}\|^2}{(\sigma_r^2(YV_r)-\sigma_{r+1}^2(Y))^2}.$$
Note that $\frac{x^{2}}{(x^{2}-y^{2})^{2}}$ is a decreasing function of $x$ and increasing function of $y$ when $x>y\geq0$. Here we use the even $Q$ that
$$
\sigma_{r}^{2}(YV_r)\geq\sigma_{r}^{2}(X)+p_{1}\tau^2-\frac{\sigma_{r}^{2}(X)}{k}$$, $$\sigma_{r+1}^{2}(Y)\leq \sigma_{r+1}^{2}(X)+p_{1}\tau^2+\frac{1}{k}\sigma_{r}^{2}(X).$$
After using the above two inequalities, and also noting that the gap condition $$\sigma_r^2(X)>(1+\epsilon)\sigma_{r+1}^2(X)$$ is equivalent to $$\sigma_{r+1}^2(X)<c\sigma_r^2(X)$$ for some $c\in (0,1)$, we have for the denominator, it holds that:
$$
     \sigma_{r}^{2}(X)+p_{1}\tau^2-\frac{\sigma_{r}^{2}(X)}{k} - \sigma_{r+1}^{2}(X)-p_{1}\tau^2-\frac{1}{k}\sigma_{r}^{2}(X) \geq \frac{k-2}{k} \sigma_r^2(X)-\sigma_{r+1}^2(X)\geq c\sigma_r^2(X).$$
The last step is guaranteed by assumption and previous choice of large enough $k$, which is related to the constant $\epsilon$ in the eigen-gap condition.
We have 
$$
\|\sin\Theta(\hat{V},V_r)\|^2 \leq C\frac{(\sigma_r^2(X)+p_1\tau^2)\|\mathbb{P}_{YV_r}YV_{r\perp}\|^2}{\sigma_r^4(X)}.$$
Next we  note that $\|\sin\Theta ( \hat{V} , V) \| \leq 1$, for any $\hat{V},V_r\in\mathbb{O}_{p_{2},r}.$ Therefore $$\mathbb{E}\left\|\sin\Theta(\hat{V},V_r)\right\|^{2}=\mathbb{E}\left\|\sin\Theta(\hat{V},V_r)\right\|^{2}1_{Q}+\mathbb{E}\left\|\sin\Theta(\hat{V},V_r)\right\|^{2}1_{Q^{c}}\leq\frac{C(\sigma_{r}^{2}(X)+p_{1}\tau^2)}{\sigma_{r}^{4}(X)}\mathbb{E}\|\mathbb{P}_{YV_r}YV_{r\perp}\|^{2}1_{Q}+\mathbb{P}(Q^{c}).$$
By basic property of exponential function $$\mathbb{P}(Q^c) \leq C\exp\left\{-c\frac{\sigma_r^4(X)}{\sigma_r^2(X)+p_1\tau^2}\right\}\leq C\frac{\sigma_r^2(X)+p_1\tau^2}{\sigma_r^4(X)}\leq C\frac{p_2(\sigma_r^2(X)+p_1\tau^2)}{\sigma_r^4(X)}.$$ It remains to consider $\mathbb{E}\left\|\mathbb{P}_{YV_r}YV_{r\perp}\right\|^{2}1_{Q}.$ Denote $T=\|\mathbb{P}_{YV_r}YV_{r\perp}\|.$ Applying Lemma \ref{lemmageneralize} again, we have for some constant $C_x$ that
\begin{align*}
\mathbb{E}T^{2}1_{Q}&\leq\mathbb{E}T^{2}1_{\{T^{2}\leq\sigma_{r}^{2}(X)+p_{1}\tau^2\}}\leq C_{x}p_{2}+\int_{C_{x}p_{2}}^{\sigma_{r}^{2}(X)+p_{1}\tau^2}\mathbb{P}\left(T^{2}1_{\{T^{2}\leq\sigma_{r}^{2}(X)+p_{1}\tau^2\}}\geq t\right)dt\\
&\stackrel{(\ref{randommatrix3})}{\leq} C_{x}p_{2}+\int_{C_{x}p_{2}}^{\sigma_{r}^{2}(X)+p_{1}\tau^2}C\left(\exp\left\{C_{2}-ct\right\}dt+\exp\left\{-c(\sigma_{r}^{2}(X)+p_{1}\tau^2)\right\}\right)dt\\
&\leq C_{x}p_{2}+C(\sigma_{r}^{2}(X)+p_{1}\tau^2)\exp\left\{-c(\sigma_{r}^{2}(X)+p_{1}\tau^2)\right\} +C\exp(Cp_{2})\cdot\exp\left(-cC_{x}p_{2}\right)\frac{1}{c}\\
&\leq C_{x}p_{2}+C+\frac{C}{c}\exp((C-cC_{x})p_{2}).\end{align*}
We can choose $C_x$ large enough and only relying on constants $C,c$ in the inequalities above, to ensure that $$\mathbb{E}T^21_Q\leq Cp_2,$$
for large constant $C>0$ as long as $p_{2}\geq1.$ Combining inequalities above, we obtain
$$\mathbb{E}\left\|\sin\Theta(\hat{V},V_r)\right\|^2\leq\frac{Cp_2(\sigma_r^2(X)+p_1\tau^2)}{\sigma_r^4(X)}\wedge1,$$
as long as $$\sigma_r^2(X)\geq C_{\mathrm{gap}}\left((p_1p_2)^{\frac12}+p_2\right)$$ for some large enough $C_\mathrm{gap}.$
Finally when $$\sigma_r^2(X)<\dot{C}_{\mathrm{gap}}\big((p_1p_2)^{{\stackrel{1}{2}}}+p_2\big),$$ we have
$$
    \frac{p_{2}(\sigma_{r}^{2}(X)+p_{1}\tau^2)}{\sigma_{r}^{4}(X)}\geq\frac{p_{2}(p_{1}\tau^2+C_{\mathrm{gap}}(p_{1}p_{2})^{\frac{1}{2}}+C_{\mathrm{gap}}p_{2})}{C_{\mathrm{gap}}^{2}(p_{1}p_{2}+p_{2}^{2}+2p_{1}^{\frac{1}{2}}p_{2}^{3/2})}=\frac{p_{1}+C_{\mathrm{gap}}(p_{1}p_{2})^{\frac{1}{2}}+C_{\mathrm{gap}}p_{2}}{C_{\mathrm{gap}}^{2}\left(p_{1}\tau^2+2(p_{1}p_{2})^{\frac{1}{2}}+p_{2}\right)}\geq\frac{c}{C_{\mathrm{gap}}},$$
for some constant $c$ that is related to $\delta$. Then
$$\mathbb{E}\left\|\sin\Theta(\hat{V},V_r)\right\|^2\leq1\leq\frac{Cp_2(\sigma_r^2(X)+p_1\tau^2)}{\sigma_r^4(X)}\wedge1,$$
when $C\geq \frac{C_{gap}}{c}$. We the finish the proof by taking $\tau=1$ without loss of generality.
\end{proof}

\begin{proof}[PROOF OF PROPOSITION \ref{generalizacainzhanglower}]
    For brevity, we denote the minimum singular value corresponding to vectors in $U_r$ as $\gamma$. First consider $\gamma^2\gtrsim \tau^2(p_{V1}+p_{V2})$. For $X \in \mathscr{H}_{r,\gamma}^{(u)}$, construct:
    $$U_r=(U_{0}^T~{\bf 0}_{r\times (r_{1*}+r_{2*})})^T,~~U_0\in \mathbb{O}_{p_1-r_{1*}-r_{2*},r},$$
    $$
    U_{1*}=({\bf 0}_{r\times p_1-r_{1*}-r_{2*}}~I_{r_{1*}}~{\bf 0}_{r\times r_{2*}})^T,~U_{2*}=({\bf 0}_{r\times p_1-r_{1*}-r_{2*}}~{\bf 0}_{r\times r_{1*}}~I_{r_{2*}})^T,$$
    $$
    X_{U_r}=\begin{pmatrix}
        U_r&U_{1*}&U_{2*}
    \end{pmatrix}\mathrm{diag}(
        \sqrt{\alpha^2+\beta^2}I_r,\Sigma_{1*},\Sigma_{2*})\begin{pmatrix}
        \frac{\alpha}{\sqrt{\alpha^2+\beta^2}}V_{1r} &V_{1*} & 0\\
         \frac{\beta}{\sqrt{\alpha^2+\beta^2}}V_{2r}& 0 & V_{2*}
    \end{pmatrix}^T,$$
    where orthonormal $V_{1r},V_{2r},V_{1*},V_{2*}$ are orthogonal to each other.  $\Sigma_{1*}$ and $\Sigma_{2*}$  can be selected to satisfy the gap condition. (e.g. select $\frac{\sqrt{\alpha^2+\beta^2}}{2}I$).
    Denote $P_{U_r}(Y)$ as the probability measure of $Y=X_{U_r}+Z$,
    with the entries of $Z \in \mathbb{R}^{p_1 \times p_2}\stackrel{iid}{\sim}N(0,\tau^2)$. By Gaussian KL:
    $$
    D(N(\mu_0,\Sigma_0)||N(\mu_1,\Sigma_1))=\frac{1}{2}\bigg(tr(\Sigma_0^{-1}\Sigma_1)+(\mu_1-\mu_0)^T\Sigma_1^{-1}(\mu_1-\mu_0)-p+\log(\frac{\det\Sigma_1}{\det\Sigma_0})\bigg).
    $$
    Set $\alpha^2+\beta^2=\gamma^2$, we have:
    \begin{align*}D(P_{U_r}||P_{U_r'})=\frac{1}{2\tau^2}\bigg\|\begin{pmatrix}
        U_r-U_r'&{\bf 0}_{1*}&{\bf 0}_{2*}
    \end{pmatrix}
    \begin{pmatrix}
        \sqrt{\alpha^2+\beta^2}I_r &&\\&\Sigma_{1*}&\\&&\Sigma_{2*}
    \end{pmatrix}\begin{pmatrix}
        \frac{\alpha}{\sqrt{\alpha^2+\beta^2}}V_{1r}^T & \frac{\beta}{\sqrt{\alpha^2+\beta^2}}V_{2r}^T\\
        V_{1*}^T & 0\\ 0 & V_{2*}^T
    \end{pmatrix}\bigg\|_F^2 
     \\=\frac{1}{2\tau^2}tr((\alpha^2+\beta^2)(U_r-U_r')(U_r-U_r')^T)=\frac{\gamma^2}{2\tau^2}
    \|U_r-U_r'\|_F^2,\end{align*}
    where ${\bf 0}_{i*}$ denotes the zero block matrix with the dimension of ${\bf 0}_{p_1\times r_{i*}}$. Construct a  radius $\epsilon$ ball centered at $U_{00} \in \mathbb{O}(p_1-r_{1*}-r_{2*} ,r)$:
$$B(U_{00},\epsilon)=\{U_{0}':\|\sin\Theta(U_{00},U_0')\|_F^2 \leq \epsilon\}.$$
    By Lemma 1 in \cite{Cai2013}, for any $0<a<1$ and $0<\epsilon<\sqrt{2r}$, there exists $\{U'_{01},...,U'_{0m}\} \subseteq B(U_{00},\epsilon)$, such that:
$$
    m\geq (\frac{c}{a})^{r(p_1-r_{1*}-r_{2*}-r)},$$ and $$\underset{1\leq i\neq j \leq m}{\min}\|\sin\Theta(U'_{0i},U'_{0j})\|_F \geq a\epsilon.$$
    Also recall that how we construct $U_r$. We have
    $$
    \|\sin\Theta(U_{r1},U_{r2})\|_F=\|\sin\Theta((U_{01}^T~{\bf 0}_{r\times (r_{1*}+r_{2*})})^T,(U_{02}^T~{\bf 0}_{r\times (r_{1*}+r_{2*})})^T)\|_F=\|\sin\Theta(U_{01},U_{02})\|_F.
    $$
    Thus Construct radius $\epsilon$ ball centered at $U_{r0} \in \mathbb{O}(p_1, r)$:
    \[
    B(U_{r0},\epsilon)=\{U_{r}':\|\sin\Theta(U_{r0},U_r')\|_F^2 \leq \epsilon\}.
    \]
    Also for any $0<a<1$ and $0<\epsilon<\sqrt{2r}$, there exists $\{U'_{r1},...,U'_{rm}\} \subseteq B(U_{r0},\epsilon)$, where $U_{ri}=(U_{0i}^T~{\bf 0}_{r\times (r_{1*}+r_{2*})})^T$ , such that: 
    \[
    m\geq (\frac{c}{a})^{r(p_1-r_{1*}-r_{2*}-r)},\] and \[
    \underset{1\leq i\neq j \leq m}{\min}\|\sin\Theta(U'_{ri},U'_{rj})\|_F \geq a\epsilon.
    \]
    By Lemma 1 in \cite{caizhang2018}, there exists $O_i \in \mathbb{O}(r,r)$ such that
    \[
    \|U_{00}-U'_{0i}O_i\|_F \leq \sqrt{2}\|\sin\Theta(U_{00},U'_{0i})\|_F \leq \sqrt{2}\epsilon.
    \]
    Denote $U_{0i}=U'_{0i}O_i \in \mathbb{O}(p_1,r)$, $U_{ri}$ are similarly constructed with $U_{0i}$. we have
    $$\underset{1\leq i\neq j \leq m}{\max}D(P_{U_{ri}}||P_{U_{rj}})=\underset{1\leq i\neq j \leq m}{\max}\frac{\gamma^2}{2\tau^2}\|U_{ri}-U_{rj}\|_F^2\leq \frac{\gamma^2}{2\tau^2}\underset{1\leq i\neq j \leq m}{\max}2(\|U_{r0}-U_{ri}\|_F^2+\|U_{r0}-U_{rj}\|_F^2) \leq \frac{4\gamma^2\epsilon^2}{\tau^2}.$$
    Define $\mathbb{O}$ as the class of $U_r$ constructed by the above means.
    Then, by Fano's Lemma \citep{Yu1997}, we have 
    $$
    \underset{\hat{U}_r}{\inf}\underset{U_r \in \mathbb{O}}{\sup}\mathbb{E}_{P_{U_r}}\|\sin\Theta(\hat{U}_r,U_r)\|_F^2\geq a^2\epsilon^2\bigg(1-\frac{4\gamma^2\epsilon^2/\tau^2+\log(2)}{r(p_1-r_{1*}-r_{2*}-r)\log(c/a)}\bigg).
    $$
    Select $a=c\exp(-(1+\log2))$ and $$\epsilon=\sqrt{\frac{r(p_1-r_{1*}-r_{2*}-r))\tau^2}{4(1+log2)\gamma^2}} \wedge \sqrt{2r},$$ when $$\sqrt{\frac{r(p_1-r_{1*}-r_{2*}-r)\tau^2}{4(1+\log2)\gamma^2}} \leq \sqrt{2r}.$$ We have:
    \begin{align*}
        \underset{\hat{U_r}}{\inf}\underset{U_r \in \mathbb{O}}{\sup}\mathbb{E}_{P_{U_r}}\|\sin\Theta(\hat{U}_r,U_r)\|_F^2 \geq \frac{c^2\exp(-2(1+\log2)\tau^2}{4(1+\log2)^2\gamma^2}\cdot\bigg((1+\log2-\frac{1}{1+\log2})r(p_1-r_{1*}-r_{2*}-r)-\log2\bigg)\\
        \stackrel{\text{(i)}}{\geq} \frac{c^2\exp(-2(1+\log2)\tau^2}{4(1+\log2)^2\gamma^2}\cdot\bigg((1+\log2-\frac{1}{1+\log2})\frac{p_1r}{3}-\log2\bigg)\stackrel{\text{(ii)}}{\geq} c\frac{p_1r\tau^2}{\gamma^2}, \end{align*}
    where (i) is because of the assumption $$p_1 \geq 2(r+r_{1*}+r_{2*})$$ and (ii) is because  $$(1+\log2-\frac{1}{1+\log2})\frac{p_1r}{3}-\log2 \geq \frac{p_1r}{c_1}$$ for large enough constant $c_1>0$.
    Considering $\gamma^2 \geq C\tau^2(p_{V1}+p_{V2})$,
    $$
    \frac{p_1r\tau^2}{\gamma^2} \geq c\frac{p_1r\tau^2(\gamma^2+\tau^2(p_{V1}+p_{V2}))}{\gamma^4}.$$
    As a result 
    $$
    \underset{\hat{U}_r}{\inf}\underset{U_r \in \mathbb{O}}{\sup}\mathbb{E}_{P_{U_r}}\|\sin\Theta(\hat{U}_r,U_r)\|_F^2 \geq c\bigg(\frac{p_1r\tau^2(\gamma^2+\tau^2(p_{V1}+p_{V2}))}{\gamma^4}
    \wedge r\bigg).
    $$
    By Lemma 6 in \cite{caizhang2018} and taking $\tau=1$ without loss of generality, we further have
    $$
    \underset{\hat{U}_r}{\inf}\underset{X \in \mathscr{H}_{r,\gamma}^{(u)}}{\sup}\mathbb{E}\|\sin\Theta(\hat{U}_r,U_r)\|_F^2 \geq c\bigg(\frac{p_1r(\gamma^2+(p_{V1}+p_{V2}))}{\gamma^4}
    \wedge r\bigg).
    $$
Next we consider the circumstance when $p_{V1}\asymp p_{V2}$ holds. In this case, when $\gamma^2 \geq \frac{\tau^2}{16}p_{V1}\wedge p_{V2}$, the proof can be finished as the same procedure in the previous case. We have
    $$
    \frac{p_1r\tau^2}{\gamma^2} \geq \frac{p_1r\tau^2(\gamma^2/2+\tau^2(p_{V1}\wedge p_{V2})/32)}{\gamma^4} \geq c\frac{p_1r\tau^2(\gamma^2+\tau^2(p_{V1}+p_{V2}))}{\gamma^4},
    $$
    considering $p_{V1}\asymp p_{V2}$  and $\gamma^2 \geq \frac{\tau^2}{16}p_{V1}\wedge p_{V2}$.

 Then we consider the case when $p_{V1}\asymp p_{V2}$ holds and $\gamma^2 < \frac{\tau^2}{16}p_{V1}\wedge p_{V2}$. We provide the proof when $p_{V1}\leq p_{V2}$, and the proof for $p_{V2}\leq p_{V1}$ is symmetric. Define the following class of density $P_Y$ : 
    $$
    \mathscr{P}_{U_r,\gamma}=\{P_Y:~Y=X+Z \in \mathbb{R}^{p_1\times p_2},~Z\sim N(0,\tau^2), X \in \mathscr{H}_{r,\gamma}^{(u)},~denoted~U_r~of~X~in~\mathscr{H}_{r,\gamma}^{(u)}~is~U_rO~here, O \in \mathbb{O}(r,r)\}.
    $$
     Construct measure
    \begin{align*}\Bar{P}_{U_r,\gamma}(Y)&=C_{U_r,\gamma}\int_{W \in \mathbb{R}^{(p_{V1}-r_*) \times r}:\sigma_{\min}(W) \geq \frac{1}{2\sqrt{2}}}\frac{1}{(2\pi\tau^2)^{p_1p_2/2}}\bigg(\frac{p_{V1}}{2\pi}\bigg)^{\frac{r(p_{V1}-r_{1*})}{2}}\\&
    \exp\bigg(-\frac{p_{V1}\|W\|_F^2}{2}\bigg)\exp\bigg(-\bigg\|Y-2\gamma \begin{pmatrix}
    U_r~U_{1*}\Sigma_{1*}~U_{2*}\Sigma_{2*}
    \end{pmatrix}\begin{pmatrix}
        W^T& 0_{r_{1*}}& W^T &0_d& 0_{r_{2*}}\\ 0 & I_{r_{1*}} & 0 &0& 0\\0&0&0&0& I_{r_{2*}}
    \end{pmatrix}\bigg\|_F^2/2\tau^2\bigg)
    dW,\end{align*}
    where $C_{U_r,\gamma}$ is a constant to make it a valid probability measure, $0_d={\bf 0}_{r\times (p_{V2}-p_{V1}-r_{2*})}$ and $0_{r_{i*}}$ denotes the zero block matrix with the dimension of ${\bf 0}_{r\times r_{i*}}$. Here $\Sigma_{1*}$ and $\Sigma_{2*}$  can be selected to satisfy the gap condition. (e.g., select them as $\gamma\sigma_{\min}(W)I$). We can verify that
    \begin{align*}2\gamma \begin{pmatrix}
    U_r~U_{1*}\Sigma_{1*}~U_{2*}\Sigma_{2*}
    \end{pmatrix}\begin{pmatrix}
        W^T& 0_{r_{1*}}& W^T &0_d& 0_{r_{2*}}\\ 0 & I_{r_{1*}} & 0 &0& 0\\0&0&0&0& I_{r_{2*}}
    \end{pmatrix}=\begin{pmatrix}
        U_rO&U_{1*}&U_{2*}
    \end{pmatrix}\mathrm{diag}(
        2\gamma \sqrt{\Sigma^2+\Sigma^2} ,\Sigma_{1*},\Sigma_{2*}) \\ \cdot\begin{pmatrix}
        \frac{2\gamma\Sigma}{\sqrt{4\gamma^2\Sigma^2+4\gamma^2\Sigma^2}}V^T&0_{r_{1*}}&\frac{2\gamma\Sigma}{\sqrt{4\gamma^2\Sigma^2+4\gamma^2\Sigma^2}}(V^T~0_d)&0_{r_{2*}}\\
        0&I&0&0\\0&0&0&I
    \end{pmatrix}.\end{align*}
    Under the event $\{\sigma_{\min}(W) \geq \frac{1}{2\sqrt{2}}\}$, the minimum singular value corresponding to $U_r$ is greater than $\gamma$, which means that $\Bar{P}_{U_r,\gamma}$ is indeed a mixture density in $\mathscr{P}_{U_r,\gamma}$.
    For bounding the KL-divergence, we have the following lemma.
    \begin{lemma}\label{unshareboundKL}
        Under the assumption of the theorem and $\gamma^2 \leq \frac{\tau^2p_{V1}}{16}$, for any $U_r,U_r' \in \mathbb{O}(p_1,r)$, we have $$
        D(\Bar{P}_{U_r,\gamma}||\Bar{P}_{U_r',\gamma}) \leq \frac{(p_{V1}-r_{1*})128\gamma^4}{(8\gamma^2+\tau^2p_{V1})(\tau^2p_{V1}-8\gamma^2)}\|\sin\Theta(U_r,U_r')\|_F^2+C_{KL},
        $$
        where $C_{KL}$ is a uniform constant.
    \end{lemma}
    When $\gamma^2 \leq \frac{\tau^2p_{V1}}{16}$, $$\tau^2p_{V1}-8\gamma^2 \geq \frac{\tau^2p_{V1}}{2}.$$ Further we have
    $$
    D(\Bar{P}_{U_r,\gamma}||\Bar{P}_{U_r',\gamma}) \leq \frac{256\gamma^4}{\tau^2(8\gamma^2+\tau^2p_{V1})}\|\sin\Theta(U_r,U_r')\|_F^2+C_{KL}.
    $$
    We still construct $U_r$ as $$U_r=(U_{0}^T~{\bf 0}_{r\times (r_{1*}+r_{2*})})^T,\qquad U_0\in \mathbb{O}(p_1-r_{1*}-r_{2*},r)$$
    Then we again construct a ball of radius $\epsilon$ centered at $U_{00} \in \mathbb{O}(p_1-r_{1*}-r_{2*},r)$:
    \[
    B(U_{00},\epsilon)=\{U_{0}':\|\sin\Theta(U_{00},U_0')\|_F^2 \leq \epsilon\}.
    \]
    Based on Lemma 1 in \cite{Cai2013}, for any $0<a<1$ and $0<\epsilon<\sqrt{2r}$, $\exists$ $\{U'_{01},...,U'_{0m}\} \subseteq B(U_{00},\epsilon)$ s.t.
    \[
    m\geq (\frac{c}{a})^{r(p_1-r_{1*}-r_{2*}-r)},\] and \[ 
    \underset{1\leq i\neq j \leq m}{\min}\|\sin\Theta(U'_{0i},U'_{0j})\|_F \geq a\epsilon.
    \]
    Also,
    \[\|\sin\Theta(U_{r1},U_{r2})\|_F=\|\sin\Theta((U_{01}^T~0_{r\times (r_{1*}+r_{2*})})^T,(U_{02}^T~0_{r\times (r_{1*}+r_{2*})})^T)\|_F=\|\sin\Theta(U_{01},U_{02})\|_F.\]
    Thus we can construct a ball of radius $\epsilon$ centered at $U_{r0} \in \mathbb{O}(p_1,r)$:
    \[
    B(U_{r0},\epsilon)=\{U_{r}':\|\sin\Theta(U_{r0},U_r')\|_F^2 \leq \epsilon\}.
    \]
    Also for any $0<a<1$ and $0<\epsilon<\sqrt{2r}$, $\exists$ $\{U'_{r1},...,U'_{rm}\} \subseteq B(U_{r0},\epsilon)$, where $U_{ri}=(U_{0i}^T~0_{r\times (r_{1*}+r_{2*})})^T$ s.t.
    \[
    m\geq (\frac{c}{a})^{r(p_1-r_{1*}-r_{2*}-r)}, \] and \[
    \underset{1\leq i\neq j \leq m}{\min}\|\sin\Theta(U'_{ri},U'_{rj})\|_F \geq a\epsilon.
    \]
    So 
    $$
        D(\Bar{P}_{U_r,\gamma}||\Bar{P}_{U_r',\gamma}) \leq \frac{256\gamma^4}{\tau^2(8\gamma^2+\tau^2p_{V1})}4\epsilon^2+C_{KL}.
        $$
    Define a class $\mathbb{O}$ as the class of $U_r$ constructed by the above means. Then, by Fano's Lemma \citep{Yu1997}, we have 
    $$\underset{\hat{U}_r}{\inf}\underset{U_r \in \mathbb{O}}{\sup}\mathbb{E}_{\hat{P}_{U_r,\gamma}}\|\sin\Theta(\hat{U}_r,U_r)\|_F^2 \geq a^2\epsilon^2\bigg(1-\frac{\frac{1024\gamma^4}{\tau^2(8\gamma^2+\tau^2p_2)}\epsilon^2+C_{KL}+\log(2)}{r(p_1-r_{1*}-r_{2*}-r)\log(c/a)}\bigg).$$
    Select $a=c\exp(-(1+C_{KL}+log2))$ and choose \[\epsilon=\sqrt{\frac{r(p_1-r)\tau^2(8\gamma^2+\tau^2p_2)}{2\cdot1024\gamma^4}} \wedge \sqrt{2r},\] when it holds that $$\sqrt{\frac{r(p_1-r_{1*}-r_{2*}-r)\tau^2(8\gamma^2+\tau^2p_2)}{2\cdot1024\gamma^4}} \leq \sqrt{2r},$$ we have 
    \begin{align*}
    \underset{\hat{U}_r}{\inf}\underset{U_r \in \mathbb{O}}{\sup}\mathbb{E}_{\hat{P}_{U_r,\gamma}}\|\sin\Theta(\hat{U}_r,U_r)\|_F^2
     \geq c^2\exp(-2(1+C_{KL}+\log2))
    \cdot\frac{r(p_1-r_{1*}-r_{2*}-r)(\tau^2(4\gamma^2+\tau^2p_2))}{2\cdot 1024\gamma^4}
    \\\cdot \bigg(r(p_1-r_{1*}-r_{2*}-r)(1+C_{KL}+\log2)-\frac{r(p_1-r_{1*}-r_{2*}-r)}{2}
    -(C_{KL}+\log2)\bigg)\\
    /(r(p_1-r_{1*}-r_{2*}-r)(1+C_{KL}+\log2))
    \geq c\frac{p_1r\tau^2(\gamma^2+\tau^2p_2)}{\gamma^4}
    .\end{align*}
    With Lemma 6 in \cite{caizhang2018}, 
    $$\underset{\hat{U}_r}{\inf}\underset{X \in \mathscr{H}_{r,\gamma}^{(u)}}{\sup}\mathbb{E}_{\hat{P}_{U_r,\gamma}}\|\sin\Theta(\hat{U}_r,U_r)\|_F^2\geq \underset{\hat{U}_r}{\inf}\underset{U_r \in \mathbb{O}}{\sup}\mathbb{E}_{\hat{P}_{U_r,\gamma}}\|\sin\Theta(\hat{U}_r,U_r)\|_F^2\geq c\bigg(\frac{p_1r\tau^2(\gamma^2+\tau^2p_{V1})}{\gamma^4} \wedge r\bigg).$$
    Noting $p_{V1}\asymp p_{V2}$ and taking $\tau=1$ without loss of generality,
    $$
    \underset{\hat{U}_r}{\inf}\underset{X \in \mathscr{H}_{r,\gamma}^{(u)}}{\sup}\mathbb{E}_{\hat{P}_{U_r,\gamma}}\|\sin\Theta(\hat{U}_r,U_r)\|_F^2
     \geq c\bigg(\frac{p_1r(\gamma^2+(p_{V1}+p_{V2}))}{\gamma^4} \wedge r\bigg).
    $$ We thus finish the proof.
\end{proof}

\subsection{Proof of technical lemmas}

\begin{proof}[PROOF OF LEMMA \ref{singularvaluerelation}]
$$
X=U(\Sigma_1,...,\Sigma_k)\mathrm{diag}(
    V_1^T,V_2^T,\\\cdots,V_k^T ),
$$
and 
\[
\sigma_i(X)=\sigma_i((\Sigma_1,...,\Sigma_k)),~ i=1,2,...,r
.\]
Suppose singular values of $X_k$ are $\sigma_{k1}, \sigma_{k2},...,\sigma_{kr}$ respectivel. Then 
$
\Sigma = 
(\Sigma_1,...,\Sigma_k).
$
We have 
$$
\Sigma\Sigma^T=\mathrm{diag}(\sum_{j=1}^k\sigma_{j(1)}^2,\sum_{j=1}^k\sigma_{j(2)}^2, ... , \sum_{j=1}^k\sigma_{j(r)}^2).$$
It is obvious that the $r$-th eigenvalue of $\Sigma\Sigma^T$ is the smallest one among the values $\sum_{j=1}^k\sigma_{j(i)}^2$ $(i=1,2,...,r)$. Then, the lemma follows.
\end{proof}

\begin{proof}[PROOF OF LEMMA \ref{shareboundKL}]
First consider the Gaussian
$$
    \Tilde{P}_{U,\gamma}(Y)=\int_{W \in \mathbb{R}^{\frac{p_2}{2} \times r}}\frac{1}{(2\pi\tau^2)^{p_1p_2/2}}\exp\bigg(-\frac{\|Y-2\gamma U(W^T~W^T~{\bf 0}_{r \times (p_{V2}-p_{V1})})\|_F^2}{2\tau^2}\bigg)
    \cdot \bigg(\frac{p_{V1}}{2\pi}\bigg)^{\frac{p_{V1}r}{2}}\exp\bigg(-\frac{p_{V1}\|W\|_F^2}{2}\bigg)dW.
$$
This is indeed a valid density. Split $Y$ column-wise, denoting the first $p_{V1}$ columns as $Y_1$, the $p_{V1}+1$ to $2p_{V1}$ columns as $Y_2$ and the last $p_{V2}-p_{V1}$ columns as $Y_3$. Then
\begin{align*}
    \Tilde{P}_{U,\gamma}=\frac{1}{(2\pi \tau^2)^{p_1p_2/2}}\bigg(\frac{p_{V1}}{2\pi}\bigg)^{p_{V1}r/2}\int_{W}\exp\bigg(-\frac{1}{2\tau^2}tr(Y_1^T(I-\frac{8\gamma^2}{8\gamma^2+\tau^2p_{V1}}UU^T)Y_1+Y_2^T(I-\frac{8\gamma^2}{8\gamma^2+\tau^2p_{V1}}UU^T)Y_2\\
     -Y_1^T\frac{8\gamma^2}{8\gamma^2+\tau^2p_{V1}}UU^TY_2-Y_2^T\frac{8\gamma^2}{8\gamma^2+\tau^2p_{V1}}UU^TY_1\\
     +\frac{(8\gamma^2+\tau^2p_2)}{2}(W-\frac{4\gamma}{8\gamma^2+\tau^2p_{V1}}(Y_1+Y_2)^TU)(W-\frac{4\gamma}{8\gamma^2+\tau^2p_{V1}}(Y_1+Y_2)^TU)^T\bigg)dW\\
     \cdot \exp\bigg(-\frac{1}{2\tau^2}tr(Y_3^TY_3)\bigg).
\end{align*}
Denote $i$-th column of $Y_1,Y_2$ as $Y_{1i},Y_{2i}$ respectively. Calculating the integral, 
\begin{align*}
    \Tilde{P}_{U,\gamma}=\frac{1}{(2\pi \tau^2)^{\frac{p_1p_2}{2}}}\bigg(\frac{p_{V1}}{2\pi}\bigg)^{\frac{p_{V1}r}{2}}\bigg(\frac{2\pi 2\tau^2}{8\gamma^2+\tau^2p_2}\bigg)^{\frac{p_{V1}r}{2}}\exp\bigg(-\frac{1}{2\tau^2}tr(Y_3^TY_3)\bigg)\exp\bigg(-\frac{1}{2\tau^2}\sum_{i=1}^{p_{V1}}\begin{pmatrix}
        Y_{1i}^T&Y_{2i}^T
    \end{pmatrix}\\\begin{pmatrix}
        I_{p_1}-\frac{8\gamma^2}{8\gamma^2+\tau^2p_{V1}}UU^T & -\frac{8\gamma^2}{8\gamma^2+\tau^2p_{V1}}UU^T\\ -\frac{8\gamma^2}{8\gamma^2+\tau^2p_{V1}}UU^T & I_{p_1}-\frac{8\gamma^2}{8\gamma^2+\tau^2p_{V1}}UU^T
    \end{pmatrix}\begin{pmatrix}
        Y_{1i} \\ Y_{2i}
    \end{pmatrix}\bigg).
\end{align*}
From the expression, $(Y_{1i}^T~Y_{2i}^T)$ are Gaussian and are independent of $Y_3$. For $i=1,2,...,p_{V1}$, $(Y_{1i}^T~Y_{2i}^T)$ iid follows $$
    N\bigg(0,\tau^2\begin{pmatrix}
        I_{p_1}-\frac{8\gamma^2}{8\gamma^2+\tau^2p_{V1}}UU^T & -\frac{8\gamma^2}{8\gamma^2+\tau^2p_{V1}}UU^T\\ -\frac{8\gamma^2}{8\gamma^2+\tau^2p_{V1}}UU^T & I_{p_1}-\frac{8\gamma^2}{8\gamma^2+\tau^2p_{V1}}UU^T
    \end{pmatrix}^{-1}\bigg)=N\bigg(0,\frac{1}{\tau^2}\begin{pmatrix}
        I_{p_1}+\frac{8\gamma^2}{\tau^2p_{V1}-8\gamma^2}UU^T & \frac{8\gamma^2}{\tau^2p_{V1}-8\gamma^2}UU^T\\ \frac{8\gamma^2}{\tau^2p_{V1}-8\gamma^2}UU^T & I_{p_1}+\frac{8\gamma^2}{\tau^2p_{V1}-8\gamma^2}UU^T 
    \end{pmatrix}\bigg).$$
Then by KL-divergence for Gaussian:
\begin{align*}
    D(\Tilde{P}_{U,\gamma}||\Tilde{P}_{U',\gamma})=\frac{p_{V1}}{2}\bigg\{tr(\begin{pmatrix}
        I_{p_1}-\frac{8\gamma^2}{8\gamma^2+\tau^2p_{V1}}UU^T & -\frac{8\gamma^2}{8\gamma^2+\tau^2p_{V1}}UU^T\\ -\frac{8\gamma^2}{8\gamma^2+\tau^2p_{V1}}UU^T & I_{p_1}-\frac{8\gamma^2}{8\gamma^2+\tau^2p_{V1}}UU^T
    \end{pmatrix}  \\
    \cdot \begin{pmatrix}
        I_{p_1}+\frac{8\gamma^2}{\tau^2p_{V1}-8\gamma^2}U'U'^T & \frac{8\gamma^2}{\tau^2p_{V1}-8\gamma^2}U'U'^T\\ \frac{8\gamma^2}{\tau^2p_{V1}-8\gamma^2}U'U'^T & I_{p_1}+\frac{8\gamma^2}{\tau^2p_{V1}-8\gamma^2}U'U'^T 
    \end{pmatrix})-2p_1\bigg\}\\
    =\frac{p_{V1}128\gamma^4}{(8\gamma^2+\tau^2p_{V1})(\tau^2p_{V1}-8\gamma^2)}\|\sin\Theta(U,U')\|_F^2.\end{align*}
Next, by integral, we have 
$$
    \frac{\Bar{P}_{U,\gamma}(Y)}{\Tilde{P}_{U,\gamma}(Y)}
    =C_{U,\gamma}\cdot \mathbb{P}\{\sigma_{\min}(W) \geq \frac{1}{2\sqrt{2}}\;|\;W \in \mathbb{R}^{p_{V1} \times r},W\sim N\bigg(\frac{4\gamma}{8\gamma^2+\tau^2p_{V1}}(Y_1+Y_2)^TU,\frac{2\tau^2}{8\gamma^2+\tau^2p_{V1}}\bigg)\}.
$$
For any $Y_1+Y_2 \in \mathbb{R}^{p_1 \times p_{V1}}$,$(Y_1+Y_2)^TU\in \mathbb{R}^{p_{V1} \times r}$, we can find $Q \in \mathbb{O}(p_{V1},p_{V1}-r)$ s.t. $Q^T(Y_1+Y_2)^TU=0$. So the entries of $Q^TW$ iid follow  $N(0,\frac{2\sigma^2}{8\gamma^2+\tau^2p_{V1}})$. Then by Corollary 35 in \cite{Vershynin2012}, $$\sigma_{\min}(W)=\sigma_r(W)\geq \sigma_r(Q^TW)\geq \frac{\sqrt{2}\tau}{\sqrt{8\gamma^2+\tau^2p_{V1}}}(\sqrt{p_{V1}-r}-\sqrt{r}-x) ,$$
with probability at least $1-2\exp(-\frac{x^2}{2})$.
Under $r\leq \frac{p_{V1}}{16} \wedge \frac{p_{V2}}{16}$, $\gamma^2 \leq \frac{\tau^2}{16}p_{V1}\wedge p_{V2}$, we can find a small enough constant $k$. Set $x=k\sqrt{p_{V1}}$, we have
$$
    \sigma_{\min}(W)\geq \frac{\sqrt{2}\sigma}{\sqrt{8\gamma^2+\tau^2p_{V1}}}(\sqrt{p_{V1}-r}-\sqrt{r}-k\sqrt{p_{V1}})\geq \frac{1}{2\sqrt{2}}
$$
holds with probability at least $1-2\exp(-cp_{V1})$. Then,  $\exists c>0$ s.t. 
$$\mathbb{P}\bigg\{\sigma_{\min}(W) \geq \frac{1}{2\sqrt{2}}\;|\;W \in \mathbb{R}^{p_{V1} \times r},W\sim N\bigg(\frac{4\gamma}{8\gamma^2+\tau^2p_{V1}}(Y_1+Y_2)^TU,\frac{2\tau^2}{8\gamma^2+\tau^2p_{V1}}\bigg)\bigg\}\geq 1-2\exp(-cp_{V1}).$$
By definition, $$
    C_{U,\gamma}^{-1}=\int_YP_{U,\gamma}(Y)dY=\mathbb{P}\bigg\{\sigma_{\min}(W) \geq \frac{1}{2\sqrt{2}} \; |\;W \in \mathbb{R}^{p_{V1}\times r}, W \stackrel{\text{iid}}{\sim}N\bigg(0,\frac{1}{p_{V1}}\bigg)\bigg\}.$$
So, $C_{U,\gamma} \geq 1$ and with the Corollary 35 in \cite{Vershynin2012}, for large constant $k$, 
\[
1 \leq C_{U_r,\gamma} \leq 1+k\exp(-cp_{V1}).
\]
Combining above, for a large constant $k$,  $$1-k\exp(-cp_{V1}) \leq \frac{\Bar{P}_{U,\gamma}(Y)}{\Tilde{P}_{U,\gamma}(Y)} \leq 1+k\exp(-cp_{V1}).$$
Finally,
$$D(\Bar{P}_{U,\gamma}||\Bar{P}_{U',\gamma})\leq\int_{Y}\Bar{P}_{U,\gamma}(Y)\bigg(\log\bigg(\frac{\Tilde{P}_{U,\gamma}(Y)}{\Tilde{P}_{U',\gamma}(Y)}\bigg)+2\log(1+k\exp(-cp_{V1}))\bigg)dY.$$
As $2\log(1+k\exp(-cp_{V1}))$ can be bounded, $$D(\Bar{P}_{U,\gamma}||\Bar{P}_{U',\gamma})\leq C+\frac{p_{V1}128\gamma^4}{(8\gamma^2+\tau^2p_{V1})(\tau^2p_{V1}-8\gamma^2)}\|\sin\Theta(U,U')\|_F^2+k\exp(-cp_{V1})\int_{Y}\bigg|\Tilde{P}_{U,\gamma}(Y)\log\bigg(\frac{\Tilde{P}_{U,\gamma}(Y)}{\Tilde{P}_{U',\gamma}(Y)}\bigg)\bigg|dY.$$
Then,
\begin{align*}
    &\int_{Y}\bigg|\Tilde{P}_{U,\gamma}(Y)\log\bigg(\frac{\Tilde{P}_{U,\gamma}(Y)}{\Tilde{P}_{U',\gamma}(Y)}\bigg)\bigg|dY=\int_{Y_1,Y_2}\Tilde{P}_{U,\gamma}(Y_1,Y_2)\bigg|\frac{8\gamma^2}{2\tau^2(8\gamma^2+\tau^2p_{V1})}\sum_{i=1}^{p_{V1}}\begin{pmatrix}
        Y_{1i}^T&Y_{2i}^T
    \end{pmatrix}\\
    &\qquad\qquad\begin{pmatrix}
        UU^T-U'U'^T & UU^T-U'U'^T\\UU^T-U'U'^T & UU^T-U'U'^T
    \end{pmatrix}\begin{pmatrix}
        Y_{1i}\\Y_{2i}
    \end{pmatrix}\bigg|dY_1dY_2\\
    &\leq \mathbb{E}\bigg\{\frac{8\gamma^2}{2\tau^2(8\gamma^2+\tau^2p_{V1})}\bigg(2\sum_{i=1}^{p_{V1}}Y_{1i}(UU^T+U'U'^T)Y_{1i}^T+2\sum_{i=1}^{p_{V1}}Y_{2i}(UU^T+U'U'^T)Y_{2i}^T\bigg) \mid ~(Y_{1i}^T~Y_{2i}^T)\\
    &\stackrel{\text{iid}}{\sim}N\bigg(0,\frac{1}{\tau^2}\begin{pmatrix}
        I_{p_1}+\frac{8\gamma^2}{\tau^2p_{V1}-8\gamma^2}UU^T & \frac{8\gamma^2}{\tau^2p_{V1}-8\gamma^2}UU^T\\ \frac{8\gamma^2}{\tau^2p_{V1}-8\gamma^2}UU^T & I_{p_1}+\frac{8\gamma^2}{\tau^2p_{V1}-8\gamma^2}UU^T 
    \end{pmatrix}\bigg)\bigg\}\leq \frac{32\gamma^2p_{V1}r}{\tau^2(\tau^2p_{V1}-8\gamma^2)} \leq \frac{32\gamma^2p_{V1}r}{\tau^2(\tau^2p_{V1})/2} \leq \frac{p_{V1}^2}{4\tau^2}.\end{align*}
Because $k\exp(-cp_{V1})\frac{p_{V1}^2}{4\tau^2}$ can be bounded by a uniform constant, we then complete the proof.
\end{proof}

\begin{proof}[PROOF OF LEMMA \ref{lemmageneralize}]
Suppose the SVD of X is $$X=U\Sigma V^T=U\Sigma (V_r ~ V_*)^T,$$ where $V_r\in \mathbb{O}(p_2,r)$, $V_*\in \mathbb{O}(p_2,r_*)$, where $\text{rank}(X)=r+r_*$, $U\in \mathbb  O(p_1,r+r_*)$ Extending V into full $p_2 \times p_2$ orthogonal matrix, we have $$(V ~ V_{\perp}) \in \mathbb{O}(p_2,p_2) \quad V_\perp \in \mathbb{O}(p_2,p_2-\text{rank}(X)).$$
Denote $YV_r:= Y_1$, then
$$\mathbb{E}Y^{\intercal}Y=X^{\intercal}X+p_{1}\tau^2I_{p_{2}}=V\Sigma^{2}V^{\intercal}+p_{1}\tau^2I_{p_{2}},$$
$$\mathbb{E}V_r^{\intercal}Y^{\intercal}YV_r=V_r^T\begin{pmatrix}
    V_r & V_*
\end{pmatrix}\Sigma^2 \begin{pmatrix}
    V_r & V_*
\end{pmatrix}^TV_r + p_1\tau^2 V_r^TV_r=\Sigma_r^2+p_1 \tau^2 I_r.$$
We introduce fixed normalization matrix $M$ as
$$M=\mathrm{diag}((\sigma_1^2(X)+p_1\tau^2)^{-\frac{1}{2}},\cdots,(\sigma_r^2(X)+p_1\tau^2)^{-\frac{1}{2}}).$$
This design yields
$$M^\intercal \mathbb{E}Y_1^\intercal Y_1M=I_r.$$ The first and last inequality can be derived following similar procedure as the proof of Lemma 4 in \cite{caizhang2018}. By best rank-$r$
approximation,
$$\sigma_{r}(Y)=\max_{\mathrm{rank}(B)\leq r}\|Y-B\|\geq\|Y-Y\cdot[V_r\:0]\|=\sigma_{\mathrm{max}}(YV_{r\perp}),$$ where $V_{r\perp}=(V_*~V_{\perp})$,
to switch our focus from $\sigma_r(Y)$ to $\sigma_\mathrm{max}(YV_{r\perp}).$ Next note that $$
\mathbb{E}V_{r\perp}^\intercal Y^\intercal YV_{r\perp}=V_{r\perp}^\intercal X^\intercal XV_{r\perp}+ \mathbb{E}V_{r\perp}^\intercal Z^\intercal ZV_{r\perp}=\mathrm{diag}(
    \Sigma^2_*,~0)+p_1\tau^2 I_{p_2-r}$$
and 
$$
\left\|\mathrm{diag}(
    \Sigma^2_* ,~ 0)\right\|\leq \sigma_{r+1}^2(X),
$$
We have
\[\sigma_{\max}^{2}(YV_{r\perp})=\sigma_{\max}(V_{r\perp}^{\intercal}Y^{\intercal}YV_{r\perp})\leq \sigma_{r+1}^2(X)+p_{1}\tau^2+\left\|V_{r\perp}^{\intercal}Y^{\intercal}YV_{r\perp}-EV_{r\perp}^{\intercal}Y^{\intercal}YV_{r\perp}\right\|.\]
For $\left\|p_{1}^{-1}\tau^{-2}(V_{r\perp}^TY^{\intercal}YV_{r\perp}-\mathbb{E}V_{r\perp}Y^{\intercal}YV_{r\perp})\right\|$, we have:
$$u^{\mathsf{T}}V_{r\perp}^TY^TYV_{r\perp}u-u^{\mathsf{T}}\mathbb{E}V_{r\perp}^TY^TYV_{r\perp}u=2u^TV_{r\perp}X^TZ^TV_{r\perp}u+u^TV_{r\perp}(Z^TZ-p_1\tau^2I)V_{r\perp}u.$$
For the second term, based on essentially the same procedure as the previous proof. of Lemma 4 in \cite{caizhang2018}, one can show that $$\mathbb{P}\left(\left\|p_{1}^{-1}\tau^{-2}u^TV_{r\perp}(Z^TZ-p_1\tau^2I)V_{r\perp}u\right\|\geq x\right)\leq C\exp\left\{-cp_{1}\left(x^{2}\wedge x\right)\right\}.$$
Also,
\[
p_1^{-1}\tau^{-2}u^TV_{r\perp}X^TZ^TV_{r\perp}u=\vec{\mathbf{z}}^{\mathsf{T}}\text{vec}(p_1^{-1}\tau^{-2}V_{r\perp}u(XV_{r\perp}u)^T).
\]
By assumption on $\sigma_{r+1}^2(X)$, we further have:
$$\left\|\text{vec}(p_1^{-1}\tau^{-2}V_{r\perp}u(XV_{r\perp}u)^T)\right\|_2^2=p_1^{-2}\tau^{-2}\left\|V_{r\perp}u\right\|_2^2\left\|XV_{r\perp}u\right\|_2^2\leq p_1^{-2}\tau^{-2}\left\|XV_{r\perp}u\right\|_2^2=p_1^{-2} \tau^{-2}\sigma_{r+1}^2(X)\leq c\frac{\sigma_r^2(X)}{p_1^2\tau^2}.$$
Then by basic property of iid sub-Gaussian random variables
$$\mathbb{P}\left(\left\|p_1^{-1}\tau^{-2}u^TV_{r\perp}X^TZ^TV_{r\perp}u\right\|\geq x\right)\leq C\exp\left(-cx^2\frac{p_1^2\tau^2}{\sigma_r^2(X)}\right).$$
Similarly, combining above,
\begin{align*}&\mathbb{P}\left(\left\|p_{1}^{-1}\tau^{-2}u^T(V_{r\perp}Y^{\intercal}YV_{r\perp}-EV_{r\perp}Y^{\intercal}YV_{r\perp})u\right\|\geq x\right) \\
&\leq \mathbb{P}\bigg(\{\|p_{1}^{-1}\tau^{-2}2u^TV_{r\perp}X^TZ^TV_{r\perp}u\| \geq \frac{x}{2}\} \cup \{\|p_{1}^{-1}\tau^{-2}u^TV_{r\perp}(Z^TZ-p_1\tau^2I)V_{r\perp}u\|\geq \frac{x}{2}\}\bigg) \\
&\leq 2\cdot\mathbb{P}\left(\{\|p_{1}^{-1}\tau^{-2}2u^TV_{r\perp}X^TZ^TV_{r\perp}u\| \geq \frac{x}{2}\}\right) \vee \mathbb{P}\left(\{\|p_{1}^{-1}\tau^{-2}u^TV_{r\perp}(Z^TZ-p_1\tau^2I)V_{r\perp}u\|\geq \frac{x}{2}\}\right) \\
&\leq C\exp\left(-c (p_1x^2) \wedge (p_1x) \wedge \bigg(x^2\frac{p_1^2\tau^2}{\sigma_r^2(X)}\bigg) \right).\end{align*}
Thus,
$$\mathbb{P}\left(\left\|p_{1}^{-1}\tau^{-2}(V_{r\perp}Y^{\intercal}YV_{r\perp}-EV_{r\perp}Y^{\intercal}YV_{r\perp})\right\|\geq x\right)\leq C\exp\left(Cp_{2}-c (p_1x^2) \wedge (p_1x) \wedge \bigg(x^2\frac{p_1^2\tau^2}{\sigma_r^2(X)}\bigg)\right).$$
We obtain the second target inequality by combining inequalities.
\end{proof}

\begin{proof}[PROOF OF LEMMA \ref{unshareboundKL}]
Consider the Gaussian
\begin{align*}
    \Bar{P}_{U_r,\gamma}(Y)=\int_{W \in \mathbb{R}^{(\frac{p_2}{2}-r_*) \times r}}\frac{1}{(2\pi\tau^2)^{p_1p_2/2}}\bigg(\frac{p_{V1}}{2\pi}\bigg)^{\frac{r(p_{V1}-r_{1*})}{2}}\exp\bigg(-\frac{p_{V1}\|W\|_F^2}{2}\bigg) \exp\bigg(-\bigg\|Y-2\gamma \begin{pmatrix}
    U_r~U_{1*}\Sigma_{1*}~U_{2*}\Sigma_{2*}
    \end{pmatrix}\\ \cdot\begin{pmatrix}
        W^T& 0_{r_{1*}}& W^T &0_d& 0_{r_{2*}}\\ 0 & I_{r_{1*}} & 0 &0& 0\\0&0&0&0& I_{r_{2*}}
    \end{pmatrix}\bigg\|_F^2/2\tau^2\bigg)
    dW.\end{align*}
This is valid density. Split $Y$ column-wise, denoting the first $p_{V1}-r_{1*}$ columns as $Y_1$, the $p_{V1}-r_{1*}+1$ to $p_{V1}$ columns as $Y_{1*}$, the $p_{V1}+1$ to $2p_{V1}-r_{1*}$ columns as $Y_{2}$ , the $2p_{V1}-r_{1*}+1$ to $p_2-r_{2*}$ as $Y_3$ and the last $r_{2*}$ columns as $Y_{2*}$. Then
\begin{align*}
    &\Tilde{P}_{U_r,\gamma}=C_{\gamma}\int_{W}\exp\bigg(-\frac{1}{2\tau^2}tr(Y_1^T(I-\frac{8\gamma^2}{8\gamma^2+\tau^2p_{V1}}U_rU_r^T)Y_1+Y_2^T(I-\frac{8\gamma^2}{8\gamma^2+\tau^2p_{V1}}U_rU_r^T)Y_2\\
    &\qquad\qquad\qquad-Y_1^T\frac{8\gamma^2}{8\gamma^2+\tau^2p_{V1}}U_rU_r^TY_2-Y_2^T\frac{8\gamma^2}{8\gamma^2+\tau^2p_{V1}}U_rU_r^TY_1\\
    &+\frac{(8\gamma^2+\tau^2p_{V1})}{2}\bigg(W-\frac{4\gamma}{8\gamma^2+\tau^2p_{V1}}(Y_1+Y_2)^TU_r\bigg)\bigg(W-\frac{4\gamma}{8\gamma^2+\tau^2p_{V1}}(Y_1+Y_2)^TU_r\bigg)^T\bigg)dWP(Y_{1*})P(Y_2*)P(Y_3),\end{align*}
where $P(Y_{1*}),P(Y_2*),P(Y_3)$ are the probability measure for $Y_{1*}$, $Y_{2*}$ and $Y_3$ respectively and $C_{\gamma}$ is a constant left after extracting $P(Y_{1*}),P(Y_2*),P(Y_3)$. Denote the probability measure for $Y_1,Y_2$ as $\Tilde{P}_{U_r,\gamma}(Y_1,Y_2)$.
Denote i-th column of $Y_1,Y_2$ as $Y_{1i},Y_{2i}$ respectively. By calculating integral, \begin{align*}\Tilde{P}_{U_r,\gamma}=P(Y_{1*})P(Y_2*)C_{\gamma}\bigg(\frac{2\pi 2\tau^2}{8\gamma^2+\tau^2p_{V1}}\bigg)^{\frac{(p_{V1}-r_{1*})r}{2}}\exp\bigg(-\frac{1}{2\tau^2}\sum_{i=1}^{p_{V1}-r_{1*}}\begin{pmatrix}
        Y_{1i}^T&Y_{2i}^T
    \end{pmatrix}\\ \begin{pmatrix}
        I_{p_1}-\frac{8\gamma^2}{8\gamma^2+\tau^2p_{V1}}U_rU_r^T & -\frac{8\gamma^2}{8\gamma^2+\tau^2p_{V1}}U_rU_r^T\\ -\frac{8\gamma^2}{8\gamma^2+\tau^2p_{V1}}U_rU_r^T & I_{p_1}-\frac{8\gamma^2}{8\gamma^2+\tau^2p_{V1}}U_rU_r^T
    \end{pmatrix}\begin{pmatrix}
        Y_{1i} \\ Y_{2i}
    \end{pmatrix}\bigg).\end{align*}
From the form of the expression, it is obvious that the stacked $(Y_{1i}^T~Y_{2i}^T)$ independently follow Gaussian distribution. For $i=1,2,...,p_{V1}-r_{1*}$ we have $(Y_{1i}^T~Y_{2i}^T)$ iid follows  \begin{align*}
    &N\bigg(0,\tau^2\begin{pmatrix}
        I_{p_1}-\frac{8\gamma^2}{8\gamma^2+\tau^2p_{V1}}U_rU_r^T & -\frac{8\gamma^2}{8\gamma^2+\tau^2p_{V1}}U_rU_r^T\\ -\frac{8\gamma^2}{8\gamma^2+\tau^2p_{V1}}U_rU_r^T & I_{p_1}-\frac{8\gamma^2}{8\gamma^2+\tau^2p_{V1}}U_rU_r^T
    \end{pmatrix}^{-1}\bigg)\\
    &=N\bigg(0,\frac{1}{\tau^2}\begin{pmatrix}
        I_{p_1}+\frac{8\gamma^2}{\tau^2p_{V1}-8\gamma^2}UU^T & \frac{8\gamma^2}{\tau^2p_{V1}-8\gamma^2}UU^T\\ \frac{8\gamma^2}{\tau^2p_{V1}-8\gamma^2}UU^T & I_{p_1}+\frac{8\gamma^2}{\tau^2p_{V1}-8\gamma^2}UU^T 
    \end{pmatrix}\bigg).\end{align*} By Gaussian KL
\begin{align*}&D(\Tilde{P}_{U_r,\gamma}(Y)||\Tilde{P}_{U_r',\gamma}(Y))=\int_{Y_1,Y_2}\Tilde{P}_{U_r,\gamma}(Y_1,Y_2)\log\frac{\Tilde{P}_{U_r,\gamma}(Y_1,Y_2)}{\Tilde{P}_{U_r',\gamma}(Y_1,Y_2)}dY_1dY_2\int_{Y_{1*},Y_{2*}}P(Y_{1*})\\&~P(Y_{2*})P(Y_3)dY_{1*}dY_{2*}dY_3=D(\Tilde{P}_{U_r,\gamma}(Y_1,Y_2)||\Tilde{P}_{U_r',\gamma}(Y_1,Y_2))=\frac{(p_{V1}-r_{1*})128\gamma^4}{(8\gamma^2+\tau^2p_{V1})(\tau^2p_{V1}-8\gamma^2)}\|\sin\Theta(U_r,U_r')\|_F^2.\end{align*}
Next from the simplified expression of $\Tilde{P}_{U_r,\gamma}$, we have 
\begin{align*}
    &\frac{\Bar{P}_{U_r,\gamma}(Y)}{\Tilde{P}_{U_r,\gamma}(Y)}=C_{U_r,\gamma}\bigg(\frac{8\gamma^2+\tau^2p_{V1}}{2\pi 2\tau^2}\bigg)^{\frac{(p_{V1}-r_{1*})r}{2}}\cdot \int_{\sigma_{\min}(W) \geq \frac{1}{2\sqrt{2}}}\\
    &\qquad\qquad\exp\bigg(\frac{8\gamma^2+\tau^2p_{V1}}{2\tau^2 \cdot 2}\bigg\|W-\frac{4\gamma}{8\gamma^2+\tau^2p_{V1}}(Y_1+Y_2)^TU_r\bigg\|_F^2\bigg)dW\\
    &=C_{U_r,\gamma}\mathbb{P}\bigg\{\sigma_{\min}(W) \geq \frac{1}{2\sqrt{2}}\;|\;W \in \mathbb{R}^{(p_{V1}-r_{1*}) \times r}, W\sim N\bigg(\frac{4\gamma}{8\gamma^2+\tau^2p_{V1}}(Y_1+Y_2)^TU_r,\frac{2\tau^2}{8\gamma^2+\tau^2p_{V1}}\bigg)\bigg\}.
\end{align*}
For any fixed $Y_1+Y_2 \in \mathbb{R}^{p_1 \times (p_{V1}-r_{1*})}$, $(Y_1+Y_2)^TU_r\in \mathbb{R}^{(p_{V1}-r_{1*}) \times r}$, we can find $Q \in \mathbb{O}(p_{V1}-r_{1*},p_{V1}-r-r_{1*})$, which is orthogonal to $(Y_1+Y_2)^TU_r$ such that $$Q^T(Y_1+Y_2)^TU_r=0.$$ So the entries of $Q^TW$ iid follow normal distribution $N(0,\frac{2\tau^2}{8\gamma^2+\tau^2p_{V1}})$. Then by Corollary 35 in \cite{Vershynin2012},
\begin{align*}
    \sigma_{\min}(W)=\sigma_r(W)\geq \sigma_r(Q^TW)=\frac{\sqrt{2}\tau}{\sqrt{8\gamma^2+\tau^2p_{V1}}}\sigma_r\bigg(\frac{\sqrt{8\gamma^2+\tau^2p_{V1}}}{\sqrt{2}\tau}Q^TW\bigg)\\\geq \frac{\sqrt{2}\tau}{\sqrt{8\gamma^2+\tau^2p_{V1}}}(\sqrt{p_{V1}-r-r_{1*}}-\sqrt{r}-x),\end{align*}
with probability at least $1-2\exp(-\frac{x^2}{2})$.
Under the assumption $p_{V1} \geq 16(r+r_{1*})$, $\gamma^2 \leq \frac{\tau^2p_{V1}}{16}$, we can find a small enough constant $k$. Set $x=k\sqrt{p_{V1}}$, we have $$
    \sigma_{\min}(W)\geq \frac{1}{\sqrt{p_{V1}}}\bigg(\sqrt{\frac{15p_{V1}}{16}}-\sqrt{\frac{p_{V1}}{16}}-k\sqrt{{p_{V1}}}\bigg)\geq \frac{1}{2\sqrt{2}}$$
holds with probability at least $1-2\exp(-cp_{V1})$. As a result, there exists a constant $c$ such that $$
\mathbb{P}\bigg\{\sigma_{\min}(W) \geq \frac{1}{2\sqrt{2}}\;|\;W \in \mathbb{R}^{(p_{V1}-r_{1*}) \times r}, W\sim N\bigg(\frac{4\gamma}{8\gamma^2+\tau^2p_{V1}}(Y_1+Y_2)^TU_r,\frac{2\tau^2}{8\gamma^2+\tau^2p_{V1}}\bigg)\bigg\}\geq 1-2\exp(-cp_{V1}).$$
Recall the definition of $C_{U_r,\gamma}$, we have $$
    C_{U_r,\gamma}^{-1}=\int_YP_{U_r,\gamma}(Y)dY=\mathbb{P}\bigg\{\sigma_{\min}(W) \geq \frac{1}{2\sqrt{2}} \; |\;W \in \mathbb{R}^{(p_{V1}-r_{1*})\times r}, W \stackrel{\text{iid}}{\sim}N\bigg(0,\frac{1}{p_{V1}}\bigg)\bigg\}.$$
So, $C_{U_r,\gamma} \geq 1$ and again with the Corollary 35 in \cite{Vershynin2012}, we have for large constant $k$, 
\[
1 \leq C_{U_r,\gamma} \leq 1+k\exp(-cp_{V1}).
\]
Combining above, for a large enough constant $k$, we have the following bound $$
    1-k\exp(-cp_{V1}) \leq \frac{\Bar{P}_{U_r,\gamma}(Y)}{\Tilde{P}_{U_r,\gamma}(Y)}=\frac{\Bar{P}_{U_r,\gamma}(Y_1,Y_2)}{\Tilde{P}_{U_r,\gamma}(Y_1,Y_2)} \leq 1+k\exp(-cp_{V1}).$$
Finally we are able to bound the target KL-divergence: \begin{align*}
    &D(\Bar{P}_{U_r,\gamma}(Y)||\Bar{P}_{U_r',\gamma}(Y))=D(\Bar{P}_{U_r,\gamma}(Y_1,Y_2)||\Bar{P}_{U_r',\gamma}(Y_1,Y_2))\\&\leq\int_{Y_1,Y_2}\Bar{P}_{U_r,\gamma}(Y_1,Y_2)\bigg(\log\bigg(\frac{\Tilde{P}_{U_r,\gamma}(Y_1,Y_2)}{\Tilde{P}_{U_r',\gamma}(Y_1,Y_2)}\bigg)+2\log\bigg(1+k\exp(-cp_{V1})\bigg)\bigg)dY_1dY_2.\end{align*}
As $2\log(1+k\exp(-cp_{V1}))$ can be bounded with a uniform constant $C$:
\begin{align*}
    &D(\Bar{P}_{U_r,\gamma}||\Bar{P}_{U_r',\gamma})
    \leq C+\int_{Y_1,Y_2}\Tilde{P}_{U_r,\gamma}(Y_1,Y_2)\log\bigg(\frac{\Tilde{P}_{U_r,\gamma}(Y_1,Y_2)}{\Tilde{P}_{U_r',\gamma}(Y_1,Y_2)}\bigg)dY_1dY_2\\
    &\qquad\qquad+ \int_{Y_1,Y_2}|\Tilde{P}_{U_r,\gamma}(Y_1,Y_2)-\Bar{P}_{U_r,\gamma}(Y_1,Y_2)|\log\bigg(\frac{\Tilde{P}_{U_r,\gamma}(Y_1,Y_2)}{\Tilde{P}_{U_r',\gamma}(Y_1,Y_2)}\bigg)dY_1dY_2
    \\&\leq C+\frac{(p_{V1}-r_{1*})128\gamma^4}{(8\gamma^2+\tau^2p_{V1})(\tau^2p_{V1}-8\gamma^2)}\|\sin\Theta(U_r,U_r')\|_F^2\\
    &\qquad\qquad+k\exp(-cp_{V1}) \int_{Y_1,Y_2}\bigg|\Tilde{P}_{U_r,\gamma}(Y_1,Y_2)\log\bigg(\frac{\Tilde{P}_{U_r,\gamma}(Y_1,Y_2)}{\Tilde{P}_{U_r',\gamma}(Y_1,Y_2)}\bigg)\bigg|dY_1dY_2.\end{align*}
We then focus on the last term, similar to the analysis in the proof of Lemma \ref{shareboundKL},
$$
    \int_{Y_1,Y_2}\bigg|\Tilde{P}_{U_r,\gamma}(Y_1,Y_2)\log\bigg(\frac{\Tilde{P}_{U_r,\gamma}(Y_1,Y_2)}{\Tilde{P}_{U_r',\gamma}(Y_1,Y_2)}\bigg)\bigg|dY_1dY_2\leq \frac{p_{V1}^2}{4\tau^2}.$$
Considering $k\exp(-cp_2)\frac{p_{V1}^2}{4\tau^2}$ goes to $0$ when $p_{V1}$ go to $\infty$, this value can be bounded by a uniform constant. This completes our proof of Lemma \ref{unshareboundKL}.
\end{proof}


%

\ifCLASSOPTIONcaptionsoff
  \newpage
\fi



\bibliographystyle{IEEEtran}
\bibliography{Transactions-Bibliography/bibliography}
%



%

\begin{IEEEbiographynophoto}{Zhengchi Ma}
received the B.S. degree in
Data Science and Big Data Technologies from
Nankai University, Tianjin, China, in 2025. He is currently pursuing
the Ph.D. degree in the Department of Electrical \& Computer Engineering at Duke University. His research
interests include statistical inference for random matrices, high-dimensional statistics, and synthetic sampling  in statistics. He received the National
Scholarship in 2022 and 2024 and was recognized as an Excellent Graduate of Nankai University in 2025.

\end{IEEEbiographynophoto}

\begin{IEEEbiographynophoto}{Rong Ma} is an Assistant Professor of biostatistics at Harvard T.H. Chan School of Public Health, and an Associate Member at the Eric and Wendy Schmidt Center at the Broad Institute of MIT and Harvard. He received his Ph.D. in biostatistics from the University of Pennsylvania and was a postdoctoral scholar in statistics at Stanford University. His current research focuses on statistical inference for large random matrices, manifold learning, geometric inference, and data integration methods for biomedical applications, such as single-cell genomics and multiomics. He was a recipient of the 2022 Lawrence D. Brown Ph.D. Student Award from the Institute of Mathematical Statistics.
\end{IEEEbiographynophoto}





\end{document}